\documentclass[11pt]{article}
\usepackage[utf8]{inputenc}
\usepackage[T1]{fontenc}
\usepackage[british]{babel}
\usepackage{amsmath,amssymb,amsthm}
\usepackage[margin=2.6cm]{geometry}
\usepackage[dvipsnames]{xcolor}
\usepackage{enumitem,booktabs,array}
\usepackage{tikz}
\usetikzlibrary{arrows.meta,positioning,calc}
\usepackage{hyperref}
\hypersetup{colorlinks=true,linkcolor=NavyBlue,citecolor=NavyBlue,urlcolor=NavyBlue,pdftitle={A Cobham theorem for expansive bases in orders of number fields},pdfauthor={Francesco L. Chiera, Taoufik Safer}}
\newtheorem{theorem}{Theorem}[section]
\newtheorem{lemma}[theorem]{Lemma}
\newtheorem{proposition}[theorem]{Proposition}
\newtheorem{corollary}[theorem]{Corollary}
\theoremstyle{definition}
\newtheorem{definition}[theorem]{Definition}
\theoremstyle{remark}
\newtheorem{remark}[theorem]{Remark}
\newtheorem{example}[theorem]{Example}
\newtheorem{question}[theorem]{Question}
\newtheorem*{theorem*}{Theorem}
\newcommand{\Z}{\mathbb Z}\newcommand{\Q}{\mathbb Q}\newcommand{\R}{\mathbb R}
\newcommand{\C}{\mathbb C}\newcommand{\N}{\mathbb N}\newcommand{\F}{\mathbb F}
\newcommand{\cO}{\mathcal O}\newcommand{\cK}{\mathcal K}\newcommand{\cG}{\mathcal G}
\newcommand{\cA}{\mathcal A}\newcommand{\cB}{\mathcal B}\newcommand{\cH}{\mathcal H}
\newcommand{\cI}{\mathcal I}\newcommand{\cW}{\mathcal W}\newcommand{\cX}{\mathcal X}
\newcommand{\cC}{\mathcal C}\newcommand{\cR}{\mathcal R}\newcommand{\bA}{\mathbb A}
\newcommand{\dig}{\mathrm{dig}}
\DeclareMathOperator{\supp}{supp}
\newcommand{\norm}[1]{\lVert #1\rVert}
\newcommand{\Per}{\mathrm{Per}}
\newcommand{\Perf}{\mathrm{Per}^{\infty}}
\newcommand{\Sinf}{\Sigma_\infty}
\newcommand{\val}[1]{[#1]}
\newcommand{\pr}{\mathrm{pr}}

\title{A Cobham theorem for expansive bases\\ in orders of number fields}
\author{Francesco L. Chiera \and Taoufik Safer}
\begin{document}
\maketitle
\begin{quote}\small
\noindent\textbf{Abstract.} Let $K$ be a number field of degree $d\ge2$, $R$ its ring of integers and $\alpha,\beta\in R$ two expansive bases (of absolute value $>1$ at every archimedean place). If $\Q(\alpha^r,\beta^s)=K$ and $\alpha^rR\ne\beta^sR$ for all $r,s\ge1$, a function from $R$ to a finite alphabet is automatic in both bases if and only if it is eventually periodic, that is, equal outside a finite set to a function periodic under an additive subgroup of $R$ of rank $d$. Both conditions are necessary. The proof distinguishes between bases with distinct and with equal prime supports; it uses the recognisability of constant-shape substitutions and comparisons of exact affine symmetries. The result extends to orders, to full-rank lattices in $K$ stable under multiplication by the bases and to commuting expansive matrices in $M_d(\Z)$ that generate a field of degree $d$ over $\Q$, under the corresponding arithmetic conditions.
\end{quote}



\section{Introduction}\label{sec:intro}

We use $\N$ for the set of non-negative integers. Let $k\ge2$ be an integer and let $\Omega$ be a finite alphabet. A sequence $f:\N\to\Omega$ is \emph{$k$-automatic} if a finite automaton, reading the base-$k$ digits of $n\in\N$, computes $f(n)$. It is \emph{eventually periodic} if there is an integer $p\ge1$ such that $f(n+p)=f(n)$ for all sufficiently large $n$. An eventually periodic sequence is automatic in every base; the characteristic function of the powers of $2$, on the other hand, is $2$-automatic but not $3$-automatic.

\begin{theorem*}[Cobham, 1969 \cite{Cobham}]
Let $k,l\ge2$ be multiplicatively independent integers, that is, $k^m\ne l^n$ for all positive integers $m,n$. For a finite alphabet $\Omega$, a sequence $f:\N\to\Omega$ is both $k$- and $l$-automatic if and only if it is eventually periodic.
\end{theorem*}

The hypothesis is arithmetic and cannot be dropped. The bases may share prime factors, as do $6$ and $10$; what matters is that no positive power of one is a positive power of the other. We ask what becomes of this principle when $\N$ is replaced by the ring of integers $R$ of a number field $K$ and the bases are two elements $\alpha,\beta\in R$.

\subsection{From Gaussian integers to number fields}
The Gaussian integers provide a natural first test case. For a base of the form $-p+i$, with $p$ a positive integer, every element of $\Z[i]$ has an expansion with digits in $\{0,\ldots,p^2\}$. Radix representations of algebraic integers were studied in \cite{Saf98} (valid number systems, normalisation by finite automata and decidability of validity) and \cite{Saf99} (polygonal number systems). Hansel and Safer \cite{HS} considered two such bases, $-p+i$ and $-q+i$, with $p\ne q$. Assuming the four exponentials conjecture, they showed that an infinite set automatic in both bases is syndetic, that is, its translates by a suitable finite set cover $\Z[i]$. They conjectured that it is eventually periodic: it agrees outside a finite set with a set periodic under an additive subgroup of $\Z[i]$ of rank $2$.

Bosma, Fokkink and Krebs \cite{BFK} showed that the question is not vacuous: for multiplicatively independent bases $a,b\in\Z[i]$ with $|a|,|b|\ge\sqrt5$, there is an $a$-automatic set which is not $b$-automatic. Bustos-Gajardo, Fokkink and Yassawi \cite{BFY} subsequently proved the conjecture unconditionally and in a wider form, following Krebs's proof of Cobham's theorem \cite{Krebs}, which avoids establishing syndeticity. For two multiplicatively independent Gaussian bases, at least one of which is not a root of a rational integer, every set automatic in both bases is eventually periodic. When both bases are roots of rational integers, there are counterexamples.

The approach of \cite{BFY} was extended in two unpublished preprints of the authors \cite{IQO,CMW} to imaginary quadratic orders and to Weil bases in orders of CM fields, that is, bases each of which has the same absolute value at every archimedean place. For such bases, suitable powers of the two bases have comparable absolute values at every archimedean place at once. This need no longer be possible in a general number field, for instance for $3$ and $3+\sqrt2$ in $\Z[\sqrt2]$, and a different approach is needed.

There is also an arithmetic difficulty: powers of the bases may generate the same ideal without being equal, or may both lie in a proper subfield. Multiplicative independence alone is therefore no longer sufficient.

Here we go beyond these restrictions on the field and on the archimedean sizes of the bases. We work in a number field of arbitrary signature and degree, with bases which need only be \emph{expansive}, under two arithmetic conditions on the pair. The proof uses the ring of integers and its prime ideals; a final reduction extends the result to orders and to lattices stable under the bases.

\subsection{The main theorem}
Let $K$ be a number field of degree $d=[K:\Q]$, and let $R=\cO_K$ be its ring of integers. An element $\gamma\in R$ is expansive if $|\sigma\gamma|>1$ at every archimedean place $\sigma$ of $K$. Every element of $R$ then has a finite base-$\gamma$ expansion in a suitable finite set of digits.

For a finite alphabet $\Omega$, a function $f:R\to\Omega$ is \emph{$\gamma$-automatic} if, for some such digit set, a finite automaton with output computes $f(x)$ from every expansion of $x\in R$, read least significant digit first. Section~\ref{sec:numeration} gives an equivalent characterisation in terms of the finiteness of kernels and shows that automaticity does not depend on the choice of digits. A subset of $R$ is automatic if its characteristic function is automatic; this is a special case of the notion of automaticity for sequences indexed by a group with an expansive endomorphism studied by von Haeseler \cite{vonHaeseler}.

The function $f$ is \emph{eventually periodic} if it agrees outside a finite set with a function periodic under an additive subgroup of $R$ of rank $d$.

The following two conditions concern the pair of bases, and are to hold for all positive integers $r,s$:
\[
 \text{(H)}\quad\Q(\alpha^r,\beta^s)=K;\qquad
 \text{(Ind)}\quad\alpha^rR\ne\beta^sR.
\]
Condition (Ind) implies multiplicative independence of the elements, but in general is stronger because of the presence of units.

\begin{theorem*}[Theorem~\ref{thm:main}]
Let $K$ be a number field of degree $d\ge2$, let $R=\cO_K$, and let $\alpha,\beta\in R$ be expansive elements satisfying \textup{(H)} and \textup{(Ind)}. For a finite alphabet $\Omega$, a function $f:R\to\Omega$ is automatic in both bases if and only if it is eventually periodic.
\end{theorem*}

The two conditions rule out distinct obstructions.

\emph{A common subfield.} In $R=\Z[\sqrt2]$, the bases $2$ and $3$ are expansive and multiplicatively independent, but the set $\Z$ is automatic in both bases without being eventually periodic: it is an additive subgroup of rank $1$, not $2$. More generally, if a proper subfield $F$ of $K$ contains a positive power of each base, then $R\cap F$ is automatic in both bases (Proposition~\ref{prop:subfield}). Condition (H) concerns the pair rather than each base separately: in $K=\Q(\sqrt2,\sqrt3)$, the powers of $22(1+\sqrt2)$ generate only $\Q(\sqrt2)$ and those of $2(4+\sqrt3)$ only $\Q(\sqrt3)$, but the two bases together satisfy (H).

\emph{A common filtration by ideals.} In $R=\Z[\sqrt2]$, the bases $3$ and $3(1+\sqrt2)$ satisfy (H) and are multiplicatively independent as elements. Since $1+\sqrt2$ is a unit, both bases generate the same ideal $3R$. An automaton reading a non-zero element $x\in R$ in either base can determine the parity of the largest integer $n\ge0$ such that $x\in3^nR$. The set of non-zero elements for which this integer is even is thus automatic in both bases, but is not eventually periodic (\S\ref{ssec:necessity}).

Neither condition implies the other. Each is necessary in the following sense: if either fails, there is a subset of $R$ which is automatic in both bases but not eventually periodic. The theorem therefore characterises, in degree at least two, the pairs of expansive bases for which Cobham's conclusion holds (Corollary~\ref{cor:characterisation}). An equivalent formulation of (H) is $K_\alpha K_\beta=K$, where, for $\gamma\in K^\times$, $K_\gamma$ is the smallest of the fields $\Q(\gamma^n)$ for positive integers $n$ (Lemma~\ref{lem:fieldpowers}).

The theorem extends to every order in $K$ containing the bases and, more generally, to every full-rank lattice $M\subset K$ stable under multiplication by them (Theorem~\ref{thm:orders}). It also has a matrix form for two commuting expansive matrices in $M_d(\Z)$ which generate a field of degree $d$ over $\Q$, under the corresponding forms of (H) and (Ind) (Corollary~\ref{cor:matrices}).

For imaginary quadratic fields, (Ind) is equivalent to multiplicative independence, and (H) says that the bases are not both roots of rational integers. For Weil bases, multiplicative independence implies (Ind) (Remark~\ref{rem:weil}). Our result thus recovers the theorem of \cite{BFY}, and covers in particular totally real fields and non-Weil bases in CM fields.

\subsection{Idea of the proof}
The proof rests on a common principle: under (H), a single non-zero translation leaving a function automatic in both bases unchanged outside a finite set forces eventual periodicity in degree at least two. We refer to this as the principle of periods at infinity. Its proof combines the finiteness of kernels at mixed scales, a combinatorial argument on their graph and the arithmetic condition (H).

The rest of the proof distinguishes two cases, according as the sets of prime ideals containing the bases (their \emph{supports}) differ or coincide.

\emph{Different supports.} We use symbolic dynamics, in particular the recognisability of primitive aperiodic constant-shape substitutions, to show that all minimal subsystems of the orbit closure of the function under translations are finite orbits. Suppose that a prime ideal contains $\alpha$ but not $\beta$. An aperiodic minimal subsystem would come from a substitution subshift mapping equivariantly onto all the groups $R/\nu^nR$, for a power $\nu$ of $\alpha$; automaticity in base $\beta$ bounds the image of $\beta^kR$ in these groups, whereas $\beta$ is invertible modulo the powers of that prime ideal. Taking differences reduces the problem to the study of sets whose closures in the $\alpha$-adic completion have empty interior. Using local logarithms and a theorem on fibres in this completion, we show that certain affine symmetries of kernels have a common centre. This rules out branching within certain components of the graph of kernels and leads to sparse sets. The theorem of Evertse, Schlickewei and Schmidt on unit equations \cite{ESS} rules out infinite sets of this kind automatic in both bases.

\emph{Equal supports.} Here every power of $\alpha$ divides a power of $\beta$, and the bound on the image of $\beta^kR$ used in the first case gives nothing. We study instead quotients of powers of the bases which belong to $R$, in particular the two extremal ones, the \emph{ray elements}, which have smaller prime support. When one of them is expansive, a transfer of automaticity brings us back to the first case. Otherwise each contracts at some archimedean place, and a comparison of two affine symmetries with opposite behaviour at such a place yields a common centre. After recentring, we establish periodicity on certain residue classes and deduce eventual periodicity.

Thus the comparison of centres of affine symmetries is a common ingredient of the two proofs, although it takes place in different settings: at a prime ideal dividing one base but not the other in the first case and at an archimedean place in the second (Remark~\ref{rem:twomethods}).

\emph{What we cite and what we prove.} Besides standard textbook results and the theorem on unit equations, we cite the kernel criterion \cite{Eilenberg,AS} and its Gaussian analogues \cite{BFY}, and Cobham's theorem (in degree one only). We prove the principle of periods at infinity in arbitrary signature, replacing the pumping lemma of \cite{BFY}, and establish the minimality and the first step of recognisability of primitive constant-shape substitutions, in the case needed here. The recognisability property is stated in \cite{Cab,CabezasLeroy}. The two cases of the main theorem are proved here.

\begin{example}[A running example]
The distinction between the two cases is illustrated in $R=\Z[\sqrt2]$ by the pairs $(7,3+\sqrt2)$ and $(14,14u)$, where $u=8+5\sqrt2$. Both pairs consist of expansive bases and satisfy (H) and (Ind). The former has different supports; the latter has equal supports, and the two elements $u$ and its conjugate $\bar u=8-5\sqrt2$ arising in the second case each contract at one real place. We return to this example in \S\ref{ssec:cone}.
\end{example}

\subsection{Degree one: back to Cobham's theorem}
For a function $f:R\to\Omega$, a non-zero $v\in R$ is a \emph{period at infinity} if $f(x+v)=f(x)$ for all but finitely many $x\in R$. Together with $0$, these periods form an additive subgroup $\Perf(f)$ of $R$.

The theorem requires $d\ge2$. On $\Z$, the function $\mathbf 1_{\Z_{\ge0}}$ is automatic in every expansive base and its group of periods at infinity is $\Z$, yet the function is not eventually periodic: the two half-lines carry different patterns. In degree at least two, the complement of a large cube in the lattice is connected (in the nearest-neighbour graph), and a full-rank group of periods at infinity gives a single periodic function outside a finite set (Proposition~\ref{prop:rank2}). When formulated in terms of periods at infinity, the conclusion holds in every degree.

\begin{theorem*}[Corollary~\ref{cor:farperiods}]
Let $K$ be a number field of degree $d\ge1$, let $R=\cO_K$, and let $\alpha,\beta\in R$ be expansive elements satisfying (H) and (Ind). For a finite alphabet $\Omega$, a function $f:R\to\Omega$ is automatic in both bases if and only if $\Perf(f)$ has rank $d$.
\end{theorem*}
For $d=1$, after replacing $\alpha,\beta$ by $\alpha^2,\beta^2$ (Lemma~\ref{lem:powers}), this is Cobham's theorem applied to the sequences $n\mapsto f(n)$ and $n\mapsto f(-n)$ on $\N$.

\medskip\noindent\emph{Organisation of the paper.} Sections~\ref{sec:numeration} and~\ref{sec:onefar} deal with numeration, kernels and the principle of periods at infinity. Sections~\ref{sec:orbits}--\ref{sec:different} treat different supports, and Sections~\ref{sec:transfer}--\ref{sec:equal} equal supports. Section~\ref{sec:main} assembles the results, proves the necessity of the hypotheses and gives criteria for determining which case applies. Section~\ref{sec:orders} treats orders, lattices and matrix bases. Section~\ref{sec:questions} proves two further results, on periods at infinity without (H) and on morphisms between substitution subshifts, and discusses modules over $\cO_K$, several bases, groups with expansive endomorphisms and decidability.

\section{Numeration and kernels}\label{sec:numeration}

Throughout this section, let $K$ be a number field of degree $d=[K:\Q]$, let $R=\cO_K$ be its ring of integers, and let $\Omega$ be a finite alphabet.

A \emph{kernel} of a function $f:R\to\Omega$ in an expansive base $\gamma\in R$ is a function $y\mapsto f(c+\gamma^ny)$ on $R$, where $n\in\N$ and $c$ ranges over a finite set meeting every residue class modulo $\gamma^nR$, as specified in \S\ref{ssec:automata}. Automaticity is then equivalent to finiteness of the set of kernels, and it is almost always in this form that we shall use it. Only expansiveness of the bases is needed in this section.

\subsection{Places and sizes}\label{ssec:places}

Let $K_{\R}:=K\otimes_{\Q}\R$ be the Minkowski space of $K$, and let $\Sinf$ be the set of archimedean places of $K$. We identify $R$ with its image under the canonical embedding into $K_{\R}$. For $\sigma\in\Sinf$ and $x\in K$, we write $|\sigma x|$ for the common absolute value of the images of $x$ under the embeddings $K\hookrightarrow\C$ inducing $\sigma$. We denote the local degree by $d_\sigma\in\{1,2\}$, so that $\sum_{\sigma\in\Sinf}d_\sigma=d$. Set $\norm x:=\max_{\sigma\in\Sinf}|\sigma x|$. We say that an estimate holds \emph{place by place} if it holds at every archimedean place; most of the estimates below are of this kind.

\begin{lemma}[Discreteness of $R$]\label{lem:proper}
The function $x\mapsto\norm x$ extends to a norm on $K_{\R}$, and $|\sigma(xy)|=|\sigma x||\sigma y|$ for all $x,y\in K$ and $\sigma\in\Sinf$. For every real number $B$, the set $\{x\in R:\norm x\le B\}$ is finite. There is a constant $c_R>0$, the \emph{covering radius}, such that every point of $K_{\R}$ is at distance at most $c_R$ from a point of $R$. Finally, for every $x\in R$, $|N_{K/\Q}(x)|=\prod_\sigma|\sigma x|^{d_\sigma}\in\Z$.
\end{lemma}
\begin{proof}
Choose one embedding for each archimedean place. The map $x\mapsto(\sigma x)_\sigma$ induces an identification of $K_{\R}$ with $\R^{r_1}\times\C^{r_2}$, where $r_1$ and $r_2$ are the numbers of real and complex places respectively. We equip this space with the supremum norm. The image of $R$ is a full-rank lattice (see \cite[Ch.~I, \S\S2, 5]{Neukirch}). The expression for the field norm follows by taking the product over all embeddings; the norm of an algebraic integer is an integer.
\end{proof}

We shall also use the identity $[R:\eta R]=|N_{K/\Q}(\eta)|$ for every $\eta\in R\setminus\{0\}$. Indeed, multiplication by $\eta$ is an injective endomorphism of the lattice $R$, and its determinant on the $\Q$-vector space $K$ is $N_{K/\Q}(\eta)$.

\begin{definition}\label{def:expansive}
An element $\gamma\in R$ is \emph{expansive} if $|\sigma\gamma|>1$ for every $\sigma\in\Sinf$. We set $\rho_\gamma:=\max_{\sigma\in\Sinf}|\sigma\gamma|^{-1}<1$.
\end{definition}

Division by $\gamma$ then contracts at every archimedean place. Moreover, $\gamma$ is not a unit, since $|N_{K/\Q}(\gamma)|>1$. For a non-zero prime ideal $P$ of $R$ and $x\in R\setminus\{0\}$, we write $v_P(x)$ for the exponent of $P$ in the prime ideal factorisation of $xR$, and set $v_P(0):=\infty$.

\subsection{Digit sets}\label{ssec:digits}

Throughout \S\S\ref{ssec:digits}--\ref{ssec:carries}, let $\gamma\in R$ be expansive. Given a \emph{digit set}, that is, a finite set $D\subseteq R$ containing $0$, the \emph{value} of a word $w=w_0\cdots w_{n-1}\in D^n$, where $n\in\N$, is $\val w_\gamma=\sum_{0\le j<n}\gamma^jw_j$ (least significant digit first). We call $w$ an \emph{expansion} of this value. Let $D^{(n)}$ be the set of values of words of length $n$, with $D^{(0)}=\{0\}$. For all $n,m\in\N$,
\begin{equation}\label{eq:nesting}
D^{(n+m)}=D^{(n)}+\gamma^nD^{(m)}.
\end{equation}
We call $D$ \emph{admissible} if every $z\in R$ has an expansion over $D$. The following lemma, which recalls the existence of valid numeration systems in the sense of \cite{Saf98}, provides such sets.

\begin{lemma}[Short expansions]\label{lem:short}
Let $r\ge2c_R$ be a real number and set $D_\gamma(r):=\{t\in R:|\sigma t|\le r|\sigma\gamma|\ \text{for all }\sigma\in\Sinf\}$. Let $n\ge1$ be an integer. Every $z\in R$ such that $|\sigma z|\le r|\sigma\gamma|^n$ at each archimedean place has an expansion of length at most $n$ over $D_\gamma(r)$. In particular, $D_\gamma(r)$ is admissible.
\end{lemma}
\begin{proof}
We argue by induction on $n$, as in \cite[Lemma~11]{BFY} (the case of $\Z[i]$, with $D=B(0,r|\gamma|)$). The case $n=1$ follows from the definition of $D_\gamma(r)$.

Suppose that $|\sigma z|\le r|\sigma\gamma|^{n+1}$ at every place. In $K_\R\simeq\R^{r_1}\times\C^{r_2}$ (proof of Lemma~\ref{lem:proper}), move each coordinate of $z/\gamma$ towards zero by a distance $c_R$, or to zero if its modulus is less than $c_R$. The resulting point has coordinate moduli at most $r|\sigma\gamma|^n-c_R$. By the covering radius property, we can choose $q\in R$ within distance $c_R$ of this point. Thus $|\sigma(q-z/\gamma)|\le2c_R$ and $|\sigma q|\le r|\sigma\gamma|^n$ at each place. It follows that $z-\gamma q\in D_\gamma(r)$. Applying the induction hypothesis to $q$ gives the desired expansion of $z$. Finally, expansiveness ensures that every $z\in R$ satisfies the size condition for some $n\ge1$.
\end{proof}

Expansions need not be unique (for example, if $\gamma,1\in D$, the one-letter word $\gamma$ and the word $0\,1$ have the same value). We shall never use uniqueness. Nor is a complete set of residue representatives modulo $\gamma R$ necessarily admissible: in base $2$ on $\Z$, with a non-zero odd digit $e$ in place of $1$, the integer $-e$ has no finite expansion over $\{0,e\}$.

\begin{lemma}[Size bound]\label{lem:fixedlength}
Let $D\subseteq R$ be finite and contain $0$. Set
\[
 m_D:=\max_{e\in D}\max_{\sigma\in\Sinf}\frac{|\sigma e|}{|\sigma\gamma|},\qquad C_D:=\frac{m_D}{1-\rho_\gamma}.
\]
Then $|\sigma c|\le C_D|\sigma\gamma|^n$ for every $n\in\N$, every $c\in D^{(n)}$, and every $\sigma\in\Sinf$. If $D$ is admissible, every residue class in $R/\gamma^nR$ meets $D^{(n)}$.
\end{lemma}
\begin{proof}
Writing $c=\sum_{0\le j<n}\gamma^je_j$ with $e_j\in D$, we obtain
\[
|\sigma c|\le m_D|\sigma\gamma|\sum_{j<n}|\sigma\gamma|^j
\le\frac{m_D|\sigma\gamma|^{n+1}}{|\sigma\gamma|-1}
\le C_D|\sigma\gamma|^n.
\]
If $D$ is admissible and $z\in R$, the first $n$ letters of an expansion of $z$ (padded with zeros if necessary) have a value $c\in D^{(n)}$ satisfying $z-c\in\gamma^nR$.
\end{proof}

This bound is the starting point for our carry arguments. When an element $x\in R$ is written as $x=c+\gamma^n b$ with $c\in D^{(n)}$ and $b\in R$, we call $b$ the carry. If $x$ and $c$ satisfy uniform place-by-place bounds after division by $\gamma^n$, then so does the carry, which therefore belongs to a finite set (Lemma~\ref{lem:proper}).

\begin{definition}[Window set]\label{def:window}
A finite set $D\subseteq R$ containing $0$ is a \emph{window set} for $\gamma$ if, for every $n\in\N$, every residue class in $R/\gamma^nR$ meets $D^{(n)}$.
\end{definition}

\begin{lemma}\label{lem:windows}
\textup{(i)} Every admissible digit set is a window set. \textup{(ii)} If $E\subseteq R$ is a complete set of representatives of $R/\gamma R$ containing $0$, then, for every $n\in\N$, $E^{(n)}$ is a complete set of representatives of $R/\gamma^nR$. Distinct words of length $n$ over $E$ have distinct values modulo $\gamma^nR$.
\end{lemma}
\begin{proof}
Part (i) is Lemma~\ref{lem:fixedlength}. For (ii), suppose two words of length $n\ge1$ have congruent values modulo $\gamma^nR$. Their first letters are congruent modulo $\gamma R$ and hence equal. Cancelling these letters and dividing by $\gamma$, we obtain two words of length $n-1$ whose values are congruent modulo $\gamma^{n-1}R$. Induction, starting with the empty word, proves the assertion about distinct words. There are thus $|E|^n=|N_{K/\Q}(\gamma)|^n=[R:\gamma^nR]$ distinct residue classes represented.
\end{proof}

By (ii), a finite set $D\subseteq R$ containing $0$ is a window set for $\gamma$ if and only if it meets every residue class modulo $\gamma R$. The condition is necessary (take $n=1$). It is sufficient because $D$ then contains a complete set $E$ of representatives of $R/\gamma R$ with $0\in E$, and $E^{(n)}\subseteq D^{(n)}$ for every $n\in\N$.

\subsection{Automata and kernels}\label{ssec:automata}

The facts in this subsection are classical over $\Z$ \cite{AS} and over $\Z[i]$ \cite{BFY,BFK,HS}; only the place-by-place estimates are specific to the general setting. For an admissible digit set $D$, a function $f:R\to\Omega$ is \emph{$\gamma$-automatic with respect to $D$} if a deterministic finite automaton with output, over the input alphabet $D$ and reading expansions least significant digit first, computes $f(x)$ on every expansion of every $x\in R$. For a window set $D$, the \emph{$\gamma$-kernel} of $f$ at $(n,c)$, with $n\in\N$ and $c\in D^{(n)}$, is the function $R\to\Omega$ given by $y\mapsto f(c+\gamma^ny)$. We denote the set of these functions by $\cK_{\gamma,D}(f)$. By \eqref{eq:nesting}, a kernel of a kernel (with respect to the same window set) is again a kernel.

\begin{lemma}[Kernel criterion]\label{lem:kernel}
Let $D$ be admissible. A function $f:R\to\Omega$ is $\gamma$-automatic with respect to $D$ if and only if $\cK_{\gamma,D}(f)$ is finite.
\end{lemma}
\begin{proof}
This is Eilenberg's criterion \cite{Eilenberg}, \cite[\S6.6]{AS}; see also \cite[Lemma~12]{BFY}. Suppose first that $f$ is $\gamma$-automatic with respect to $D$. Two words leading to the same state give the same kernel: if $w$ has length $n$ and value $c$, and $w'$ is an expansion of $y\in R$, the concatenation $ww'$ is an expansion of $c+\gamma^ny$, and since $D$ is admissible, this tests the two kernels at every $y\in R$. Hence there are at most as many kernels as states.

Conversely, if the set of kernels is finite, take these kernels as states, with initial state $f$. Reading a digit $e\in D$ sends a state $g$ to the function $y\mapsto g(e+\gamma y)$, and the output at $g$ is $g(0)$. These transitions stay within the set of kernels by \eqref{eq:nesting}, and reading a word of value $x$ gives output $f(x)$.
\end{proof}

\begin{proposition}[Independence of the digit set]\label{prop:indep}
Let $D,D'$ be window sets for $\gamma$, and let $f:R\to\Omega$. Then $\cK_{\gamma,D}(f)$ is finite if and only if $\cK_{\gamma,D'}(f)$ is finite.
\end{proposition}
\begin{proof}
Suppose that $\cK_{\gamma,D}(f)$ is finite. For $n\in\N$ and $c'\in D'^{(n)}$, choose $c\in D^{(n)}$ with $c\equiv c'\pmod{\gamma^nR}$, and set $b=(c'-c)/\gamma^n\in R$. By Lemma~\ref{lem:fixedlength}, $|\sigma b|\le C_D+C_{D'}$ at each place. Thus $b$ belongs to a finite set $B\subseteq R$ independent of $n$ and $c'$. Since $f(c'+\gamma^ny)=f(c+\gamma^n(b+y))$ for every $y\in R$, each kernel with respect to $D'$ is a translate $y\mapsto h(b+y)$, with $b\in B$, of a kernel $h$ with respect to $D$. The reverse implication is symmetric.
\end{proof}

\begin{corollary}[Window criterion]\label{cor:windowcrit}
Let $D$ be a window set for $\gamma$. A function $f:R\to\Omega$ is $\gamma$-automatic if and only if $\cK_{\gamma,D}(f)$ is finite, independently of the choice of $D$.
\end{corollary}
From now on, \emph{$f$ is $\gamma$-automatic if its set of $\gamma$-kernels is finite with respect to one, or equivalently every, window set} (Lemmas~\ref{lem:short}, \ref{lem:kernel} and Proposition~\ref{prop:indep}). When a window set has been fixed or is clear from the context, we write $\cK_\gamma(f)$ for $\cK_{\gamma,D}(f)$; the set itself still depends on $D$.

\begin{lemma}[Powers of the base]\label{lem:powers}
Let $k\ge1$ be an integer. A function $f:R\to\Omega$ is $\gamma$-automatic if and only if it is $\gamma^k$-automatic.
\end{lemma}
\begin{proof}
If $D$ is a window set for $\gamma$, then $D^{(k)}$ is a window set for $\gamma^k$, by \eqref{eq:nesting}. The $\gamma^k$-kernels for this window set are $\gamma$-kernels for $D$. Conversely, write the level of a $\gamma$-kernel as $jk+s$, with integers $j\ge0$ and $0\le s<k$. Its representative can be written as $c'+\gamma^{jk}c''$, with $c'\in D^{(jk)}$ and $c''\in D^{(s)}$. The kernel is then $y\mapsto H(c''+\gamma^sy)$, where $H$ is a $\gamma^k$-kernel. If $f$ is $\gamma^k$-automatic, there are only finitely many choices of $H,s,c''$. (For the classical case over $\Z$, see \cite[Theorem~6.6.4]{AS}.)
\end{proof}

\begin{example}\label{ex:parity}
Automaticity depends on the base, even when one base divides the other. Let $R=\Z[\sqrt2]$, $\pi=3+\sqrt2$, and $P=\pi R$. Define $\chi:R\to\{0,1\}$ by $\chi(x)=1$ if $x\ne0$ and $v_P(x)$ is even, and $\chi(x)=0$ otherwise. Then $\chi$ is $\pi$-automatic but not $7$-automatic.

Indeed, fix a complete set $D_\pi\subseteq R$ of representatives modulo $\pi R$ containing $0$. For every $k\in\N$, the set $D_\pi^{(k)}$ is a complete set of representatives modulo $\pi^kR$ by Lemma~\ref{lem:windows}. Every non-zero representative in this set has $P$-valuation less than $k$, since the class $P^k$ is represented by $0$. The corresponding kernel is constant, while the kernel with representative $0$ is determined by the parity of $k$ (with value $0$ at the origin). There are therefore only finitely many kernels.

On the other hand, $7R=P\bar P$, where $\bar P=(3-\sqrt2)R$ is a distinct prime ideal. Take $E=\{a+b\sqrt2:a,b\in\Z,\ 0\le a,b\le6\}$. For each $k\ge1$, the set $E^{(k)}\cap P^k$ consists of $7^k$ representatives, by the Chinese remainder theorem. The corresponding kernels are pairwise distinct. To see this, take distinct $c,c'$ in this set and write $c=\pi^ka$, $c'=\pi^ka'$, with $a,a'\in R$. Put $\varpi=(3-\sqrt2)^k$, so that $7^k=\pi^k\varpi$. The element $\varpi$ is a unit modulo every power of $P$. If $m=v_P(a'-a)$, one can choose $y\in R$ with $v_P(a+\varpi y)=m+1$ by imposing a suitable congruence modulo $P^{m+2}$. Then $v_P(a'+\varpi y)=m$, and the two kernel values at $y$ differ. Hence $\chi$ has at least $7^k$ distinct $7$-kernels for every $k\ge1$, and therefore infinitely many.
\end{example}

\subsection{Carries and affine maps}\label{ssec:carries}

For a function $F:R\to\Omega$ and elements $t,\eta\in R$ with $\eta\ne0$, the \emph{section} of $F$ on $t+\eta R$ is the function $R\to\Omega$ given by $y\mapsto F(t+\eta y)$. Kernels are particular sections on residue classes $c+\gamma^nR$.

\begin{proposition}[Affine stability]\label{prop:affine}
Let $f:R\to\Omega$ be $\gamma$-automatic and let $t,\eta\in R$, with $\eta\ne0$. Then $g:R\to\Omega$, defined by $g(y):=f(t+\eta y)$, is $\gamma$-automatic. More precisely, for every window set $D$, there is a finite set $B\subseteq R$ such that every $\gamma$-kernel of $g$ with respect to $D$ is of the form $y\mapsto h(b+\eta y)$, where $h\in\cK_{\gamma,D}(f)$ and $b\in B$.
\end{proposition}
\begin{proof}
For $n\in\N$ and $c\in D^{(n)}$, choose $b_0\in D^{(n)}$ with $b_0\equiv t+\eta c\pmod{\gamma^nR}$. Set $b=(t+\eta c-b_0)/\gamma^n\in R$ (the \emph{carry}). Then $g(c+\gamma^ny)=h(b+\eta y)$ for all $y\in R$, where $h(z)=f(b_0+\gamma^nz)$. By Lemma~\ref{lem:fixedlength}, at each archimedean place,
\[
|\sigma b|\le |\sigma t|+(|\sigma\eta|+1)C_D.
\]
Thus $b$ belongs to a finite set independent of $n$ and $c$ (Lemma~\ref{lem:proper}).
\end{proof}

\begin{corollary}\label{cor:closure}
The class of $\gamma$-automatic functions is closed under translations, affine substitutions and pointwise composition $\Psi(f_1,\dots,f_k)$, where $k\ge1$ is an integer, $f_i:R\to\Omega_i$ are $\gamma$-automatic functions with finite alphabets $\Omega_i$, and $\Psi:\prod_{i=1}^k\Omega_i\to\Omega$ is a map to a finite set. In particular, it is closed under Boolean operations on characteristic functions.
\end{corollary}
\begin{proof}
The first two assertions follow from Proposition~\ref{prop:affine} (constant functions cover the case of zero slope). A kernel of the pointwise composition is obtained by applying $\Psi$ to the corresponding kernels of the $f_i$, so there are only finitely many possibilities.
\end{proof}

The next lemma, used in \S\S\ref{ssec:cyclic} and~\ref{ssec:recentre}, places an automatic function on a single residue class.

\begin{lemma}[Insertion on a residue class]\label{lem:insert}
Let $g:R\to\Omega$ be $\gamma$-automatic, let $t,\eta\in R$ with $\eta\ne0$, and choose a symbol $\bot\notin\Omega$. Define $\hat g:R\to\Omega\cup\{\bot\}$ by $\hat g(x):=g((x-t)/\eta)$ for $x\in t+\eta R$, and $\hat g(x):=\bot$ otherwise. Then $\hat g$ is $\gamma$-automatic. No coprimality assumption on $\eta$ and $\gamma$ is made.
\end{lemma}
\begin{proof}
Fix a window set $D$ for $\gamma$ and a complete set $T\subseteq R$ of representatives of $R/\eta R$. Let $n\in\N$ and $c\in D^{(n)}$. Write $y\in R$ uniquely as $s+\eta z$, with $s\in T$ and $z\in R$. Whether $c+\gamma^ny$ belongs to $t+\eta R$ depends only on $s$. If it does not, the kernel is identically $\bot$ on $s+\eta R$. Otherwise set $w=(c+\gamma^ns-t)/\eta\in R$. Choose $e\in D^{(n)}$ congruent to $w$ modulo $\gamma^nR$, and write $w=e+\gamma^nv$, with $v\in R$.

The value of the kernel at $y$ is $G(v+z)$, where $G(u)=g(e+\gamma^nu)$. At each archimedean place,
\[
|\sigma v|\le (C_D+|\sigma s|+|\sigma t|)/|\sigma\eta|+C_D.
\]
Since $T$ is finite, $v$ belongs to a finite set $V\subseteq R$ independent of $n,c,s$. The kernel is therefore determined by a finite collection of choices: for each $s\in T$, either the symbol $\bot$ or a pair $(G,v)\in\cK_{\gamma,D}(g)\times V$.
\end{proof}

\subsection{Two bases: mixed kernels}\label{ssec:mixed}

From now on, let $\alpha,\beta\in R$ be expansive. A function $f:R\to\Omega$ is \emph{doubly automatic} if it is both $\alpha$- and $\beta$-automatic.

\begin{proposition}[Mixed kernels]\label{prop:mixed}
Let $f:R\to\Omega$ be doubly automatic and let $A>0$ be a real number. The family of functions $R\to\Omega$ given by $y\mapsto f(c+\Gamma y)$, where $\Gamma=\alpha^k\beta^l$ with $k,l\in\N$, and $c\in R$ satisfies $|\sigma c|\le A|\sigma\Gamma|$ \emph{at every archimedean place}, is finite.
\end{proposition}
\begin{proof}
Fix window sets $D_\alpha,D_\beta$ with size constants $C_\alpha,C_\beta$ as in Lemma~\ref{lem:fixedlength}. Choose $c_0\in D_\alpha^{(k)}$ with $c_0\equiv c\pmod{\alpha^kR}$, and set $c_1=(c-c_0)/\alpha^k\in R$. Then $f(c+\Gamma y)=Z(c_1+\beta^ly)$, where $Z\in\cK_\alpha(f)$ is the kernel $z\mapsto f(c_0+\alpha^kz)$ and $|\sigma c_1|\le A|\sigma\beta|^l+C_\alpha$ at every place.

By Proposition~\ref{prop:affine}, applied with $\gamma=\beta$ and $\eta=\alpha^k$, every such $Z$ is $\beta$-automatic. Since $\cK_\alpha(f)$ is finite, so is the set $\cK:=\bigcup_{Z\in\cK_\alpha(f)}\cK_\beta(Z)$, where the $\beta$-kernels are taken with respect to the fixed window set $D_\beta$.

Choose $e_1\in D_\beta^{(l)}$ with $e_1\equiv c_1\pmod{\beta^lR}$, and set $c_2=(c_1-e_1)/\beta^l\in R$. We obtain $f(c+\Gamma y)=W(c_2+y)$, where $W\in\cK$ is the function $z\mapsto Z(e_1+\beta^lz)$ and $|\sigma c_2|\le A+C_\alpha+C_\beta$ at every place. Hence $c_2$ belongs to a finite set, and the family in question is finite.
\end{proof}

The condition at every archimedean place is essential to this estimate. If the absolute value of $\Gamma$ is much larger at one place than at another, the bound $\norm c\le A\norm\Gamma$ alone does not give a uniform bound for the quotient $c/\Gamma$ at the latter place. We call the sections with the place-by-place bound above the \emph{mixed kernels} of $f$ (for the chosen constant $A$).

\section{The principle of periods at infinity}\label{sec:onefar}

\subsection{Periods}\label{ssec:periods}
As in Section~\ref{sec:numeration}, let $K$ be a number field of degree $d=[K:\Q]$, let $R=\cO_K$, and let $\Omega$ be a finite alphabet. Let $f:R\to\Omega$. Set
\begin{align*}
\Per(f)&:=\{v\in R:\ f(x+v)=f(x)\text{ for all }x\in R\},\\
\Perf(f)&:=\{v\in R:\ f(x+v)=f(x)\text{ for all but finitely many }x\in R\}.
\end{align*}
These are additive subgroups of $R$. We call the non-zero elements of $\Perf(f)$ the \emph{periods at infinity} of $f$. The function $f$ is \emph{eventually periodic} if it agrees, outside a finite subset of $R$, with a function $g:R\to\Omega$ whose period group has rank $d$ (or, equivalently, finite index in $R\simeq\Z^d$). If $f=g$ outside a finite set $X_0\subseteq R$ and $\Lambda\subseteq\Per(g)$ is a subgroup of rank $d$, every non-zero $v\in\Lambda$ is a period at infinity of $f$: the possible exceptions lie in $X_0\cup(X_0-v)$. In particular, $\Perf(f)$ has rank $d$. Changing $f$ at a single point may destroy all its non-zero exact periods, but none of its periods at infinity.

\subsection{Statement}\label{ssec:onefarstatement}
Let $\alpha,\beta\in R$ be expansive elements and suppose that a doubly automatic function $f:R\to\Omega$ has a period at infinity $v\in R$. The idea is to find, among the rescaled copies $\alpha^{-k}\beta^{-l}v$ ($k,l\in\N$), a family belonging to the rational span of the periods at infinity of a single mixed kernel of $f$. We shall see that suitable families of this kind span a space of the form $\lambda F$, where $\lambda\in K^\times$ and $F\subseteq K$ is a subfield containing $\alpha^r$ and $\beta^r$ for some $r\ge1$ that we do not control. This motivates assumption (H) of the Introduction:
\[
\Q(\alpha^r,\beta^s)=K\qquad\text{for all }r,s\ge1. \tag{H}
\]
The assumption concerns both the powers and the pair of bases. For example, in $R=\Z[\sqrt2]$, the element $\alpha=14$ is rational, but $\beta=14(8+5\sqrt2)$ has different absolute values at the two real places. No positive power of $\beta$ is therefore rational, and (H) holds. An element may generate $K$ even though its square does not: in $K=\Q(\zeta_8)$, where $\zeta_8=\exp(2\pi i/8)$, the element $\zeta_8(1+2i)$ generates $K$, whereas its square $-4-3i$ lies in $\Q(i)$.

\begin{lemma}[Fields generated by powers]\label{lem:fieldpowers}
Let $\gamma\in K^\times$. Among the fields $\Q(\gamma^n)$, $n\ge1$, there is one, denoted by $K_\gamma$, which is contained in all the others. Moreover, $K_{\gamma^m}=K_\gamma$ for every $m\ge1$. For $\alpha,\beta\in K^\times$, every field $\Q(\alpha^r,\beta^s)$, $r,s\ge1$, contains the compositum $K_\alpha K_\beta$ in $K$, with equality for suitable $r,s$. Thus (H) is equivalent to $K_\alpha K_\beta=K$. Under this assumption, $\Q[\alpha^r,\beta^s]=K$ for all $r,s\ge1$.
\end{lemma}
\begin{proof}
Choose $h\ge1$ such that $[\Q(\gamma^h):\Q]$ is minimal. For any $n\ge1$, put $l=\operatorname{lcm}(h,n)$. We have $\Q(\gamma^l)\subseteq\Q(\gamma^h)$, with equality by minimality, and $\Q(\gamma^l)\subseteq\Q(\gamma^n)$. Hence $K_\gamma:=\Q(\gamma^h)\subseteq\Q(\gamma^n)$. For every $m\ge1$,
\[
K_\gamma\subseteq K_{\gamma^m}\subseteq\Q(\gamma^{mh})\subseteq K_\gamma.
\]
The assertions about the compositum follow at once. Finally, $\Q[\alpha^r,\beta^s]$ is a finite-dimensional integral domain over $\Q$, hence a field.
\end{proof}

\begin{theorem}[The principle of periods at infinity]\label{thm:onefar}
Let $\alpha,\beta\in R$ be expansive and satisfy \textup{(H)}, and let $f:R\to\Omega$ be doubly automatic. If $f$ has a period at infinity, then $\Perf(f)$ has rank $d$. If, in addition, $d\ge2$, then $f$ is eventually periodic.
\end{theorem}

For $d=1$, the first assertion is immediate. The argument generalises that of \cite{BFY}; a version for Weil bases in CM fields appears in the preprint \cite{CMW}. In the proof of \cite{BFY}, the pigeonhole principle and Krebs's pumping lemma \cite{Krebs} provide two independent periods from ``return numbers'' (\cite[Definition~13 and Lemmas~14--15]{BFY}). Here the graph of mixed kernels takes the place of these return numbers, and (H) replaces the assumption that no positive power of $\beta$ is a rational integer. The argument requires neither a relation between the absolute values of a base at different archimedean places nor assumption (Ind).

\subsection{Full rank and eventual periodicity}\label{ssec:fareventual}

\begin{proposition}[Converse]\label{prop:converse}
Let $\gamma\in R$ be expansive. Every eventually periodic function $f:R\to\Omega$ is $\gamma$-automatic.
\end{proposition}
\begin{proof}
Suppose that $f=g$ outside a finite subset of $R$, where $g:R\to\Omega$ is periodic with respect to a full-rank additive subgroup $\Lambda\subseteq R$. Choose an integer $q\ge1$ such that $qR\subseteq\Lambda$. The kernels of $g$ are all $qR$-periodic and hence belong to a finite set of functions, since $R/qR$ and $\Omega$ are finite.

Fix a window set $D$ for $\gamma$. For $z\in R$, the kernel of the characteristic function $\mathbf 1_{\{z\}}$ at $(n,c)$ is either the zero function (if $z-c\notin\gamma^nR$) or $\mathbf 1_{\{(z-c)/\gamma^n\}}$. In the latter case the point $(z-c)/\gamma^n$ lies in $R$ and satisfies
\[
|\sigma((z-c)/\gamma^n)|\le |\sigma z|+C_D
\]
at each archimedean place. There are only finitely many such points, so $\mathbf 1_{\{z\}}$ is automatic. Since $f$ is obtained from $g$ by changing its values at finitely many points, Corollary~\ref{cor:closure} gives the result. (For $\Z[i]$, see \cite[Lemma~21]{BFY}, where the automaton is constructed directly.)
\end{proof}

\begin{proposition}[Full rank]\label{prop:rank2}
Suppose that $d\ge2$ and that $\Perf(f)$ has rank $d$. Then $f$ agrees outside a finite subset of $R$ with a $\Perf(f)$-periodic function.
\end{proposition}
\begin{proof}
Put $L=\Perf(f)$ and choose a $\Z$-basis $(v_1,\ldots,v_d)$ of $L$. There is a finite set $X_0\subseteq R$ such that $f(x+v_j)=f(x)$ for every $x\in R\setminus X_0$ and every $1\le j\le d$.

Fix a coset $y_0+L$, with $y_0\in R$, and identify it with $\Z^d$ using the chosen basis. In its nearest-neighbour graph, only edges of the form $\{x,x+v_j\}$ with an endpoint in $X_0$ can join vertices with different values of $f$. All such edges have both endpoints in some cube $[-N,N]^d$, with $N\in\N$. Since $d\ge2$, the complement of the cube in $\Z^d$ is connected by nearest-neighbour edges. Thus $f$ is constant on that complement. Assigning this constant value to the whole coset gives a function agreeing with $f$ except at finitely many points of the coset. There are only finitely many cosets of $L$ in $R$, so the resulting function on $R$ has the required properties. (Together with Lemma~\ref{lem:hered}, this is the analogue of the transfer of periods in \cite[Lemmas~18--19]{BFY}.)
\end{proof}

\begin{remark}\label{rem:degreeone}
The preceding proposition fails for $d=1$: the function $f=\mathbf 1_{\Z_{\ge0}}$ on $\Z$ satisfies $\Perf(f)=\Z$ but is not eventually periodic in our sense. The complement of a finite interval has two connected components. This explains why the version of the main theorem valid in all degrees $d\ge1$ is stated in terms of periods at infinity (Corollary~\ref{cor:farperiods}).
\end{remark}

It remains to prove the first assertion of Theorem~\ref{thm:onefar}. We organise the mixed kernels into a finite directed graph (\S\ref{ssec:graph}), whose paths from $f$ record scales at which the rescaled period lies in the rational span of the periods at infinity of the kernel reached (\S\ref{ssec:heredity}). A vertex reached by a path that admits insertions of closed paths in two independent directions has a full-rank group of periods at infinity, by (H); the exponents of the paths to the other vertices lie in bounded neighbourhoods of finitely many rays (\S\ref{ssec:goodbad}). At a single scale outside these neighbourhoods, all residue classes therefore lead to kernels with full-rank groups of periods at infinity, which we combine (\S\ref{ssec:proofonefar}).

\subsection{The graph of mixed kernels}\label{ssec:graph}
Throughout the rest of \S\ref{sec:onefar}, let $\alpha,\beta\in R$ be expansive elements and let $f:R\to\Omega$ be doubly automatic. Fix window sets $D_\alpha,D_\beta$ for $\alpha,\beta$ respectively, with size constants $C_\alpha,C_\beta$ as in Lemma~\ref{lem:fixedlength}, and put $C:=C_\alpha+C_\beta$.

\begin{definition}[The graph $\cG$]\label{def:graph}
A pair $(\Gamma,c)\in R\times R$ is called \emph{bounded} if $\Gamma=\alpha^k\beta^l$ for some $k,l\in\N$, and $|\sigma c|\le C|\sigma\Gamma|$ for every $\sigma\in\Sinf$. The vertices of $\cG$ are the functions $f_{\Gamma,c}:R\to\Omega$, $y\mapsto f(c+\Gamma y)$, with $(\Gamma,c)$ bounded. By Proposition~\ref{prop:mixed}, their number, denoted by $|\cG|$, is finite. We regard $\cG$ as a directed multigraph. For each vertex $g$ and digit $e\in D_\alpha$, draw an edge of the first type, labelled $(\alpha,e)$, to the function $y\mapsto g(e+\alpha y)$; for each $e\in D_\beta$, draw an edge of the second type, labelled $(\beta,e)$, to $y\mapsto g(e+\beta y)$. The two types remain formally distinct even if $\alpha=\beta$, and distinct digits give distinct edges even if their endpoints coincide. Lemma~\ref{lem:graph}\,(a) below shows that the endpoint is again a vertex.
\end{definition}
The \emph{exponent} of a directed path is the pair $(k,l)\in\N^2$ counting its edges of the first and second types, respectively. For a vertex $g$, let $E_g\subseteq\N^2$ be the set of exponents of directed paths from $f$ to $g$. Paths may repeat vertices and edges, and we allow the path of length zero.

\begin{lemma}[Paths]\label{lem:graph}
\textup{(a)} If $(\Gamma,c)$ is bounded, $g=f_{\Gamma,c}$, $\gamma\in\{\alpha,\beta\}$ and $e\in D_\gamma$, then $g(e+\gamma y)=f_{\gamma\Gamma,c+\Gamma e}(y)$ for $y\in R$, with $(\gamma\Gamma,c+\Gamma e)$ bounded. In particular, $f=f_{1,0}$ is a vertex.

\textup{(b)} A path of exponent $(k,l)$ starting at $f_{\Gamma,c}$, with $(\Gamma,c)$ bounded, ends at a vertex $f_{\Gamma',c'}$ with $(\Gamma',c')$ bounded, $\Gamma'=\alpha^k\beta^l\Gamma$ and $c'\equiv c\pmod{\Gamma R}$.

\textup{(c)} For all $k,l\in\N$ and every $c^*\in R$, there is a path from $f$ of exponent $(k,l)$ ending at $f_{\Gamma,c}$, where $\Gamma=\alpha^k\beta^l$, $(\Gamma,c)$ is bounded and $c\equiv c^*\pmod{\Gamma R}$.
\end{lemma}
\begin{proof}
(a) The identity follows by composition. At each archimedean place,
\[
|\sigma(c+\Gamma e)|\le(C+m_{D_\gamma}|\sigma\gamma|)|\sigma\Gamma|
\le C|\sigma\gamma\Gamma|,
\]
since $C\ge C_\gamma\ge m_{D_\gamma}|\sigma\gamma|/(|\sigma\gamma|-1)$. The pair $(1,0)$ is bounded.

(b) Induction on the length of the path, using (a).

(c) Choose $c_1\in D_\alpha^{(k)}$ congruent to $c^*$ modulo $\alpha^kR$, and then $c_2\in D_\beta^{(l)}$ congruent to $(c^*-c_1)/\alpha^k$ modulo $\beta^lR$. Set $c=c_1+\alpha^kc_2\in R$. By construction $c\equiv c^*\pmod{\Gamma R}$. Reading first a word over $D_\alpha$ of length $k$ and value $c_1$, then a word over $D_\beta$ of length $l$ and value $c_2$, gives a path of exponent $(k,l)$ ending at $f_{\Gamma,c}$. Part (a) ensures that $(\Gamma,c)$ is bounded.
\end{proof}
In part (c), the exponent can be prescribed independently of the residue class. No relation between $\alpha$ and $\beta$, not even coprimality, is needed.

\subsection{Inheritance of periods}\label{ssec:heredity}
For a function $g:R\to\Omega$, let $\Pi(g):=\operatorname{span}_{\Q}\Perf(g)\subseteq K$. This is a $\Q$-vector subspace of $K$. The group $\Perf(g)$ has rank $d$ if and only if $\Pi(g)=K$.

\begin{lemma}[Inheritance]\label{lem:hered}
Let $v\in\Perf(f)$, $\Gamma\in R\setminus\{0\}$ and $c\in R$. If $g(y)=f(c+\Gamma y)$ for $y\in R$, then $\Gamma^{-1}v\in\Pi(g)$.
\end{lemma}
\begin{proof}
Let $m=|N_{K/\Q}(\Gamma)|=[R:\Gamma R]$. Since the finite additive group $R/\Gamma R$ is annihilated by $m$, we have $m\in\Gamma R$. Thus $w=m\Gamma^{-1}v\in R$. As $mv\in\Perf(f)$,
\[
g(y+w)=f(c+\Gamma y+mv)=g(y)
\]
for all but finitely many $y\in R$ (the map $y\mapsto c+\Gamma y$ is injective). Hence $w\in\Perf(g)$, and $\Gamma^{-1}v=w/m\in\Pi(g)$.
\end{proof}
It follows from Lemma~\ref{lem:graph}\,(b) that, for every $v\in\Perf(f)$, every vertex $g$ and every $(k,l)\in E_g$, we have $\alpha^{-k}\beta^{-l}v\in\Pi(g)$.

\subsection{Good and bad vertices}\label{ssec:goodbad}
For a strongly connected component $\cC$ of $\cG$, let $U(\cC)\subseteq\N^2$ be the set of exponents of closed paths lying in $\cC$. We call $\cC$ \emph{two-dimensional} if $U(\cC)$ spans $\R^2$, \emph{one-dimensional} if it spans a line, and \emph{trivial} otherwise (a single vertex without a loop). The direction of a one-dimensional component is the ray in $\R_{\ge0}^2$ containing $U(\cC)\setminus\{(0,0)\}$. Let $\mathcal U$ be the finite set of these directions.

We call a vertex $g$ \emph{good} if some path from $f$ to $g$ meets either a two-dimensional component or two one-dimensional components of different directions. Otherwise we call it \emph{bad}.

\begin{lemma}[Good vertices]\label{lem:good}
Assume \textup{(H)}. If $f$ has a period at infinity $v$ and $g$ is a good vertex of $\cG$, then $\Perf(g)$ has rank $d$.
\end{lemma}
\begin{proof}
Choose a path $\pi$ from $f$ to $g$ witnessing goodness, and denote its exponent by $(k_0,l_0)\in\N^2$.

\emph{Two closed paths with independent exponents.} Suppose first that $\pi$ meets a vertex $Z$ of a two-dimensional component $\cC$. For any closed path $W$ based at a vertex $Z'\in\cC$, strong connectivity gives paths $A_1$ from $Z$ to $Z'$ and $A_2$ from $Z'$ to $Z$, both in $\cC$. The exponents of the closed paths $A_1A_2$ and $A_1WA_2$ at $Z$ differ by the exponent of $W$. Consequently the real span of the exponents of closed paths based at $Z$ contains $U(\cC)$, and so is all of $\R^2$. We can therefore choose two closed paths based at $Z$ with linearly independent exponents.

If instead $\pi$ meets vertices $Z_1,Z_2$ in two one-dimensional components of different directions, the same construction gives a closed path of positive length at each $Z_i$. Its exponent is non-zero and lies on the direction of that component. The two exponents are again linearly independent. In either case, inserting arbitrary numbers $j_1,j_2\in\N$ of copies of the chosen closed paths into $\pi$ gives
\[
E_g\supseteq(k_0,l_0)+\N(k_1,l_1)+\N(k_2,l_2),
\qquad \Delta:=k_1l_2-k_2l_1\in\Z\setminus\{0\},
\]
where $(k_i,l_i)\in\N^2$ are the exponents of the two closed paths.

\emph{The generated field.} Set $\Gamma_0=\alpha^{k_0}\beta^{l_0}\in R$ and $\eta_i=\alpha^{-k_i}\beta^{-l_i}\in K^\times$ for $i=1,2$. For $j_1,j_2\in\N$, the preceding inclusion gives a path to $g$ with scale $\Gamma_0\eta_1^{-j_1}\eta_2^{-j_2}$. By inheritance (Lemma~\ref{lem:hered}),
\[
 \Gamma_0^{-1}\eta_1^{j_1}\eta_2^{j_2}v\in\Pi(g).
\]
Taking $\Q$-linear combinations gives
\[
\Pi(g)\supseteq\Gamma_0^{-1}v\,\Q[\eta_1,\eta_2].
\]
The algebra $\Q[\eta_1,\eta_2]\subseteq K$ is an integral domain of finite dimension over $\Q$, hence a field. In particular, it contains
\[
\eta_1^{k_2}\eta_2^{-k_1}=\beta^{\Delta},\qquad
\eta_1^{l_2}\eta_2^{-l_1}=\alpha^{-\Delta}.
\]
It thus contains $\Q(\alpha^{|\Delta|},\beta^{|\Delta|})=K$ by (H). Since $v\ne0$, multiplication by $\Gamma_0^{-1}v$ is an automorphism of the $\Q$-vector space $K$. Hence $\Pi(g)=K$.
\end{proof}

\begin{lemma}[Bad vertices]\label{lem:bad}
Let $\cR:=\{(0,0)\}\cup\bigcup_{u\in\mathcal U}u\subseteq\R^2$. For every bad vertex $g$, every point of $E_g$ lies at Euclidean distance at most $|\cG|$ from $\cR$.
\end{lemma}
\begin{proof}
Take a path $\pi$ from $f$ to $g$. Successively remove closed subpaths of positive length until no vertex is repeated. The remaining path has length at most $|\cG|$, and its exponent $(k_0,l_0)\in\N^2$ has Euclidean norm at most $|\cG|$. Each removed subpath lies in a non-trivial strongly connected component visited by $\pi$. Since $g$ is bad, all these components are one-dimensional with the same direction $u\in\mathcal U$. Thus the exponent of $\pi$ is the sum of $(k_0,l_0)$ and a point of $u$, or just $(k_0,l_0)$ if there were no closed subpaths to remove. The asserted bound follows.
\end{proof}

\subsection{Proof of Theorem~\ref{thm:onefar}}\label{ssec:proofonefar}
Let $v\in R\setminus\{0\}$ be a period at infinity of $f$. Choose $w_0\in\N^2\setminus\{(0,0)\}$ lying on none of the finitely many rays in $\mathcal U$. The set $\cR$ is closed and invariant under positive dilations, so $\operatorname{dist}(w_0,\cR)>0$ and
\[
\operatorname{dist}(tw_0,\cR)=t\operatorname{dist}(w_0,\cR)\qquad(t>0).
\]
Choose an integer $t\ge1$ such that $\operatorname{dist}(tw_0,\cR)>|\cG|$, and put $(k,l)=tw_0\in\N^2$ and $\Gamma=\alpha^k\beta^l\in R$.

By Lemma~\ref{lem:graph}\,(c), there is a complete set $\mathcal S\subseteq R$ of representatives of $R/\Gamma R$ such that, for every $c\in\mathcal S$, the function $g_c:y\mapsto f(c+\Gamma y)$ is a vertex reached by a path of exponent $(k,l)$. By Lemma~\ref{lem:bad}, this vertex must be good. Thus $\Perf(g_c)$ has rank $d$ by Lemma~\ref{lem:good}. Since $\mathcal S$ is finite, the additive subgroup $L:=\bigcap_{c\in\mathcal S}\Perf(g_c)$ has finite index in $R$.

Take $w\in L$. On each residue class, writing $x=c+\Gamma y$ with $y\in R$, we have
\[
f(x+\Gamma w)=g_c(y+w)=g_c(y)=f(x)
\]
except for $y$ in a finite set $X_c\subseteq R$ (which may depend on $w$). Since there are only finitely many classes, the exceptional points $x$ form a finite set. Hence $\Gamma w\in\Perf(f)$. We conclude that $\Perf(f)$ contains the full-rank subgroup $\Gamma L$ of $R$, as required. \qed

The proof uses only finiteness of the mixed kernels and (H), the latter in Lemma~\ref{lem:good}. It therefore applies in both cases of the proof of the main theorem, distinguished by the \emph{supports} of the bases, that is, the sets of prime ideals of $R$ containing them.

\section{Orbit closures and the rigid class}\label{sec:orbits}

Let $K$ be a number field of degree $d=[K:\Q]$, let $R=\cO_K$, and let $\Omega$ be a finite alphabet. Under (H), with $d\ge2$, a doubly automatic set that is not eventually periodic has no period at infinity (Theorem~\ref{thm:onefar}). To study its behaviour at infinity, we consider the orbit closure of its characteristic function. More generally, for a function $f:R\to\Omega$, we set
\[
 \cX_f:=\overline{\{S^q f:q\in R\}}\subseteq\Omega^R,
 \qquad (S^qu)(x):=u(x+q)\quad(q,x\in R).
\]
Here $\Omega^R$ carries the product topology, with $\Omega$ discrete. A \emph{subsystem} is a non-empty closed subset invariant under all translations $S^q$; it is \emph{minimal} if it contains no proper subsystem. This is the point of view adopted by Durand \cite{DurandSemenov,Durand}. We use only elementary notions of symbolic dynamics, in particular minimal subsystems and constant-shape substitutions. As before, $\alpha,\beta\in R$ are expansive. Condition (Ind) is not used in this section.

\subsection{Minimal subsystems}\label{ssec:shifts}
A \emph{pattern} over a finite alphabet is a function on a finite subset of $R$, called its support. A pattern $w$ occurs in a configuration $u\in\Omega^R$ at $q\in R$ if $u(q+x)=w(x)$ for every $x$ in its support. We use balls in $R$ for the norm of \S\ref{ssec:places}.

A minimal subsystem $Y\subseteq\Omega^R$ is \emph{uniformly recurrent}: for each pattern $w$ occurring in a point of $Y$, there is a radius $r=r(w)>0$ such that $w$ occurs in the restriction of every point of $Y$ to every translated ball of radius $r$. Indeed, the points in which $w$ never occurs form a closed invariant subset, which must be empty by minimality. The open sets specifying an occurrence of $w$ at a given position thus cover $Y$. By compactness, finitely many positions suffice. Enlarging the radius to include the supports at all these positions, and using translation invariance, gives the assertion.

Set
\[
 \Lambda_Y:=\{q\in R:S^qu=u\text{ for all }u\in Y\}.
\]
This is the stabiliser of any point of $Y$: if a translation fixes one point, it fixes its dense orbit and hence all of $Y$. We say that $Y$ is \emph{aperiodic} if $\Lambda_Y=0$.

\begin{lemma}\label{lem:rankfinite}
Let $Y\subseteq\Omega^R$ be minimal. Then $\Lambda_Y$ has rank $d$ if and only if $Y$ is a finite orbit.
\end{lemma}
\begin{proof}
If $\Lambda_Y$ has rank $d$, it has finite index in $R$. The orbit of any point of $Y$ is then finite, hence closed, and by minimality it equals $Y$. Conversely, the stabiliser of a point with finite orbit has finite index in $R$.
\end{proof}

\subsection{Adapted windows and constant-shape substitutions}\label{ssec:window}
Let $\gamma\in R$ be expansive. An \emph{adapted window} for $\gamma$ is a pair $(D,\lambda)$, where $\lambda$ is a positive power of $\gamma$ and $D\subseteq R$ is a complete set of representatives of $R/\lambda R$ containing $0$, admissible for $\lambda$. Here and below $D^{(n)}$ is formed using the base $\lambda$. By Lemma~\ref{lem:windows}, for each $n\in\N$ the sets $\lambda^ny+D^{(n)}$ ($y\in R$), called \emph{level-$n$ tiles}, partition $R$. The sets $D^{(n)}$ are nested and their union is $R$. Moreover, if $(D,\lambda)$ is adapted, so is $(D^{(j)},\lambda^j)$ for every $j\ge1$: the grouped digits give the sets $D^{(jn)}$ at level $n$.

\begin{proposition}\label{prop:window}
Every expansive element $\gamma\in R$ admits an adapted window.
\end{proposition}
\begin{proof}
Choose a half-open fundamental parallelepiped $\mathcal P$ for the lattice $R$ in the Minkowski space $K_{\R}$, with $0$ in its interior. Put $C_0:=\sup_{x\in\mathcal P}\norm x$ and choose $r_0>0$ such that the closed ball of radius $r_0$ about $0$ is contained in $\mathcal P$; thus $r_0\le C_0$. Set $B:=4C_0+1$. Since $\gamma$ is expansive, there is a $k\ge1$ such that $\lambda:=\gamma^k$ satisfies
\[
 m_\lambda:=\min_{\sigma\in\Sinf}|\sigma\lambda|>B/r_0>4.
\]
The set $D:=R\cap\lambda\mathcal P$ is a complete set of representatives of $R/\lambda R$ containing $0$. Moreover, if $z\in R$ and $\norm z\le B$, then $\norm{z/\lambda}\le B/m_\lambda<r_0$, so $z\in D$.

For $z\in R$, let $\dig z\in D$ be its representative and set $\delta(z):=(z-\dig z)/\lambda\in R$, so that $z=\dig z+\lambda\delta(z)$. Since $\dig z/\lambda\in\mathcal P$ and $m_\lambda>4$, if $\norm z>B$ then $C_0<\norm z/4$ and
\[
 \norm{\delta(z)}\le\frac{\norm z}{m_\lambda}+C_0<\frac{\norm z}{2}.
\]
The iterates of $\delta$ therefore enter the ball of radius $B$ after finitely many steps. In that ball $\dig z=z$, and the next quotient is $0$. Iterating the identity $z=\dig z+\lambda\delta(z)$ gives a finite expansion of every $z\in R$ over $D$. Hence $D$ is admissible.
\end{proof}

By Lemma~\ref{lem:powers}, a function is $\gamma$-automatic if and only if it is $\lambda$-automatic. We may therefore pass freely to the base $\lambda$.

\begin{definition}\label{def:substitution}
Let $(E,\mu)$ be an adapted window and let $\bA$ be a finite non-empty alphabet. A \emph{constant-shape substitution} with expansion $\mu$ and shape $E$ is a map $\theta:\bA\to\bA^E$. It extends to configurations in $\bA^R$ by
\[
 \theta(u)_{\mu y+b}=\theta(u_y)_b\qquad(y\in R,\ b\in E),
\]
where $u_y=u(y)$. The same rule extends $\theta$ to a pattern $P$ of finite support $C\subseteq R$, with image supported on $\mu C+E$.

The $n$-th iterate has expansion $\mu^n$ and shape $E^{(n)}$. For $a,b\in\N$, $g\in\bA$, $c\in E^{(a)}$ and $c'\in E^{(b)}$,
\[
 \theta^{a+b}(g)_{c+\mu^ac'}=\theta^a\bigl(\theta^b(g)_{c'}\bigr)_c.
\]
A pattern is \emph{legal} if it occurs in $\theta^n(g)$ for some $n\in\N$ and $g\in\bA$. Occurrence in a finite pattern requires the translated support to be contained in its support. The \emph{substitution subshift} $X_\theta\subseteq\bA^R$ consists of the configurations all of whose finite patterns are legal. The substitution is \emph{primitive} if there is a $p\ge1$ such that every letter occurs in $\theta^p(g)$ for every $g\in\bA$.
\end{definition}

\begin{lemma}[Fixed point]\label{lem:germ}
Let $\theta$ be a constant-shape substitution, not necessarily primitive. There exist a letter $g\in\bA$, an integer $M\ge1$ and a point $z\in X_\theta$ such that $\theta^M(z)=z$ and $z|_{E^{(Mn)}}=\theta^{Mn}(g)$ for every $n\in\N$.
\end{lemma}
\begin{proof}
The map $g\mapsto\theta(g)_0$ of the finite set $\bA$ has a periodic point. Choose a letter $g$ on a cycle of length $M$, so that $\theta^M(g)_0=g$. Since $E^{(M(n+1))}=E^{(Mn)}+\mu^{Mn}E^{(M)}$, the block
\[
 \theta^{M(n+1)}(g)=\theta^{Mn}(\theta^M(g))
\]
is made up of the blocks $\theta^{Mn}(\theta^M(g)_y)$ on $\mu^{Mn}y+E^{(Mn)}$, for $y\in E^{(M)}$. The block indexed by $y=0$ is $\theta^{Mn}(g)$. Thus the blocks $\theta^{Mn}(g)$ extend one another. Their supports exhaust $R$, so they define a point $z\in\bA^R$. Every finite pattern of $z$ is contained in one of these blocks; hence $z\in X_\theta$.

Finally, for $y\in E^{(Mn)}$ and $b\in E^{(M)}$, using $\theta^{M(n+1)}=\theta^M\circ\theta^{Mn}$ gives
\[
 z(\mu^My+b)=\theta^M\bigl(\theta^{Mn}(g)_y\bigr)_b
 =\theta^M(z_y)_b=\theta^M(z)_{\mu^My+b}.
\]
As the sets $E^{(Mn)}$ exhaust $R$, this proves $\theta^M(z)=z$.
\end{proof}

\begin{lemma}[Minimality]\label{lem:minimal}
If $\theta$ is primitive, $X_\theta$ is non-empty and minimal, and every legal pattern occurs in every point of $X_\theta$.
\end{lemma}
\begin{proof}
The set $X_\theta$ is non-empty by Lemma~\ref{lem:germ}. It is closed and invariant under translations by its definition in terms of finite patterns. Let $w$ be a legal pattern occurring in $\theta^n(g_1)$, and let $p$ be as in the definition of primitivity. For every $m\ge n+p$ and every letter $g$, the letter $g_1$ occurs in $\theta^{m-n}(g)$: an iterate of order at least $p$ is made up of blocks of order $p$, each containing every letter. Consequently $w$ occurs in $\theta^m(g)$. The same argument shows that a legal pattern occurring in $\theta^{N'}(g')$ also occurs in $\theta^{N'+i}(g'')$ for every $i\ge p$ and every letter $g''$.

Fix $m:=n+p$. Choose $b_m>0$ bounding the norms of the elements of $E^{(m)}$ and set $V:=\{c\in R:\norm c\le2b_m\}$. Let $x\in X_\theta$. The restriction $x|_V$ is legal. By the preceding paragraph, it occurs in a block $\theta^N(g')$ with $N\ge m$. In other words, there is a $q\in R$ such that $q+V\subseteq E^{(N)}$ and $\theta^N(g')(q+c)=x(c)$ for $c\in V$.

The support of $\theta^N(g')=\theta^m(\theta^{N-m}(g'))$ is partitioned into translates $\mu^mc'+E^{(m)}$, each of diameter at most $2b_m$ and each carrying a block $\theta^m(g'')$. The tile containing $q$ is contained in $q+V$. Thus $x$ contains a block $\theta^m(g'')$, and hence an occurrence of $w$.

Every legal pattern therefore occurs in every point of $X_\theta$. This implies that each point has a dense orbit in $X_\theta$, proving minimality.
\end{proof}

\subsection{The kernel substitution}\label{ssec:kernelsubst}
We encode the kernels of an automatic function as a constant-shape substitution; the minimal subsystems of its orbit closure lift to subshifts of primitive substitutions extracted from it.

Let $f:R\to\Omega$ be $\gamma$-automatic and let $(D,\lambda)$ be an adapted window for $\gamma$. Let $\cK$ be the finite set of its $\lambda$-kernels formed using $D$. For $G\in\cK$ and $b\in D$, the function $G_b:y\mapsto G(b+\lambda y)$ is again in $\cK$.

We record the values of all these kernels at once. Set $\bA:=\Omega^{\cK}$ and define
\[
 \zeta(g)_b(G):=g(G_b),\qquad
 \Theta(y):=(G(y))_{G\in\cK},\qquad
 \pr(g):=g(f),
\]
for $g\in\bA$, $b\in D$, $G\in\cK$ and $y\in R$. Thus $\zeta:\bA\to\bA^D$ is a substitution with expansion $\lambda$ and shape $D$, $\Theta:R\to\bA$ is a configuration, and $\pr:\bA\to\Omega$ is the projection onto the coordinate indexed by $f$. For $n\in\N$, $c\in D^{(n)}$ and $y\in R$,
\[
 \Theta(c+\lambda^ny)=\zeta^n(\Theta(y))_c.
\]
Indeed, the formula for $n=1$ follows coordinate by coordinate from the definition of $G_b$, and iteration gives the general case. In particular $\zeta^n(\Theta)=\Theta$. Extending $\pr$ coordinatewise to configurations, we have $f=\pr\circ\Theta$ and $\pr(\cX_\Theta)=\cX_f$, by continuity, equivariance and compactness.

We shall use the following desubstitution property: for each $T\ge1$, every $z\in\cX_\Theta$ has a representation
\[
 z=S^c\zeta^T(z'),\qquad c\in D^{(T)},\quad z'\in\cX_\Theta.
\]
To see this, write $z=\lim_j S^{q_j}\Theta$ and $q_j=\lambda^Tb_j+c_j$, with $c_j\in D^{(T)}$ and $b_j\in R$. Passing to a subsequence, we may suppose that $c_j=c$ is constant and $S^{b_j}\Theta\to z'\in\cX_\Theta$. The identity $\zeta^T(S^bu)=S^{\lambda^Tb}\zeta^T(u)$, together with $\zeta^T(\Theta)=\Theta$, gives the required formula by continuity.

The \emph{letter graph} $\cH$ has vertex set $\bA$ and an edge $g\to g'$ if $g'$ occurs in $\zeta(g)$. Every letter has a successor. A \emph{sink} is a strongly connected component from which no edge leads out. Let $C$ be a sink of period $h$ (the greatest common divisor of the lengths of its closed paths of positive length), and let $C_0,\ldots,C_{h-1}$ be its cyclic classes, with edges from $C_i$ to $C_{i+1}$ (indices modulo $h$). The letters occurring in $\zeta^n(g)$, for $g\in C_i$, are precisely the endpoints of paths of length $n$ starting at $g$. They all lie in $C_{i+n}$. The substitutions
\[
 \theta_i:=\zeta^h|_{C_i}:C_i\longrightarrow C_i^{D^{(h)}}
\]
are primitive by the Perron--Frobenius theorem for irreducible matrices of period $h$ \cite[Ch.~8]{HJ}.

\begin{lemma}[Lifting]\label{lem:lift}
Let $Z\subseteq\cX_\Theta$ be minimal. There exist a sink $C$ and a cyclic class $C_i$ such that $Z=X_{\theta_i}$. If $\pr(Z)$ is aperiodic, so is $Z$. Consequently $\cX_f$ has at most $|\bA|$ minimal subsystems, each of the form $\pr(X_\theta)$, with $\theta$ primitive of expansion $\lambda^h$ and shape $D^{(h)}$ for some $h\ge1$.
\end{lemma}
\begin{proof}
Fix $z\in Z$ and put $\tau:=|\bA|$.

\emph{Every pattern $w$ of $z$ occurs in a block $\zeta^n(g_0)$ with $g_0$ in a sink.} From every letter, a path of length less than $\tau$ leads to a sink, and it can be extended within the sink to length $\tau$; so every block $\zeta^\tau(g)$ contains a letter of a sink. By uniform recurrence, $w$ occurs in every translated ball of some radius $r$ in $z$. Choose $T>\tau$ such that $D^{(T-\tau)}$ contains the ball of radius $r$ about $0$, and write $z=S^c\zeta^T(z')$ (desubstitution). For a letter $g$ of $z'$, the block $\zeta^T(g)=\zeta^{T-\tau}(\zeta^\tau(g))$ occurs in $z$ and contains a block $\zeta^{T-\tau}(g_0)$, with $g_0$ in a sink, supported on a translate of $D^{(T-\tau)}$. This block contains a translated ball of radius $r$, hence an occurrence of $w$.

\emph{A single cyclic class.} The letters of a block $\zeta^n(g_0)$ lie in one cyclic class of the sink of $g_0$. Applying the preceding paragraph to the restrictions of $z$ to increasing balls, we find that all letters of $z$ lie in one cyclic class $C_i$ of one sink $C$ (all these classes contain the letter $z(0)$, and distinct cyclic classes are disjoint). In particular, every non-empty pattern of $z$ occurs in a block $\zeta^n(g_0)$ with $g_0\in C$.

\emph{Legality and minimality.} If a non-empty pattern $w$ of $z$ occurs in $\zeta^n(g_0)$ with $g_0\in C_l$, then $l+n\equiv i\pmod h$. A path from some $g'\in C_i$ to $g_0$ has length $\ell\equiv l-i\pmod h$, so $w$ occurs in $\zeta^{n+\ell}(g')$ with $h\mid n+\ell$, that is, in a block of an iterate of $\theta_i$. Thus every pattern of $z$ is legal for $\theta_i$, and $Z$, the orbit closure of $z$, is contained in $X_{\theta_i}$. As $X_{\theta_i}$ is minimal (Lemma~\ref{lem:minimal}), $Z=X_{\theta_i}$.

\emph{Consequences.} A non-zero translation fixing a point of $Z$ fixes every point of $Z$ and of $\pr(Z)$; so $Z$ is aperiodic if $\pr(Z)$ is. Distinct pairs formed by a sink and one of its cyclic classes have disjoint alphabets, so $\cX_\Theta$ has at most $|\bA|$ minimal subsystems. Finally, if $Y\subseteq\cX_f$ is minimal, the set $\pr^{-1}(Y)\cap\cX_\Theta$ is non-empty (as $\pr(\cX_\Theta)=\cX_f$), closed and invariant. It contains a minimal subsystem $Z$, and $\pr(Z)$, a non-empty closed invariant subset of $Y$, equals $Y$.
\end{proof}

To compare the periods of a minimal subsystem in the two bases (Theorem~\ref{thm:ranks}), we need an automatic point in it.

\begin{lemma}[Automatic generator]\label{lem:generator}
Let $Y=\pr(X_\theta)\subseteq\cX_f$ be a minimal subsystem, with $\theta$ of expansion $\mu$ and shape $E$ as in Lemma~\ref{lem:lift}. There is a point $y\in Y$ with orbit closure $Y$ and with only finitely many kernels $x\mapsto y(c+\mu^nx)$ ($n\in\N$, $c\in E^{(n)}$). In particular $y$ is $\gamma$-automatic, and $\Lambda_Y=\Per(y)$.
\end{lemma}
\begin{proof}
By Lemma~\ref{lem:germ}, choose $z\in X_\theta$ fixed by $\theta^M$ for some $M\ge1$. Its orbit closure is the minimal subshift $X_\theta$, so the orbit closure of $y:=\pr(z)$ is $\pr(X_\theta)=Y$. Since $\theta^{iM}(z)=z$, for $i\in\N$ and $b\in E^{(iM)}$ we have
\[
 y(\mu^{iM}x+b)=F_{i,b}(z(x)),\qquad
 F_{i,b}(g):=\pr\bigl(\theta^{iM}(g)_b\bigr)
 \quad(x\in R).
\]
There are only finitely many distinct maps $F_{i,b}$, since they map the finite alphabet of $\theta$ into $\Omega$. For an arbitrary $N\in\N$, write $N=iM+s$ with $0\le s<M$. Each element of $E^{(N)}$ has the form $b+\mu^{iM}c'$ with $b\in E^{(iM)}$ and $c'\in E^{(s)}$. Then
\[
 y(\mu^Nx+b+\mu^{iM}c')=F_{i,b}(z(\mu^sx+c')).
\]
The pairs $(s,c')$ lie in a finite set independent of $i$, and the maps $F_{i,b}$ also form a finite family. Thus $y$ has only finitely many kernels. Since $E$ is a window set for $\mu$, $y$ is $\mu$-automatic, hence $\gamma$-automatic by Lemma~\ref{lem:powers}. Finally, $\Lambda_Y$ is the stabiliser of $y$, namely $\Per(y)$.
\end{proof}

\subsection{Recognisability}\label{ssec:recog}
The odometer factor of \S\ref{ssec:odometer} rests on a property of a primitive substitution $\theta$ with aperiodic subshift: in a fixed point, membership of a position in $\mu R$ can be read from a bounded neighbourhood of that position. This is the first of the two conditions defining recognisability \cite[Definition~3.4, p.~1279]{Cab}. It is stated by Cabezas \cite[Proposition~3.7, p.~1281]{Cab}, as the special case of the identity map in a result on symbolic factors, and by Cabezas and Leroy \cite[Proposition~5.6]{CabezasLeroy}. We prove the case we need. Following their idea, we compare, at two different levels, the two block decompositions of the same point that a failure of the property would yield. The next three lemmas are used only in the proof of Proposition~\ref{prop:recog}.

\emph{Setting and patterns.} In this subsection, $(E,\mu)$ is an adapted window for an expansive element $\gamma\in R$, and $\theta:\bA\to\bA^E$ is a constant-shape substitution with expansion $\mu$ and shape $E$. The estimates below use only expansiveness of $\mu$ and the existence and uniqueness, for $n\in\N$, of the decomposition $x=\mu^nc+e$ of any $x\in R$, with $c\in R$ and $e\in E^{(n)}$.

Recall that a pattern $P$ is a map from a finite subset $\supp P\subseteq R$ to $\bA$. It occurs in $u\in\bA^R$ at $q\in R$ if $u(q+x)=P(x)$ for every $x\in\supp P$, and in a pattern $W$ at $q$ if also $q+\supp P\subseteq\supp W$. Its image $\theta^n(P)$ has support $\mu^n\supp P+E^{(n)}$ and value $\theta^n(P(c))_e$ at $\mu^nc+e$, for $c\in\supp P$ and $e\in E^{(n)}$. The formula for iterates in Definition~\ref{def:substitution} gives the following identities and occurrence properties (with $a,b,n\in\N$):
\begin{enumerate}[label=\textup{(B\arabic*)},nosep]
\item $\theta^a(\theta^b(P))=\theta^{a+b}(P)$;
\item if $P$ occurs in $u$ (or in a pattern $W$) at $q$, then $\theta^n(P)$ occurs in $\theta^n(u)$ (or in $\theta^n(W)$) at $\mu^nq$;
\item $\theta^n(S^qu)=S^{\mu^nq}\theta^n(u)$;
\item if $u=S^q\theta^n(u')$, where $u'\in\bA^R$, and $C\subseteq R$ is finite, then $\theta^n(u'|_C)$ occurs in $u$ at $-q$.
\end{enumerate}
For (B3), both sides have value $\theta^n(u(c+q))_e$ at $\mu^nc+e$, for $c\in R$ and $e\in E^{(n)}$. For (B4), we have $u(-q+\mu^nc+e)=\theta^n(u')(\mu^nc+e)=\theta^n(u'(c))_e$ for $c\in C$ and $e\in E^{(n)}$.

\emph{Balls at every level.} Put $m_\mu:=\min_{\sigma\in\Sinf}|\sigma\mu|>1$ and $M_\mu:=\max_{\sigma\in\Sinf}|\sigma\mu|$. For $n\in\N$ and $s\in\R$, set
\[
 \rho_E:=\sup_{j\ge0}\ \max_{e\in E^{(j)}}\norm{e/\mu^j},\qquad
 \mathbb B_n(s):=\{x\in R:\norm{x/\mu^n}\le s\},\qquad
 \mathbb B(s):=\mathbb B_0(s).
\]
In particular, $\mathbb B_n(s)=\varnothing$ if $s<0$. Since $e=\sum_{i<n}\mu^ie_i$ with $e_i\in E$, we have, place by place,
\[
 |\sigma(e/\mu^n)|\le\sum_{i<n}m_\mu^{\,i-n}\max_{e'\in E}\norm{e'}.
\]
Thus $\rho_E\le\max_{e\in E}\norm e/(m_\mu-1)$, and each block $\mu^nc+E^{(n)}$ is contained in $\mu^nc+\mathbb B_n(\rho_E)$. Recall that $c_R$ denotes the covering radius of $R$ (\S\ref{ssec:places}).

\begin{lemma}[Balls and blocks]\label{lem:levelballs}
For $n,h\in\N$ and $s\ge0$,
\[
 \mathbb B_n(s-\rho_E)\subseteq\mu^n\mathbb B(s)+E^{(n)},\qquad
 \mathbb B_{n+h}(s)\subseteq\mathbb B_n(M_\mu^{\,h}s).
\]
\end{lemma}
\begin{proof}
If $x\in\mathbb B_n(s-\rho_E)$, write $x=\mu^nc+e$, with $e\in E^{(n)}$. Then $c=x/\mu^n-e/\mu^n$ satisfies $\norm c\le s$. For the second inclusion, use the identity $|\sigma(x/\mu^n)|=|\sigma\mu|^h|\sigma(x/\mu^{n+h})|$ at every archimedean place.
\end{proof}

Two patterns occurring in the same point agree wherever they overlap, and substitution preserves this property. Patterns $P,P'$ are \emph{compatible at} $\Delta\in R$ if $P(x)=P'(x-\Delta)$ whenever $x\in\supp P$ and $x-\Delta\in\supp P'$.

\begin{lemma}[Compatibility]\label{lem:consistency}
\textup{(a)} If $P$ and $P'$ occur in the same point $u\in\bA^R$ at $q$ and $q'$ respectively, they are compatible at $q'-q$. \textup{(b)} If $P,P'$ are compatible at $\Delta$, then $\theta^h(P),\theta^h(P')$ are compatible at $\mu^h\Delta$ for every $h\in\N$.
\end{lemma}
\begin{proof}
(a) For $x\in\supp P$ and $x-(q'-q)\in\supp P'$,
\[
 P(x)=u(q+x)=u(q'+x-(q'-q))=P'(x-(q'-q)).
\]
(b) Let $x\in\supp\theta^h(P)$ with $x-\mu^h\Delta\in\supp\theta^h(P')$, and write $x=\mu^hc+e$, with $e\in E^{(h)}$. By uniqueness of this decomposition, $c\in\supp P$ and $c-\Delta\in\supp P'$. Hence $P(c)=P'(c-\Delta)$ and
\[
 \theta^h(P)(x)=\theta^h(P(c))_e
 =\theta^h(P'(c-\Delta))_e=\theta^h(P')(x-\mu^h\Delta).
\]
\end{proof}

By minimality (Lemma~\ref{lem:minimal}), each legal pattern occurs with bounded gaps. The proof of Proposition~\ref{prop:recog} requires a version of this property that is uniform across scales.

\begin{lemma}[Recurrence at every scale]\label{lem:allscales}
Let $\theta$ be primitive. For each real number $s\ge0$, there is an $N=N(s)\in\N$ with the following property. Let $n\in\N$ and suppose that $y\in\bA^R$ has decompositions
\[
 y=S^q\theta^n(u)=S^{q'}\theta^{n+N}(u'),
 \qquad q,q'\in R,\quad u\in X_\theta,\quad u'\in\bA^R.
\]
Then for every finite $F\subseteq\mathbb B_n(s)$ and every $x,t\in R$, there is an $x'\in R$ such that
\[
 x'+F\subseteq t+\mathbb B_{n+N}(\rho_E+c_R),\qquad
 y(x'+w)=y(x+w)\quad(w\in F).
\]
\end{lemma}
\begin{proof}
\emph{Choice of $N$.} There are only finitely many legal patterns whose support is contained in $\mathbb B(s+\rho_E+c_R)$. Each occurs in a block $\theta^{n'}(g')$. Let $N(s)$ be the sum of $p$ (from the definition of primitivity) and the largest of the finitely many levels $n'$ chosen in this way. Every one of these patterns then occurs in $\theta^N(g)$ for every letter $g$, by the proof of Lemma~\ref{lem:minimal}.

\emph{The pattern of $y$ on $x+F$ occurs in every level-$(n+N)$ block.} Put $P:=y|_{x+F}$, $x^*:=(x+q)/\mu^n\in K_{\R}$ and
\[
 C:=\{c\in R:\norm{c-x^*}\le s+\rho_E\}.
\]
For $w\in F$, write $x+w+q=\mu^nc+e$ with $e\in E^{(n)}$. Then $\norm{c-x^*}=\norm{w/\mu^n-e/\mu^n}\le s+\rho_E$, so $c\in C$. Thus $x+F\subseteq-q+\mu^nC+E^{(n)}$, and $P$ occurs in $\theta^n(u|_C)$ at $q$, by (B4).

Choose $c_0\in R$ such that $\norm{c_0-x^*}\le c_R$, and define a pattern $W$ by $W(c):=u(c+c_0)$ for $c\in C-c_0$. Its support is contained in $\mathbb B(s+\rho_E+c_R)$, and it is legal because it occurs at $c_0$ in $u\in X_\theta$. Hence $W$ occurs in $\theta^N(g)$ for every letter $g$. By (B1) and (B2), $\theta^n(W)$ occurs in every $\theta^{n+N}(g)$. Since $\theta^n(u|_C)$ is the translate of $\theta^n(W)$ by $\mu^nc_0$, the pattern $P$ also occurs in every $\theta^{n+N}(g)$.

\emph{Every level-$(n+N)$ ball of radius $\rho_E+c_R$ contains a block.} By (B4), for each $c\in R$, the configuration $y$ carries the block $\theta^{n+N}(u'(c))$ on
\[
 -q'+\mu^{n+N}c+E^{(n+N)}.
\]
Choose $c\in R$ such that $\norm{c-(t+q')/\mu^{n+N}}\le c_R$. For every $e\in E^{(n+N)}$,
\[
\begin{split}
 \norm{\bigl(-q'+\mu^{n+N}c+e-t\bigr)/\mu^{n+N}}
 &\le\norm{c-(t+q')/\mu^{n+N}}+\norm{e/\mu^{n+N}}\\
 &\le c_R+\rho_E.
\end{split}
\]
The block is therefore contained in $t+\mathbb B_{n+N}(\rho_E+c_R)$. By the preceding part of the proof, $P$ occurs in $\theta^{n+N}(u'(c))$, say at $a\in R$. This means that $a+x+F\subseteq E^{(n+N)}$ and $\theta^{n+N}(u'(c))_{a+x+w}=y(x+w)$ for $w\in F$. Set $x':=-q'+\mu^{n+N}c+a+x$. Then
\[
 x'+F\subseteq-q'+\mu^{n+N}c+E^{(n+N)}\subseteq t+\mathbb B_{n+N}(\rho_E+c_R),
\]
and by (B4) we have $y(x'+w)=\theta^{n+N}(u'(c))_{a+x+w}=y(x+w)$ for every $w\in F$.
\end{proof}

\begin{proposition}[First step of recognisability]\label{prop:recog}
Let $\theta$ be primitive with expansion $\mu$ and shape $E$, and suppose that $X_\theta$ is aperiodic. Let $z\in X_\theta$ satisfy $\theta(z)=z$. There is an $r>0$ such that, for $i,j\in R$, if $z(i+x)=z(j+x)$ for every $x\in\mathbb B(r)$ and $i\in\mu R$, then $j\in\mu R$.
\end{proposition}
\emph{Idea of the proof.} If the property fails, a limit point $y$ of translates of $z$ has, at every level $n$, two decompositions into level-$n$ blocks. In the first decomposition the block origins belong to the zero class modulo $\mu R$; in the second they belong to another class. At two levels $n<m$, the local patterns of the two desubstituted configurations agree. Applying $\theta^{m-n}$ to the patterns at level $n$ and comparing them with those at level $m$ gives the same pair of patterns in two relative positions. The difference $v$ of these positions is not in $\mu R$ and is a period of $y$ on a large ball. Lemma~\ref{lem:allscales} allows us to find a copy of every pattern of the required size inside this ball. It follows that $v$ is a period of $y$, contrary to aperiodicity.

\begin{proof}
Suppose that the property fails.

\emph{Step 1: Two decompositions at every level.} For each $k\ge1$ there are $i_k\in\mu R$ and $j_k\in R\setminus\mu R$ with $z(i_k+x)=z(j_k+x)$ for $x\in\mathbb B(k)$. Passing to a subsequence, we may suppose that the $j_k$ belong to a single residue class $j_*+\mu R$, with $j_*\in R\setminus\mu R$, and that $S^{i_k}z\to y$. Then $S^{j_k}z\to y$ as well, and $y\in X_\theta$.

Fix $n\ge1$ and write $i_k=\mu^nb_k+e_k$ and $j_k=\mu^nb'_k+e'_k$, with $e_k,e'_k\in E^{(n)}$. Since $z=\theta^n(z)$, (B3) gives
\[
 S^{i_k}z=S^{e_k}\theta^n(S^{b_k}z),\qquad
 S^{j_k}z=S^{e'_k}\theta^n(S^{b'_k}z).
\]
Along a further subsequence, $e_k=a_n$ and $e'_k=a'_n$ are constant, and $S^{b_k}z\to u_n$ and $S^{b'_k}z\to u'_n$, with $u_n,u'_n\in X_\theta$. By continuity,
\begin{equation}\label{eq:twodec}
 y=S^{a_n}\theta^n(u_n)=S^{a'_n}\theta^n(u'_n),\qquad a_n,a'_n\in E^{(n)}.
\end{equation}
Since $\mu^nR\subseteq\mu R$, we have $a_n\equiv i_k\equiv0\pmod{\mu R}$ and $a'_n\equiv j_k\equiv j_*\not\equiv0\pmod{\mu R}$. This holds for every $n\ge1$. No diagonal subsequence is needed: the rest of the proof uses only the identities \eqref{eq:twodec}, each of which involves a single level and the common limit $y$.

\emph{Step 2: Two levels with the same local data.} Let $N:=N(4\rho_E)$ be given by Lemma~\ref{lem:allscales}, and put $\bar r:=5\rho_E+M_\mu^{\,N}(\rho_E+c_R)$. The pairs of restrictions
\[
 U_n:=u_n|_{\mathbb B(\bar r)},\qquad U'_n:=u'_n|_{\mathbb B(\bar r)}
\]
can take only finitely many values. Thus there are $1\le n<m$ with $U_n=U_m=:U$ and $U'_n=U'_m=:U'$. Note that the radius $\bar r$ is fixed before $n$ and $m$ are chosen.

\emph{Step 3: A pair of patterns in two relative positions.} By \eqref{eq:twodec} and (B4), the patterns $\theta^n(U)$ and $\theta^n(U')$ occur in $y$ at $-a_n$ and $-a'_n$, respectively. By Lemma~\ref{lem:consistency}, they are compatible at $a_n-a'_n$. Applying $\theta^{m-n}$, we find that $\theta^m(U)$ and $\theta^m(U')$ are compatible at
\[
 \Delta_1:=\mu^{m-n}(a_n-a'_n).
\]
On the other hand, $\theta^m(U)$ and $\theta^m(U')$ occur in $y$ at $-a_m$ and $-a'_m$, respectively, so they are compatible at $\Delta_2:=a_m-a'_m$. Since $m-n\ge1$, we have $\Delta_1\in\mu R$, whereas $\Delta_2\equiv-j_*\not\equiv0\pmod{\mu R}$. Thus $v:=\Delta_2-\Delta_1\notin\mu R$, in particular $v\ne0$. Moreover,
\[
 \norm{\Delta_1/\mu^m}=\norm{(a_n-a'_n)/\mu^n}\le2\rho_E,
 \qquad \norm{\Delta_2/\mu^m}\le2\rho_E,
\]
so $v\in\mathbb B_m(4\rho_E)$.

\emph{Step 4: A local period.} By Lemma~\ref{lem:levelballs}, the ball $\mathbb B_m(\bar r-\rho_E)$ is contained in $\mu^m\mathbb B(\bar r)+E^{(m)}$, the common support of $\theta^m(U)$ and $\theta^m(U')$. Let $w\in\mathbb B_m(\bar r-5\rho_E)$ and set $\xi:=w+\Delta_2$. Then
\[
 \xi\in\mathbb B_m(\bar r-3\rho_E),\qquad
 \xi-\Delta_2=w,\qquad
 \xi-\Delta_1=w+v\in\mathbb B_m(\bar r-\rho_E).
\]
All three positions belong to the common support. Compatibility at $\Delta_1$ and at $\Delta_2$ therefore gives
\[
 \theta^m(U')(w+v)=\theta^m(U)(\xi)=\theta^m(U')(w).
\]
Since $\theta^m(U')$ occurs in $y$ at $-a'_m$, we obtain
\begin{equation}\label{eq:localperiod}
 y(-a'_m+w+v)=y(-a'_m+w)\qquad(w\in\mathbb B_m(\bar r-5\rho_E)).
\end{equation}

\emph{Step 5: The period propagates.} Let $x\in R$ and set $F:=\{0,v\}\subseteq\mathbb B_m(4\rho_E)$. Apply Lemma~\ref{lem:allscales} with $s=4\rho_E$, at level $m$, to the decompositions
\[
 y=S^{a_m}\theta^m(u_m)=S^{a_{m+N}}\theta^{m+N}(u_{m+N})
\]
given by \eqref{eq:twodec}, with $u_m\in X_\theta$, and take $t:=-a'_m$. There is an $x'\in R$ such that $x'+F\subseteq-a'_m+\mathbb B_{m+N}(\rho_E+c_R)$, $y(x')=y(x)$ and $y(x'+v)=y(x+v)$. By Lemma~\ref{lem:levelballs},
\[
 \mathbb B_{m+N}(\rho_E+c_R)\subseteq\mathbb B_m(M_\mu^{\,N}(\rho_E+c_R))
 =\mathbb B_m(\bar r-5\rho_E).
\]
As $0\in F$, we have $x'=-a'_m+w$ for some $w\in\mathbb B_m(\bar r-5\rho_E)$. By \eqref{eq:localperiod}, $y(x'+v)=y(x')$. Hence $y(x+v)=y(x)$ for every $x\in R$. The non-zero element $v$ is a period of $y\in X_\theta$ and thus, by minimality, of every point of $X_\theta$. This contradicts aperiodicity.
\end{proof}

The proof gives no explicit value of $r$, because Step~1 uses compactness. Only the existence of $r$ is needed here. Applying the proposition to $i-e$ and $j-e$, where $e\in E$ represents $i$ modulo $\mu R$, also shows that the residue class of a position modulo $\mu R$ is determined by the values of $z$ on the ball of radius $r+\max_{e\in E}\norm e$ about that position.

\subsection{The odometer}\label{ssec:odometer}
From Proposition~\ref{prop:recog} we deduce that an aperiodic substitution subshift factors onto the finite groups $R/\nu^nR$, for a positive power $\nu$ of the base; these factors are used in \S\ref{ssec:rankzero}.

\begin{corollary}[Odometer factor]\label{cor:odometer}
Let $\theta$ be primitive with expansion $\mu$ and shape $E$, where $(E,\mu)$ is an adapted window for $\alpha$. Suppose that $X_\theta$ is aperiodic. There is a power $\nu=\alpha^t$ ($t\ge1$) such that, for every $n\ge1$, there is a continuous map $\psi_n:X_\theta\to R/\nu^nR$ satisfying
\[
 \psi_n(S^qu)=\psi_n(u)+(q\bmod\nu^nR)\qquad(q\in R,\ u\in X_\theta).
\]
\end{corollary}
\begin{proof}
\emph{Powers of $\theta$.} By Lemma~\ref{lem:germ}, choose $z\in X_\theta$ with $\theta^M(z)=z$ for some $M\ge1$, and put $\nu:=\mu^M$, a positive power of $\alpha$. Fix $n\ge1$. The substitution $\theta^{Mn}$ has expansion $\nu^n$ and shape $E^{(Mn)}$, and $(E^{(Mn)},\nu^n)$ is again an adapted window. It is primitive, since every positive power of a primitive substitution is primitive (see the proof of Lemma~\ref{lem:minimal}).

Its legal patterns are legal for $\theta$, so $X_{\theta^{Mn}}\subseteq X_\theta$. Also $z\in X_{\theta^{Mn}}$, since the patterns of $z$ occur in the blocks $\theta^{Mnk}(g)$, where $g$ is the letter given by Lemma~\ref{lem:germ}. As $X_\theta$ is minimal by Lemma~\ref{lem:minimal} and $X_{\theta^{Mn}}$ is closed and invariant, we have $X_{\theta^{Mn}}=X_\theta$. In particular this subshift is aperiodic, and $z$ is fixed by $\theta^{Mn}$. Proposition~\ref{prop:recog}, applied to $\theta^{Mn}$, gives a radius $r_n>0$ such that if $z(i+x)=z(j+x)$ for $\norm x\le r_n$ ($x\in R$) and $i\in\nu^nR$, then $j\in\nu^nR$.

\emph{A finite partition.} Let $C_n$ be the closure of $\{S^{\nu^nw}z:w\in R\}$. For $g\in R$, the set $S^gC_n$ is the closure of $\{S^{g+\nu^nw}z:w\in R\}$ and depends only on the class of $g$ modulo $\nu^nR$. Since the orbit of $z$ is dense in the minimal subshift $X_\theta$, and $R/\nu^nR$ is finite, the union of these closed sets is $X_\theta$.

If $u\in C_n\cap S^gC_n$, there exist sequences $w_m,w'_m\in R$ such that
\[
 u=\lim_m S^{\nu^nw_m}z=\lim_m S^{g+\nu^nw'_m}z.
\]
For sufficiently large $m$ the two translates agree with $u$ on the ball of radius $r_n$. The choice of $r_n$ implies that $g+\nu^nw'_m\in\nu^nR$, hence $g\in\nu^nR$. Translating this argument shows that $S^gC_n$ and $S^{g'}C_n$ are disjoint whenever $g\not\equiv g'\pmod{\nu^nR}$. The sets $S^gC_n$ thus form a finite partition of $X_\theta$ into closed, and hence also open, sets.

\emph{The factor.} Define $\psi_n(u):=g\bmod\nu^nR$ if $u\in S^gC_n$. This is a locally constant, surjective map. With the discrete topology on $R/\nu^nR$, it is continuous, and $S^q(S^gC_n)=S^{q+g}C_n$ gives $\psi_n(S^qu)=\psi_n(u)+(q\bmod\nu^nR)$.
\end{proof}

The partitions are compatible as $n$ varies: $C_{n+1}\subseteq C_n$, so $\psi_{n+1}$ reduces to $\psi_n$ modulo $\nu^nR$. They thus define a factor map onto the inverse limit $\varprojlim_n R/\nu^nR$. Indeed, every compatible sequence of residue classes determines a nested family of non-empty closed sets, whose intersection is non-empty by compactness. In what follows we use only the finite factors $\psi_n$.

\subsection{No intermediate rank}\label{ssec:ranks}
Under (H), the period group of a minimal subsystem of the orbit closure of a doubly automatic function cannot have intermediate rank: its rational span is stable under multiplication by powers of both bases, hence is a $K$-vector subspace of $K$.

\begin{lemma}\label{lem:periods}
Let $g:R\to\Omega$ have only finitely many kernels $y\mapsto g(c+\gamma^ny)$ ($n\in\N$, $c\in E^{(n)}$), where $\gamma\in R$ is expansive and $E$ is a window set for $\gamma$. Set $W:=\Per(g)$. There exist $n<m$ in $\N$ such that $W':=W\cap\gamma^nR$ has finite index in $W$ and $\gamma^{m-n}W'\subseteq W'$.
\end{lemma}
\begin{proof}
Since $E^{(n)}$ meets every residue class modulo $\gamma^nR$,
\[
 \bigcap_{c\in E^{(n)}}\Per\bigl(y\mapsto g(c+\gamma^ny)\bigr)
 =\{w\in R:\gamma^nw\in W\}.
\]
The intersection on the left can take only finitely many values, because the family of kernels is finite. Thus it takes the same value at two levels $n<m$. The subgroup $W'=W\cap\gamma^nR$ has finite index in $W$, since $\gamma^nR$ has finite index in $R$. If $u=\gamma^nw\in W'$, equality of the two intersections implies that $\gamma^{m-n}u=\gamma^mw\in W\cap\gamma^nR=W'$.
\end{proof}

\begin{theorem}[No intermediate rank]\label{thm:ranks}
Under (H), if $f:R\to\Omega$ is doubly automatic and $Y\subseteq\cX_f$ is minimal, then $\Lambda_Y$ is either zero or of rank $d$.
\end{theorem}
\begin{proof}
Applying Lemma~\ref{lem:generator} first to $\alpha$ and then to $\beta$, we obtain points $y_\alpha,y_\beta\in Y$, automatic in the respective bases, both with period group $\Lambda_Y$. By Lemma~\ref{lem:periods}, the rational span $V:=\operatorname{span}_{\Q}\Lambda_Y\subseteq K$ satisfies
\[
 \alpha^rV\subseteq V,\qquad \beta^sV\subseteq V
\]
for some $r,s\ge1$: the finite-index subgroups in that lemma have the same rational span as $\Lambda_Y$. By (H), the algebra $\Q[\alpha^r,\beta^s]$ equals $K$ (Lemma~\ref{lem:fieldpowers}). Hence $V$ is a $K$-vector subspace of $K$. It is either $0$ or $K$ itself.
\end{proof}

\subsection{Excluding rank zero when the supports differ}\label{ssec:rankzero}
Recall that the support of a base is the set of prime ideals containing it. Suppose that a prime ideal contains $\alpha$ but not $\beta$. We exclude rank zero by a counting argument. The substitution subshift lifting an aperiodic minimal subsystem $Y$ admits equivariant maps onto $R/\nu^nR$, where $\nu$ is a positive power of $\alpha$. On the other hand, $\beta$-automaticity gives a uniform bound on the size of the image of $\beta^kR$ in such a factor, for sufficiently large $k$.

\begin{proposition}[Bound on finite factors]\label{prop:factorbound}
Let $f:R\to\Omega$ be $\beta$-automatic. There is an integer $N^*\ge1$, depending only on $f$, $\beta$ and a fixed window set for $\beta$, with the following property. For every subsystem $Z\subseteq\cX_f$, every finite abelian group $F$, every additive homomorphism $\chi:R\to F$ and every continuous map $\psi:Z\to F$ satisfying
\[
 \psi(S^qu)=\psi(u)+\chi(q)\qquad(q\in R,\ u\in Z),
\]
we have $|\chi(\beta^kR)|\le N^*$ for all sufficiently large $k\in\N$ (with a threshold that may depend on $Z$, $F$, $\chi$ and $\psi$).
\end{proposition}
\begin{proof}
Fix a window set $D_\beta$ for $\beta$, with size constant $C_\beta$ (Lemma~\ref{lem:fixedlength}), and the kernels $\cK_\beta(f)$ for $D_\beta$; put
\[
 \cW:=\{w\in R:\norm w\le2C_\beta+1\},\qquad
 N^*:=|\Omega|^{|\cK_\beta(f)||\cW|}.
\]
Since $\psi$ is continuous on the compact set $Z$ with finite discrete codomain, $\psi(u)$ depends only on the restriction of $u$ to a ball $B_\psi:=\{x\in R:\norm x\le r_\psi\}$. Let $k\in\N$ satisfy $\rho_\beta^kr_\psi\le1$, where $\rho_\beta=\max_{\sigma\in\Sinf}|\sigma\beta|^{-1}<1$ (Definition~\ref{def:expansive}).

Take $u\in Z$ and $q_i\in R$ with $S^{q_i}f\to u$, and write $q_i=\beta^kz_i+t_i$ with $z_i\in R$ and $t_i\in D_\beta^{(k)}$. Passing to a subsequence, we may assume that $t_i=t$ for all $i$. For $x\in B_\psi$, write
\[
 t+x=\beta^kw_x+c_x,\qquad
 w_x\in R,\quad c_x\in D_\beta^{(k)}.
\]
At every archimedean place,
\[
 |\sigma w_x|\le 2C_\beta+\rho_\beta^kr_\psi
 \le 2C_\beta+1,
\]
so $w_x\in\cW$. For $v\in R$,
\[
 (S^{\beta^kv+q_i}f)(x)=f\bigl(c_x+\beta^k(v+z_i+w_x)\bigr),
\]
where $y\mapsto f(c_x+\beta^ky)$ belongs to $\cK_\beta(f)$. Hence the restriction of $S^{\beta^kv+q_i}f$ to $B_\psi$ is determined by the table
\[
 \bigl(G(v+z_i+w)\bigr)_{G\in\cK_\beta(f),\ w\in\cW},
\]
which takes at most $N^*$ values.

If $|\chi(\beta^kR)|>N^*$, choose $v_0,\dots,v_{N^*}\in R$ with distinct images $\chi(\beta^kv_s)$. For each $i$, two of their $N^*+1$ tables coincide, and some fixed pair $s\ne s'$ does so for infinitely many $i$. Passing to the limit along these indices, we find that $S^{\beta^kv_s}u$ and $S^{\beta^kv_{s'}}u$, which lie in $Z$, agree on $B_\psi$. Hence $\psi$ takes the same value at both, and equivariance gives $\chi(\beta^kv_s)=\chi(\beta^kv_{s'})$, a contradiction.
\end{proof}

\begin{theorem}[Finite orbits]\label{thm:finiteorbits}
Assume (H), and suppose that a prime ideal $Q$ of $R$ contains $\alpha$ but not $\beta$. If $f:R\to\Omega$ is doubly automatic, then $\cX_f$ has finitely many minimal subsystems, each of which is a finite orbit.
\end{theorem}
\begin{proof}
Finiteness follows from Lemma~\ref{lem:lift}. In view of Theorem~\ref{thm:ranks} and Lemma~\ref{lem:rankfinite}, it remains to exclude an aperiodic minimal subsystem $Y\subseteq\cX_f$.

Take an adapted window $(D,\lambda)$ for $\alpha$ (Proposition~\ref{prop:window}). By Lemma~\ref{lem:lift}, we have $Y=\pr(Z)$, with $Z=X_\theta\subseteq\cX_\Theta$ minimal, where $\theta$ is primitive of expansion $\lambda^h$ and shape $D^{(h)}$. Moreover, $Z$ is aperiodic. Let $\nu=\alpha^t$ and $\psi_n:Z\to R/\nu^nR$ be given by Corollary~\ref{cor:odometer}.

The configuration $\Theta$ is $\beta$-automatic: its coordinates are the $\lambda$-kernels of $f$, which are $\beta$-automatic by Proposition~\ref{prop:affine}, and Corollary~\ref{cor:closure} applies to their tuple. Let $N^*$ be the constant of Proposition~\ref{prop:factorbound} for $\Theta$; it depends only on $\Theta$, $\beta$ and a window set. Choose $n\ge1$ with $2^n>N^*$ and denote reduction modulo $\nu^nR$ by $\chi_n:R\to R/\nu^nR$.

Since $\nu\in Q$, the natural map $R/\nu^nR\to R/Q^n$ is surjective. Since $\beta\notin Q$, the element $\beta$ is invertible modulo $Q^n$. Thus $\beta^kR$ maps onto $R/Q^n$, and for every $k\in\N$,
\[
 |\chi_n(\beta^kR)|\ge |R/Q^n|=|R/Q|^n\ge2^n>N^*.
\]
This contradicts Proposition~\ref{prop:factorbound}, applied to $\Theta$, the subsystem $Z\subseteq\cX_\Theta$, the group $F=R/\nu^nR$ and the map $\psi_n$, with $k$ beyond the threshold (which may depend on $n$).
\end{proof}

The order of the quantifiers in Proposition~\ref{prop:factorbound} is essential: $N^*$ is independent of $F$, $\chi$ and $\psi$, whereas the threshold for $k$ may depend on them. The lower bound above holds for every $k$, so this dependence causes no difficulty.

\begin{example}\label{ex:Qsupport}
In $R=\Z[\sqrt2]$, take $\alpha=7$ and $\beta=3+\sqrt2$. With $P=(3+\sqrt2)R$ and $\bar P=(3-\sqrt2)R$ as in Example~\ref{ex:parity}, the prime ideal $\bar P$ contains $7$ but not $\beta$. Under the Chinese remainder identification
\[
 R/7^nR\simeq R/P^n\times R/\bar P^n,
\]
the image of $\beta^kR$ is $P^{\min(k,n)}/P^n\times R/\bar P^n$, of order at least $7^n$. The second component is always filled, because $\beta$ is invertible modulo $\bar P^n$. For the pair $(14,14u)$ of the Introduction, no such prime ideal exists: every prime ideal containing $14$ also contains $\beta$ (the case of Section~\ref{sec:equal}).
\end{example}

\subsection{Reduction to the rigid class}\label{ssec:reduction}
When the supports differ, a doubly automatic set that is not eventually periodic gives rise to a member of the rigid class defined below by taking a difference along a non-zero common period of the minimal subsystems of its orbit closure.

We identify a set $X\subseteq R$ with its characteristic function, taking values in $\F_2$, and write $\mathbf 0$ for the zero function. For $v\in R$, put
\[
 \partial_vX:=\{x\in R:\mathbf 1_X(x)\ne\mathbf 1_X(x+v)\}.
\]
We use $\cX_X$ and $\Perf(X)$ for the orbit closure and the group of periods at infinity of the characteristic function of $X$.

\begin{definition}[Rigid class]\label{def:I}
The class $\cI$ consists of the doubly automatic subsets $X\subseteq R$ satisfying the following two conditions:
\begin{enumerate}[label=\textup{(\roman*)},nosep]
\item $\Perf(X)=0$;
\item $\{\mathbf 0\}$ is the only minimal subsystem of $\cX_X$.
\end{enumerate}
\end{definition}
Every member of $\cI$ is infinite and is not eventually periodic. A finite set has $\Perf(X)=R$, and eventual periodicity implies a group of periods at infinity of rank $d$.

\begin{lemma}\label{lem:bdstab}
Assume (H). If $X\subseteq R$ is doubly automatic and has no period at infinity, then for every $v\in R\setminus\{0\}$, the set $\partial_vX$ is doubly automatic, infinite and has no period at infinity.
\end{lemma}
\begin{proof}
The set $\partial_vX$ is doubly automatic by Corollary~\ref{cor:closure}. It is infinite, since otherwise $v\in\Perf(X)$. If $\partial_vX$ had a period at infinity, Theorem~\ref{thm:onefar} would imply that $\Perf(\partial_vX)$ has full rank in $R$. Hence $mv\in\Perf(\partial_vX)$ for some $m\ge1$. Outside a finite set, we would then have, in $\F_2$,
\[
\begin{split}
 \mathbf 1_X(x+2mv)+\mathbf 1_X(x)
 &=\sum_{i<m}\bigl(\mathbf 1_{\partial_vX}(x+iv)
       +\mathbf 1_{\partial_vX}(x+iv+mv)\bigr)\\
 &=0.
\end{split}
\]
The first equality is obtained by telescoping the differences in the direction $v$. The second holds outside a finite union of translates of the exceptional set for the period $mv$. It follows that $2mv\in\Perf(X)$, a contradiction.
\end{proof}

\begin{lemma}\label{lem:ifacered}
Suppose that $\cX_X$ has finitely many minimal subsystems $Y_1,\ldots,Y_s$, all finite orbits, and put $L:=\bigcap_{i=1}^s\Lambda_{Y_i}$. Then $L$ has finite index in $R$, and for every $v\in L$, the singleton $\{\mathbf 0\}$ is the only minimal subsystem of $\cX_{\partial_vX}$.
\end{lemma}
\begin{proof}
The index assertion follows from Lemma~\ref{lem:rankfinite}. The map $\Delta_v:\F_2^R\to\F_2^R$ given by $\Delta_v(u)=u+S^vu$ is continuous, commutes with translations and maps $\cX_X$ onto $\cX_{\partial_vX}$. Let $Z$ be minimal in $\cX_{\partial_vX}$. The set $\Delta_v^{-1}(Z)\cap\cX_X$ is non-empty, closed and invariant, so it contains one of the minimal subsystems $Y_i$. As $v\in\Lambda_{Y_i}$, we have $\Delta_v(Y_i)=\{\mathbf 0\}$. By minimality of $Z$, this gives $Z=\{\mathbf 0\}$.
\end{proof}

\begin{theorem}[Reduction]\label{thm:reduction}
Let $d\ge2$ and let the pair $(\alpha,\beta)$ satisfy (H). Suppose that a prime ideal of $R$ contains $\alpha$ but not $\beta$. If $X\subseteq R$ is doubly automatic but not eventually periodic, there is a finite-index subgroup $L\subseteq R$ such that $\partial_vX\in\cI$ for every $v\in L\setminus\{0\}$. Thus $\cI=\varnothing$ if and only if every doubly automatic subset of $R$ is eventually periodic.
\end{theorem}
\begin{proof}
By Theorem~\ref{thm:onefar}, $X$ has no period at infinity. Apply Theorem~\ref{thm:finiteorbits}, followed by Lemmas~\ref{lem:bdstab} and \ref{lem:ifacered}. The equivalence follows because no member of $\cI$ is eventually periodic and every doubly automatic counterexample gives a member of $\cI$.
\end{proof}

\subsection{Structure of the rigid class}\label{ssec:structure}
This subsection collects the properties of the members of $\cI$ used in Section~\ref{sec:different}: stability under infinite sections, a bound on their number of points in each tile, and the emptiness of the interior of their $\alpha$-adic closure.

For a subset $X\subseteq R$ and $t,\gamma\in R$, with $\gamma\ne0$, set $X_{\gamma,t}:=\{y\in R:t+\gamma y\in X\}$.

\begin{lemma}[Holes]\label{lem:holes}
The singleton $\{\mathbf 0\}$ is the only minimal subsystem of $\cX_X$ if and only if every $u\in\cX_X$ vanishes on a translate of every finite subset of $R$. This property passes to the sections $X_{\gamma,t}$.
\end{lemma}
\begin{proof}
For the equivalence, observe that the orbit closure of every $u\in\cX_X$ contains a minimal subsystem. The condition that $\mathbf 0$ belongs to this orbit closure is exactly that $u$ vanishes on a translate of each finite subset of $R$. If the latter condition holds for every $u\in\cX_X$, every minimal subsystem contains $\mathbf 0$ and must therefore equal $\{\mathbf 0\}$.

For the assertion about sections, note that
\[
 (S^q\mathbf 1_{X_{\gamma,t}})(y)
 =(S^{t+\gamma q}\mathbf 1_X)(\gamma y)\qquad(q,y\in R).
\]
If $z\in\cX_{X_{\gamma,t}}$, write $z=\lim_n S^{q_n}\mathbf 1_{X_{\gamma,t}}$. By passing to a subsequence such that $S^{t+\gamma q_n}\mathbf 1_X$ converges, we obtain $w\in\cX_X$ with $z(y)=w(\gamma y)$ for all $y\in R$.

Let $F'\subseteq R$ be finite and let $T\subseteq R$ be a complete set of representatives of $R/\gamma R$. To account for the residue class of the translation provided by the hypothesis, put
\[
 F:=\{\gamma\xi-s:\xi\in F',\ s\in T\}.
\]
By hypothesis, $w$ vanishes on $q+F$ for some $q\in R$. Write $q=\gamma q'+s_0$ with $q'\in R$ and $s_0\in T$. For $\xi\in F'$, we have $\gamma(q'+\xi)=q+\gamma\xi-s_0\in q+F$, hence $z(q'+\xi)=w(\gamma(q'+\xi))=0$. Thus $z$ vanishes on $q'+F'$.
\end{proof}

\begin{proposition}[Stability under sections]\label{prop:stability}
Assume $d\ge2$ and (H), and let $X\in\cI$. Every section $X_{\gamma,t}$ ($t,\gamma\in R$, $\gamma\ne0$) is either finite or belongs to $\cI$. For each $\gamma\in R\setminus\{0\}$, at least one such section is infinite. In particular every infinite $\alpha$- or $\beta$-kernel of $X$, viewed as a set, belongs to $\cI$.
\end{proposition}
\begin{proof}
The set $X':=X_{\gamma,t}$ is doubly automatic (Proposition~\ref{prop:affine}), and every point of its orbit closure vanishes on a translate of every finite subset of $R$, that is, $\{\mathbf 0\}$ is the only minimal subsystem of $\cX_{X'}$ (Lemma~\ref{lem:holes}). If $X'$ has no period at infinity, it therefore belongs to $\cI$.

Otherwise $X'$ is eventually periodic (Theorem~\ref{thm:onefar}): its characteristic function agrees outside a finite set with a function $\tilde g:R\to\F_2$ whose period group $\Lambda$ has full rank. Translating by elements of $\Lambda$ whose norms tend to infinity shows that $\tilde g\in\cX_{X'}$. Hence $\tilde g$ vanishes on a translate of a complete set of representatives of $R/\Lambda$; as such a translate meets every class, $\tilde g=0$ and $X'$ is finite.

Finally, for each $\gamma\ne0$, the set $X$ is the disjoint union of the sets $s+\gamma X_{\gamma,s}$, where $s$ runs through a complete set of representatives of $R/\gamma R$. Since $X$ is infinite and the union is finite, at least one section is infinite.
\end{proof}

For the rest of this subsection, fix $X\in\cI$ and an adapted window $(D,\lambda)$ for $\alpha$. Write $q_D:=|D|$, and let $\zeta,\Theta,\pr$ be the objects of \S\ref{ssec:kernelsubst} for $f=\mathbf 1_X$ and $\Omega=\F_2$. Put $\bA_\Theta:=\{\Theta(y):y\in R\}$ and $\tau_\Theta:=|\bA_\Theta|$. The level-$n$ tile indexed by $y\in R$ carries the block $\pr(\zeta^n(\Theta(y)))$. Let $N_n$ be the maximum number of points of $X$ in a level-$n$ tile.

\begin{lemma}[Sinks carry only zeros]\label{lem:sinkzero}
If $\{\mathbf 0\}$ is the only minimal subsystem of $\cX_X$ and $C\subseteq\bA_\Theta$ is a sink of the letter graph, then $\pr(g)=0$ for every $g\in C$.
\end{lemma}
\begin{proof}
For $g\in C\subseteq\bA_\Theta$, the blocks $\zeta^n(g)$ occur in $\Theta$. Thus every legal pattern for $\theta_i=\zeta^h|_{C_i}$ occurs in $\Theta$, and $X_{\theta_i}\subseteq\cX_\Theta$: each finite restriction of a point in $X_{\theta_i}$ can be realised by a translate of $\Theta$. By Lemma~\ref{lem:minimal}, $X_{\theta_i}$ is minimal, so its image $\pr(X_{\theta_i})$ is minimal and must equal $\{\mathbf 0\}$. Every letter of $C_i$ is legal and occurs in every point of $X_{\theta_i}$ (Lemma~\ref{lem:minimal}). Its projection is therefore zero. The classes $C_i$ partition $C$, giving the result.
\end{proof}

\begin{proposition}[Thinness]\label{prop:thin}
Let $X\in\cI$. For every $n\in\N$, each level-$(n+\tau_\Theta)$ tile contains a level-$n$ tile disjoint from $X$. Moreover,
\[
 N_n\le(q_D/\varrho)^{\tau_\Theta}\varrho^n,
 \qquad \varrho:=(q_D^{\tau_\Theta}-1)^{1/\tau_\Theta}<q_D.
\]
\end{proposition}
\begin{proof}
The set $\bA_\Theta$ is closed under successors, because the letters of $\zeta(\Theta(y))$ are $\Theta(b+\lambda y)$ for $b\in D$. Hence a sink in the letter graph restricted to $\bA_\Theta$ is also a sink of the full graph. Starting from the letter indexing a tile, there is a path of length exactly $\tau_\Theta$ to such a sink: take a shortest path to a sink and extend it within the sink. The corresponding level-$n$ subtile carries only letters of the sink, all of which project to zero by Lemma~\ref{lem:sinkzero}. It follows that
\[
 N_{n+\tau_\Theta}\le(q_D^{\tau_\Theta}-1)N_n.
\]
Since $\lambda$ is expansive, $q_D=[R:\lambda R]\ge2$, and $\varrho>0$. Write $n=j\tau_\Theta+r$, with $0\le r<\tau_\Theta$. Using $N_r\le q_D^r$, the size of a level-$r$ tile, and iterating the inequality gives
\[
 N_n\le\varrho^{j\tau_\Theta}N_r
 \le(q_D/\varrho)^r\varrho^n
 \le(q_D/\varrho)^{\tau_\Theta}\varrho^n.
\]
\end{proof}

Let $R_\alpha:=\varprojlim_n R/\alpha^nR$ be the $\alpha$-adic completion, in which $R$ is dense. Because $\lambda$ is a positive power of $\alpha$, we may equally use the inverse limit over the ideals $\lambda^mR$. The kernel of the projection to $R/\lambda^mR$ is $\lambda^mR_\alpha$: a compatible family in the kernel can be divided by $\lambda^m$ using its coordinates at levels $m+j$, $j\in\N$. Thus the cosets modulo $\lambda^mR_\alpha$ form a basis for the topology. By the Chinese remainder theorem, $R_\alpha$ is the product of the completions of $R$ at the prime ideals containing $\alpha$ (\S\ref{ssec:completion}).

\begin{proposition}[Empty interior]\label{prop:emptyint}
If $X\in\cI$, its closure in $R_\alpha$ has empty interior.
\end{proposition}
\begin{proof}
Suppose that the closure contains a coset $b+\lambda^mR_\alpha$, with $m\in\N$; by density, we may take $b\in R$. The non-zero $\lambda$-kernels of $\mathbf 1_X$ are finitely many, so we may choose a fixed finite ball $V_0\subseteq R$ meeting the support of each of them.

For $n\ge m$, there are exactly $q_D^{n-m}$ elements $c\in D^{(n)}$ with $c\equiv b\pmod{\lambda^mR}$. Each coset $c+\lambda^nR_\alpha$ is a non-empty open subset of $b+\lambda^mR_\alpha$, and therefore meets $X$. Its intersection with $R$ is $c+\lambda^nR$. The associated kernel is consequently non-zero. Choose $y_c\in V_0$ such that $c+\lambda^ny_c\in X$. These points are distinct, since the representatives $c$ belong to different classes modulo $\lambda^nR$. For at least one $y\in V_0$, the level-$n$ tile indexed by $y$ therefore contains at least $q_D^{n-m}/|V_0|$ points of $X$. As $\varrho<q_D$, this contradicts the bound in Proposition~\ref{prop:thin} for large $n$.
\end{proof}
\section{Different supports}\label{sec:different}

Throughout this section, unless otherwise stated, $K$ is a number field of degree $d\ge2$, $R=\cO_K$, and the expansive bases $\alpha,\beta\in R$ satisfy (H).

We assume that the supports of the bases differ, and more precisely that some prime ideal contains $\alpha$ but not $\beta$. We call such a prime an \emph{excess prime}, and a prime containing both bases a \emph{common prime}. By Theorem~\ref{thm:reduction}, it remains to show that $\cI=\varnothing$. Starting from a hypothetical member $X$ of $\cI$, we shall construct an infinite \emph{sparse} set (a set $S\subseteq R$ for which $\#\{x\in S:\norm x\le N\}\le C(\log(2+N))^{h_0}$ for some constants $C>0$, $h_0\in\N$ and every real $N\ge0$), and then show that every doubly automatic sparse set is finite. A \emph{symmetry} of a function $G:R\to\Omega$ is an affine map $\varphi(x)=w+\mu x$, with $w,\mu\in K$, $\mu\notin\{0,1\}$ and $\varphi(R)\subseteq R$, such that $G(\varphi(y))=G(y)$ for every $y\in R$; the identity is exact, not only outside a finite set. We call $\mu$ its \emph{multiplier}, and speak of a \emph{$\beta$-symmetry} if $\mu$ is a positive power of $\beta$; a symmetry of a set is one of its characteristic function. By the centre of an affine map $x\mapsto w+\mu x$ with $w,\mu\in K$ and $\mu\ne1$, we mean its unique fixed point $w/(1-\mu)\in K$. There are four steps:
\begin{enumerate}[label=\textup{(\arabic*)},nosep]
\item \emph{A mixed symmetry}, with multiplier $\alpha^r\beta^s$ ($r\ge1$, $s\in\N$), on a terminal component of the graph of mixed kernels (\S\ref{ssec:terminal}).
\item \emph{Equality of centres.} If a vertex $G$ of such a component also has a $\beta$-symmetry, the two symmetries have the same centre (Proposition~\ref{prop:centres}). Otherwise, rotations induced by the $\beta$-symmetry, combined with the mixed symmetry, would fill an open subset of a fibre in the $\alpha$-adic completion. By the fibre theorem this would give the closure of $\supp G$ non-empty interior, contrary to Proposition~\ref{prop:emptyint}.
\item \emph{No branching.} Two returns to $G$ through different residue classes modulo $\beta R$ would have the same centre, forcing their translation terms $u,u'\in R$ to satisfy $u\equiv u'\pmod{\beta R}$. It follows that the support is sparse (Lemma~\ref{lem:nobranch} and Proposition~\ref{prop:sparsity}).
\item \emph{Every doubly automatic sparse set is finite}, by the theorem of Evertse, Schlickewei and Schmidt on unit equations (\S\ref{ssec:sparse}). The bases are multiplicatively independent, since an excess prime contains every positive power of $\alpha$ and no positive power of $\beta$.
\end{enumerate}

\subsection{Mixed symmetries on a terminal component}\label{ssec:terminal}
Let $X\in\cI$ and let $\cG$ be the graph of mixed kernels of $\mathbf 1_X$, constructed with fixed window sets $D_\alpha,D_\beta$ as in \S\ref{ssec:graph}. A vertex is called \emph{infinite} if its support is infinite. Let $\cG_\infty$ be the subgraph induced by these vertices. A \emph{terminal component} is a strongly connected component of $\cG_\infty$ from which no edge of $\cG$ leads to an infinite vertex outside the component. Such a component exists: $\cG_\infty$ is finite and non-empty, since $X$ is infinite, and its graph of strongly connected components is acyclic.

\begin{lemma}[Terminal components]\label{lem:terminal}
Let $\cC$ be a terminal component.
\textup{(i)} Each vertex of $\cC$ has an $\alpha$-edge and a $\beta$-edge ending in $\cC$.
\textup{(ii)} Each $G\in\cC$ satisfies $G(w+\varpi y)=G(y)$ for all $y\in R$, for some $w\in R$ and $\varpi=\alpha^r\beta^s$, with $r\ge1$ and $s\in\N$.
\textup{(iii)} The support of each vertex of $\cC$ belongs to $\cI$.
\textup{(iv)} Every infinite vertex reached from $\cC$ belongs to $\cC$.
\end{lemma}
\begin{proof}
(i) For either $\mu\in\{\alpha,\beta\}$,
\[
 \supp G=\bigcup_{e\in D_\mu}\bigl(e+\mu\supp\bigl(y\mapsto G(e+\mu y)\bigr)\bigr).
\]
Since $\supp G$ is infinite, at least one of the corresponding children is infinite and hence belongs to $\cC$.

(ii) Start with an $\alpha$-edge staying in $\cC$ and complete it to a closed path at $G$, using strong connectivity. If the successive digits are $e_1,\ldots,e_m$ and the successive bases are $\mu_1,\ldots,\mu_m$, this path gives the identity in (ii), with
\[
 \varpi=\mu_1\cdots\mu_m,\qquad
 w=\sum_{i=1}^{m}(\mu_1\cdots\mu_{i-1})e_i.
\]
The first edge is of type $\alpha$, so $r\ge1$.

(iii) Each vertex is the characteristic function of a section $X_{\Gamma,c}$. Since the vertices in $\cC$ have infinite support, Proposition~\ref{prop:stability} applies.

(iv) Every child of a finite vertex is finite. Thus a path reaching an infinite vertex cannot pass through a finite vertex, and terminality prevents it from leaving $\cC$.
\end{proof}

\subsection{The completion and its two parts}\label{ssec:completion}
\emph{Normalisation.} Let $h_K$ be the class number of $K$. We replace $\alpha$ by $\alpha^{h_K}$ and keep the notation $\alpha$ for the new base. The part of the new ideal $\alpha R$ supported on the excess primes is principal, allowing us to describe the completion in terms of elements of $R$. This preserves expansivity, $\alpha$-automaticity (Lemma~\ref{lem:powers}), condition (H), the rigid class $\cI$ (automaticity is unchanged, and the other defining conditions do not involve the bases), the completion $R_\alpha$ and the sets of primes. This normalisation is already non-trivial in $K=\Q(\sqrt{10})$: for $\alpha=\sqrt{10}-2$ and $\beta=3$, one has $\alpha R=QP$, where $Q=2R+\sqrt{10}R$ is a non-principal prime above $2$ and $P$ is a prime above $3$. The class number is $2$, and $\alpha^2=2(7-2\sqrt{10})$; after normalisation, $\delta=2$ and $\gamma=7-2\sqrt{10}$ in the notation below.

After normalisation, let $\delta\in R$ generate the ideal $\prod_Q Q^{v_Q(\alpha)}$, where the product is over the excess primes $Q$. Set $\gamma:=\alpha/\delta\in R$, which generates the part of $\alpha R$ supported on the common primes. The elements $\delta$ and $\gamma$ are coprime, and $\delta$ is not a unit. For each non-zero prime ideal $Q$ of $R$, let $K_Q$ be the completion of $K$ at $Q$ and let $\cO_Q$ be its valuation ring, so that $\cO_Q=\varprojlim_n R/Q^n$ and $K\otimes_{\Q}\Q_\ell\simeq\prod_{Q\mid\ell}K_Q$ for every rational prime $\ell$ \cite[Ch.~II, \S\S4 and 8]{Neukirch}. By the Chinese remainder theorem,
\[
\begin{split}
 R_\alpha&=\cA\times\cB,\\
 \cA&=\varprojlim_n R/\delta^nR=\prod_{Q\text{ excess}}\cO_Q,\\
 \cB&=\varprojlim_n R/\gamma^nR=\prod_{P\text{ common}}\cO_P.
\end{split}
\]
We use the zero ring for the empty product in the last formula. The element $\beta$ is a unit in $\cA$ and is topologically nilpotent in $\cB$; in particular, $1-\beta^k$ is a unit in $\cB$ for $k\ge1$.

Let $\mathcal L$ be the set of rational primes $\ell$ below the excess primes. For $\ell\in\mathcal L$, put
\[
 \cA_\ell:=\prod_{\substack{Q\mid\ell\\Q\text{ excess}}}K_Q,
 \qquad n_\ell:=\dim_{\Q_\ell}\cA_\ell.
\]
Then $\cA\otimes_{\Z}\Q=\prod_{\ell\in\mathcal L}\cA_\ell$. The natural map $K\otimes_{\Q}\Q_\ell\to\cA_\ell$ is surjective, being the projection onto the indicated factors. We also set $\Z_{\mathcal L}:=\prod_{\ell\in\mathcal L}\Z_\ell$, in which $\N$ is dense under the diagonal embedding.

\begin{example}
For $\alpha=7$ and $\beta=3+\sqrt2$ in $R=\Z[\sqrt2]$, we have $\delta=3-\sqrt2$ and $\gamma=3+\sqrt2$. With $P,\bar P$ as in Example~\ref{ex:parity},
\[
 R_7=\cO_{\bar P}\times\cO_P\simeq\Z_7\times\Z_7.
\]
The affine map $x\mapsto e+\beta^kx$ ($e\in R$, $k\ge1$) acts on the first coordinate as an affine map with unit multiplier $\beta^k$, and contracts the second coordinate.
\end{example}

\subsection{Logarithms and the open mapping theorem}\label{ssec:logs}
To obtain continuous families of rotations from a single $\beta$-symmetry, we use $\ell$-adic powers of units.

Set $m_\ell=1$ for odd $\ell$ and $m_2=2$. Recall \cite[Ch.~5]{Robert} that for a prime $Q$ above $\ell$, the exponential and logarithm are mutually inverse bijections
\[
 \ell^{m_\ell}\cO_Q\ \overset{\exp}{\underset{\log}{\rightleftarrows}}\ 1+\ell^{m_\ell}\cO_Q,
\]
identifying the additive group on the left with the multiplicative group on the right. For $V\in1+\ell^{m_\ell}\cO_Q$ and $t\in\Z_\ell$, the expression $V^t:=\exp(t\log V)$ extends integer powers and satisfies $V^{t+t'}=V^tV^{t'}$. It is a \emph{restricted power series} in $t$, meaning that its coefficients tend to zero: this follows from Legendre's formula for the valuation of a factorial. If $V\ne1$, then $\log V\ne0$.

If $V\in\cA$ has all its components in the indicated groups, we define $V^t\in\cA$ for $t\in\Z_{\mathcal L}$ componentwise, using the $\ell$-coordinate of $t$ in each factor above $\ell$. The derivative of the map $t\mapsto V^t$ on a factor $K_Q$ at $t_0\in\Z_\ell$ is $u\mapsto u V^{t_0}\log V$, for $u\in\Q_\ell$.

\begin{lemma}[Open mapping]\label{lem:ift}
Let $W$ be an $n$-dimensional $\Q_\ell$-vector space and let $F:\Z_\ell^n\to W$ be given by a restricted power series. If the differential $\mathrm dF(t_0):\Q_\ell^n\to W$ at $t_0\in\Z_\ell^n$ is an isomorphism, then the image of every sufficiently small ball centred at $t_0$ contains a ball centred at $F(t_0)$.
\end{lemma}
\begin{proof}
This is the ultrametric local inverse function theorem. Equip $\Q_\ell^n$ with its usual supremum norm and $W$ with the norm induced by a basis. Write
\[
 F(t_0+\eta)=F(t_0)+J\eta+\Delta(\eta),
 \qquad J=\mathrm dF(t_0).
\]
The remainder $\Delta$ vanishes at $0$ and has no linear term. On a ball of radius $\varepsilon\le1$ it is Lipschitz with constant at most $M\varepsilon$ for a fixed $M$, by the boundedness of the coefficients of its restricted power series. Take $\varepsilon$ small enough that $\norm{J^{-1}}M\varepsilon<1$. For $y\in W$ with $\norm{J^{-1}(y-F(t_0))}\le\varepsilon$, the map
\[
 \eta\longmapsto J^{-1}\bigl(y-F(t_0)-\Delta(\eta)\bigr)
\]
is a contraction of the closed ball of radius $\varepsilon$ into itself. Its fixed point satisfies $F(t_0+\eta)=y$, proving the assertion.
\end{proof}

\subsection{Rotations induced by a \texorpdfstring{$\beta$}{beta}-symmetry}\label{ssec:rotations}
A $\beta$-symmetry of a set induces, on its closure in $R_\alpha$, rotations about its centre in the factor $\cA$, parametrised continuously by $\Z_{\mathcal L}$.

\begin{lemma}[Collapsed rotations]\label{lem:rotations}
Let $E\subseteq R$ satisfy $\mathbf 1_E(e+\beta^ky)=\mathbf 1_E(y)$ for all $y\in R$, for some $e\in R$ and $k\ge1$. Put $c=e/(1-\beta^k)\in K$ and let $\overline E$ be the closure of $E$ in $R_\alpha$. There exists $h\ge1$ such that $V:=\beta^{kh}$ has its components at the excess primes in $1+\ell^{m_\ell}\cO_Q$. For each $t\in\Z_{\mathcal L}$, the map
\[
 \Psi_t(x)=\bigl(c_{\cA}+V^t(x_{\cA}-c_{\cA}),\,c_{\cB}\bigr)
\]
is a continuous self-map of $R_\alpha$ and sends $\overline E$ into $\overline E$. Here subscripts denote projections in the product of the corresponding local fields.
\end{lemma}
\begin{proof}
The element $\beta$ is a unit in each of the finite rings $\cO_Q/\ell^{m_\ell}\cO_Q$ for $Q$ excess. A suitable positive power of $\beta^k$ is therefore $1$ in all these rings. Choose $h$ accordingly. Also $c_{\cB}\in\cB$ since $1-\beta^k$ is a unit there.

The affine map $\psi(x)=e+\beta^kx=c+\beta^k(x-c)$ is a continuous self-map of $R_\alpha$ and sends $\overline E$ into itself, as do all its iterates $\psi^{hn}(x)=c+V^n(x-c)$ for $n\in\N$. Given $t\in\Z_{\mathcal L}$, choose $n_j\in\N$ tending to $+\infty$ in the usual sense and to $t$ in $\Z_{\mathcal L}$. Such a sequence exists by the Chinese remainder theorem: at each stage, the required congruences are realised by arbitrarily large elements of $\N$. Then $\psi^{hn_j}(x)\to\Psi_t(x)$ in the product of local fields. On $\cA$ this follows from continuity of $t\mapsto V^t$; on $\cB$ it follows from $V^{n_j}\to0$.

As $R_\alpha$ is closed in that product, $\Psi_t(x)$ lies in $R_\alpha$, even if $c_{\cA}\notin\cA$. The formula is continuous in $x$, and the closedness of $\overline E$ gives the desired preservation property.
\end{proof}

\subsection{Choosing a basis}\label{ssec:select}
The following elementary lemma provides, under (H), the $\Q$-basis of $K$ used in \S\ref{ssec:centres} to show that a differential is surjective.

\begin{lemma}[Selection]\label{lem:select}
Let $x,y\in K^\times$ satisfy $\Q(x,y)=K$, and put $e=[\Q(x):\Q]$. Then the elements $x^iy^{\lfloor i/e\rfloor}$, for $0\le i<d$, form a $\Q$-basis of $K$.
\end{lemma}
\begin{proof}
Write $d=eg$ and $i=je+k$ with $0\le k<e$. For fixed $j$, the elements $x^{je+k}y^j=x^{je}y^j x^k$ form a $\Q$-basis of $y^j\Q(x)$, since $1,x,\ldots,x^{e-1}$ form a basis of $\Q(x)$ and $x^{je}$ is non-zero. As $[\Q(x)(y):\Q(x)]=g$, we have
\[
 K=\bigoplus_{j<g}y^j\Q(x).
\]
The asserted basis is the union of the bases of these summands.
\end{proof}

\subsection{From a fibre to an open set}\label{ssec:fibre}
A \emph{fibre} is a subset of $R_\alpha$ of the form $\cA\times\{b\}$, with $b\in\cB$. When there are no common primes, $\cB=0$ and a fibre is open.

\begin{lemma}[Balanced representatives]\label{lem:units}
There exists a constant $C_K>0$ depending only on $K$ such that, for every $z\in K^\times$, each coset of the lattice $zR$ in $K_{\R}$ contains an element $w$ with
\[
 \norm w\le C_K|N_{K/\Q}(z)|^{1/d}.
\]
\end{lemma}
\begin{proof}
By Dirichlet's unit theorem \cite[Ch.~I, \S7]{Neukirch}, the map $x\mapsto(\log|\sigma x|)_{\sigma\in\Sinf}$ sends $R^\times$ onto a full-rank lattice in the hyperplane $\sum_\sigma d_\sigma x_\sigma=0$; let $r_K$ be its covering radius for the supremum norm. Subtracting $d^{-1}\log|N_{K/\Q}(z)|$ from each archimedean logarithm of $z$ puts the vector in this hyperplane, and translating by a logarithmic unit vector brings all coordinates within $r_K$ of zero. This yields $\varepsilon_z\in R^\times$ with $\norm{z\varepsilon_z}\le e^{r_K}|N_{K/\Q}(z)|^{1/d}$. If $\mathcal P_0\subset K_{\R}$ is a fixed bounded fundamental domain for $R$, then $z\varepsilon_z\mathcal P_0$ is one for $zR$, and its elements have norm at most $\norm{z\varepsilon_z}\sup_{p\in\mathcal P_0}\norm p$; we may take $C_K=e^{r_K}\sup_{p\in\mathcal P_0}\norm p$.
\end{proof}

The following theorem says that if the closure of an automatic set contains a relatively open subset of a fibre, it also contains an open subset of the completion. In the proof, the relevant sections of $E$ form a finite family; for each of their supports, the parameters of the fibres contained in its closure form a closed set; these closed sets have dense union, hence cover $\cB$, and one of them therefore has non-empty interior.

\begin{theorem}[Fibre theorem]\label{thm:fibre}
Let $\alpha\in R$ be expansive and suppose that $\alpha=\delta\gamma$, with $\delta,\gamma\in R$ coprime non-units. Let $E\subseteq R$ be $\alpha$-automatic. Write $R_\alpha=\cA\times\cB$ with $\cA=\varprojlim_j R/\delta^jR$ and $\cB=\varprojlim_j R/\gamma^jR$. Write $\overline E$ for the closure of $E$ in $R_\alpha$. If $(a+\delta^n\cA)\times\{b\}\subseteq\overline E$ for some $a\in\cA$, $b\in\cB$ and $n\in\N$, then $\overline E$ has non-empty interior in $R_\alpha$.
\end{theorem}
\begin{proof}
Let $D$ be a window set for $\alpha$.

\emph{A finite family of sections.} The functions
\[
 y\longmapsto\mathbf 1_E(u+\alpha^jy),\qquad
 j\in\N,\quad u\in R,\quad |\sigma u|\le|\sigma\alpha|^j
 \quad(\sigma\in\Sinf),
\]
form a finite family $H_1,\ldots,H_M$. Indeed, writing $u=u'+\alpha^jv$ with $u'\in D^{(j)}$ and $v\in R$, Lemma~\ref{lem:fixedlength} gives $|\sigma v|\le1+C_D$ at every place, independently of $j$ and $u$. Thus $v$ lies in a finite set (Lemma~\ref{lem:proper}), and each such function is the translate by $v$ of an $\alpha$-kernel.

If $H(y)=\mathbf 1_E(u+\alpha^jy)$ for $y\in R$, with $u\in R$ and $j\in\N$, then
\[
 \overline{\supp H}=\{z\in R_\alpha:u+\alpha^jz\in\overline E\}.
\]
Indeed, $z\mapsto u+\alpha^jz$ is a homeomorphism of $R_\alpha$ onto the clopen set $U:=u+\alpha^jR_\alpha$. Since $U\cap R=u+\alpha^jR$, it carries $\supp H$ onto $E\cap U$, whose closure is $\overline E\cap U$. The sets
\[
 \mathcal F_i:=\{\tau\in\cB:\cA\times\{\tau\}\subseteq\overline{\supp H_i}\}
\]
are closed: each is the intersection, over $x\in\cA$, of the preimages of $\overline{\supp H_i}$ under the continuous maps $\tau\mapsto(x,\tau)$.

\emph{Every $\tau\in\cB$ is close to some $\mathcal F_i$.} Let $\tau\in\cB$ and $N\in\N$; we find $i$ and $\tau'\in\mathcal F_i$ with $\tau'-\tau\in\gamma^N\cB$. For $j\ge n$, the elements $u\in R$ such that
\[
 u_{\cA}\equiv a\pmod{\delta^n\cA},\qquad
 u_{\cB}+\alpha^j\tau\equiv b\pmod{\gamma^{j+N}\cB}
\]
form a coset of $\delta^n\gamma^{j+N}R$, by the Chinese remainder theorem. Dividing by $\alpha^j$ gives a coset of $zR$, where $z=\delta^{n-j}\gamma^N\in K^\times$. By Lemma~\ref{lem:units}, we may choose $u$ in the original coset so that
\[
 \norm{u/\alpha^j}\le
 C_K\bigl(|N_{K/\Q}(\delta)|^{n-j}|N_{K/\Q}(\gamma)|^N\bigr)^{1/d}\le1,
\]
the last inequality holding for $j$ large enough in terms of $n$ and $N$, since $|N_{K/\Q}(\delta)|\ge2$. Then $|\sigma u|\le|\sigma\alpha|^j$ at every place, so the function $y\mapsto\mathbf 1_E(u+\alpha^jy)$ is $H_i$ for some $i$.

Put $\tau':=\alpha^{-j}(b-u_{\cB})$. Since $\delta$ is a unit in $\cB$, the second congruence gives
\[
 \tau'-\tau\in\alpha^{-j}\gamma^{j+N}\cB=\gamma^N\cB;
\]
in particular $\tau'\in\cB$. For every $x\in\cA$, the point $u+\alpha^j(x,\tau')$ has $\cB$-coordinate $b$ and $\cA$-coordinate in $a+\delta^n\cA$, because $\alpha^j\cA=\delta^j\cA\subseteq\delta^n\cA$; hence it lies in $\overline E$ by hypothesis. By the identity above, $\cA\times\{\tau'\}\subseteq\overline{\supp H_i}$, that is, $\tau'\in\mathcal F_i$.

\emph{The covering.} Since the ideals $\gamma^N\cB$ form a neighbourhood basis at $0$, the union of the $\mathcal F_i$ is dense in $\cB$. Being a finite union of closed sets, it is closed, hence equal to $\cB$. A finite union of closed sets with empty interior has empty interior, so some $\mathcal F_i$ contains a non-empty open set $U\subseteq\cB$. Writing $H_i(y)=\mathbf 1_E(u_0+\alpha^{j_0}y)$, we conclude that the non-empty open set $u_0+\alpha^{j_0}(\cA\times U)$ is contained in $\overline E$.
\end{proof}

\subsection{Equality of centres}\label{ssec:centres}
Suppose that a mixed symmetry $\varphi$ with multiplier $\varpi$ and a $\beta$-symmetry have distinct centres $q,c\in K$. Write $\vartheta=q-c$. Combining $\varphi$ with the rotations of Lemma~\ref{lem:rotations}, in the form $\Psi_t\circ\varphi\circ\Psi_{-t}$, makes part of its translation term, $(1-\varpi)V^t\vartheta$, vary with $t$. By composing $d$ such maps, we obtain a family of points of $\overline E$ in the same fibre, depending on $d$ parameters in $\Z_{\mathcal L}$. By (H) and the selection lemma, its differential is surjective on each $\cA_\ell$. The family therefore fills a relatively open subset of the fibre. The fibre theorem would then give $\overline E$ non-empty interior, contrary to Proposition~\ref{prop:emptyint}.

\begin{proposition}[Equality of centres]\label{prop:centres}
Let $E\in\cI$ and suppose that, for every $y\in R$,
\[
 \mathbf 1_E(w+\varpi y)=\mathbf 1_E(y),\qquad
 \mathbf 1_E(e+\beta^ky)=\mathbf 1_E(y),
\]
where $w,e\in R$, $\varpi=\alpha^r\beta^s$ ($r,k\ge1$, $s\in\N$). Then $w/(1-\varpi)=e/(1-\beta^k)$.
\end{proposition}
\begin{proof}
Put $q=w/(1-\varpi)$ and $c=e/(1-\beta^k)$, and suppose for a contradiction that $\vartheta:=q-c\ne0$.

\emph{1. The composite maps.} When $\cB\ne0$, the maps $\Psi_t$ are not invertible, since they send the $\cB$-coordinate to $c_{\cB}$; their restrictions to $\cA\times\{c_{\cB}\}$ are invertible, with inverses the restrictions of $\Psi_{-t}$. The centres lie in $K$, which embeds in $R_\alpha\otimes_{\Z}\Q$; in particular, we do not assume that $c_{\cA}$ is integral. By Lemma~\ref{lem:rotations}, the maps $\Psi_t$ send $\overline E$ into itself. So does the affine map $\varphi(x)=w+\varpi x=q+\varpi(x-q)$, by continuity. Hence each composite $\varphi_t:=\Psi_t\circ\varphi\circ\Psi_{-t}$ sends $\overline E$ into itself. Its $\cB$-coordinate is constant, equal to $c_{\cB}$, and
\[
 (\varphi_t(x))_{\cA}=\varpi x_{\cA}+(1-\varpi)(c_{\cA}+V^t\vartheta_{\cA}).
\]
Here and below, elements of $K$ multiplying elements of a product of local fields are understood through their natural images.

\emph{2. A family of points in one fibre.} Fix $x_0\in E$ and, for $t=(t_0,\ldots,t_{d-1})\in\Z_{\mathcal L}^d$, set
\[
 F(t):=\varphi_{t_0}\circ\cdots\circ\varphi_{t_{d-1}}(x_0)\in\overline E.
\]
Then $F(t)_{\cB}=c_{\cB}$ and
\[
 F(t)_{\cA}=\varpi^d(x_0)_{\cA}+(1-\varpi)\sum_{i<d}\varpi^i(c_{\cA}+V^{t_i}\vartheta_{\cA}).
\]
Write $F_\ell$ for the projection of $F$ to $\cA_\ell$, and $t_{i,\ell}$ for the $\ell$-coordinate of $t_i$. Then
\[
 \frac{\partial F_\ell}{\partial t_{i,\ell}}=\xi_\ell\varpi^iV^{t_{i,\ell}},
 \qquad \xi_\ell=(1-\varpi)\vartheta\log V\in\cA_\ell,
\]
where the logarithm and the power are taken componentwise.

\emph{3. The common factor is invertible.} In the product of fields $\cA_\ell$, an element is invertible if none of its components vanishes. This holds for $\xi_\ell$: the elements $1-\varpi$ and $\vartheta$ lie in $K^\times$ ($\varpi$ is expansive, hence distinct from $1$), and $\log V\ne0$ in every $K_Q$, because $V=\beta^{kh}\ne1$ in $K$, which embeds in each $K_Q$, and the logarithm is injective on the subgroup in question.

\emph{4. The directions supplied by (H).} The field $\Q(\varpi,V)$ contains
\[
 \varpi^{kh}V^{-s}=\alpha^{rkh},\qquad V=\beta^{kh}.
\]
It is therefore $K$ by (H). By Lemma~\ref{lem:select}, the elements
\[
 \varpi^iV^{j_i},\qquad 0\le i<d,\quad
 j_i:=\left\lfloor\frac{i}{[\Q(\varpi):\Q]}\right\rfloor,
\]
form a $\Q$-basis of $K$; their images form a $\Q_\ell$-basis of $K\otimes_{\Q}\Q_\ell$ and span its quotient $\cA_\ell$. At the point $t^*=(j_0,\ldots,j_{d-1})\in\Z_{\mathcal L}^d$ (with the integers embedded diagonally), the partial derivatives are these images multiplied by $\xi_\ell$, so they span $\cA_\ell$ by Step~3. This is the only use of (H) in the proof.

\emph{5. An open subset of the fibre.} Each $F_\ell$ depends only on the $\ell$-coordinates of the parameters and is a sum of restricted power series in one variable each. For every $\ell\in\mathcal L$, choose $n_\ell$ indices $i$, possibly depending on $\ell$, whose partial derivatives at $t^*$ form a $\Q_\ell$-basis of $\cA_\ell$, and fix the other coordinates $t_{i,\ell}$ at $j_i$. This gives a restricted power series on $\Z_\ell^{n_\ell}$ with invertible differential at the corresponding point, whose image contains an open neighbourhood of $F_\ell(t^*)$ by Lemma~\ref{lem:ift}.

Since $\Z_{\mathcal L}=\prod_{\ell\in\mathcal L}\Z_\ell$, the parameters at distinct rational primes can be chosen independently. Taking the product of these neighbourhoods, the image of $t\mapsto F(t)_{\cA}$ contains an open neighbourhood of $a:=F(t^*)_{\cA}$ in $\cA\otimes_{\Z}\Q$. As $a\in\cA$ (because $F(t^*)\in\overline E$) and the ideals $\delta^n\cA$ form a neighbourhood basis at $0$, this neighbourhood contains $a+\delta^n\cA$ for $n$ large. Thus
\[
 (a+\delta^n\cA)\times\{c_{\cB}\}\subseteq\overline E.
\]

\emph{6. The contradiction.} If $\cB=0$, this set is open in $R_\alpha$; otherwise Theorem~\ref{thm:fibre} applies. In either case $\overline E$ has non-empty interior, contrary to Proposition~\ref{prop:emptyint}.
\end{proof}

\begin{remark}
The mixed symmetry supplies directions which the rotations alone need not supply. Take $K=\Q(\sqrt2,\sqrt3)$, $\alpha=22(1+\sqrt2)$ and $\beta=2(4+\sqrt3)$. Condition (H) holds, but the $11$-adic powers of $\beta$ lie in the image of $\Q(\sqrt3)\otimes_{\Q}\Q_{11}$, a subalgebra of dimension $2$ in $\cA_{11}$, which has dimension $4$. Indeed, all primes above $11$ are excess primes: $1+\sqrt2$ is a unit and $N_{\Q(\sqrt3)/\Q}(4+\sqrt3)=13$. To illustrate the algebraic mechanism of Lemma~\ref{lem:select}, take $\varpi=\alpha$ and let $V$ be a positive power of $\beta$. We have $[\Q(\alpha):\Q]=2$ and $\Q(\alpha,V)=K$ by (H). Hence the elements $1,\alpha,\alpha^2V,\alpha^3V$ form a basis of $K$ over $\Q$. This example concerns only the choice of directions; it does not assert the existence of a mixed symmetry with multiplier $\alpha$.
\end{remark}

\subsection{No branching}\label{ssec:nobranch}
For each residue class $\kappa\in R/\beta R$, fix a representative $e_\kappa\in D_\beta$. Given a function $G:R\to\F_2$, put $G_\kappa(y):=G(e_\kappa+\beta y)$. Then
\[
 \supp G=\bigsqcup_{\kappa\in R/\beta R}(e_\kappa+\beta\supp G_\kappa).
\]
For a vertex $G_0$ of $\cG$, the \emph{class graph} of $G_0$ has as vertices the functions obtained from $G_0$ by repeatedly taking $G_\kappa$, including $G_0$ itself. Draw one edge for each residue class. This is a finite directed multigraph whose edges are $\beta$-edges of $\cG$.

\begin{lemma}[No branching]\label{lem:nobranch}
Let $X\in\cI$, let $\cC$ be a terminal component and let $G_0\in\cC$. In the class graph of $G_0$ with the zero function removed, every strongly connected component is either a simple directed cycle or a single vertex without a loop.
\end{lemma}
\begin{proof}
Suppose that a vertex $G$ has two outgoing edges $G\to G_\kappa$ and $G\to G_{\kappa'}$, with $\kappa\ne\kappa'$, within its strongly connected component.

If $G$ is infinite, then $G\in\cC$ by Lemma~\ref{lem:terminal}\,(iv), so $G$ has a mixed symmetry with some centre $q\in K$, and $\supp G\in\cI$ (parts (ii) and (iii) of the same lemma). Completing each edge to a closed path at $G$ gives returns
\[
 G(u_i+\beta^{k_i}y)=G(y)\qquad(y\in R,\ i=1,2),
\]
where $u_1,u_2\in R$, $k_1,k_2\ge1$, $u_1\equiv e_\kappa$ and $u_2\equiv e_{\kappa'}\pmod{\beta R}$. By Proposition~\ref{prop:centres},
\[
 \frac{u_1}{1-\beta^{k_1}}=q=\frac{u_2}{1-\beta^{k_2}}.
\]
Cross-multiplying in $R$ and reducing modulo $\beta R$ gives $u_1\equiv u_2\pmod{\beta R}$, contradicting $\kappa\ne\kappa'$.

If $G$ has finite non-empty support, the cardinality of the support cannot increase along an edge, so it is constant, say $M_0\ge1$, on the component. But $\supp G$ contains the disjoint sets $e_\kappa+\beta\supp G_\kappa$ and $e_{\kappa'}+\beta\supp G_{\kappa'}$, so $M_0\ge2M_0$, again a contradiction.

Thus each vertex has at most one outgoing edge within its component. By strong connectivity, a component containing a closed path of positive length is a single simple directed cycle, and any other component is a single vertex without a loop.
\end{proof}

\begin{proposition}[Sparsity]\label{prop:sparsity}
Under the hypotheses of Lemma~\ref{lem:nobranch}, the support of $G_0$ is sparse. More precisely, there are constants $C>0$ and $h_0\in\N$ such that, for every real $N\ge0$,
\[
 \#\{x\in\supp G_0:\norm x\le N\}\le C(\log(2+N))^{h_0}.
\]
\end{proposition}
\begin{proof}
Let $h_0$ be the number of strongly connected components of the class graph, with the zero function removed, which are cycles. A path of length $n$ starting at $G_0$ and avoiding the zero function cannot return to a strongly connected component after leaving it. It is determined by the components visited, the edges joining them, its final vertex and the number of complete turns in each cycle. All these data except the numbers of turns range over a fixed finite set, so there are at most $C_1(n+1)^{h_0}$ such paths, for a constant $C_1$ depending only on the graph.

For $y\in\supp G_0$, put $y_0=y$ and recursively $\kappa_j=y_j+\beta R$, $y_{j+1}=(y_j-e_{\kappa_j})/\beta\in R$. This defines a path in the class graph, and $y_j$ lies in the support of its $j$-th vertex, so the path avoids the zero function. With $M_\beta:=\max_\kappa\norm{e_\kappa/\beta}$ and $\rho_\beta=\max_{\sigma\in\Sinf}|\sigma\beta|^{-1}<1$ (Definition~\ref{def:expansive}),
\[
 \norm{y_n}\le\rho_\beta^n\norm y+\frac{M_\beta}{1-\rho_\beta}.
\]
If $\norm y\le\rho_\beta^{-n}$, the element $y_n$ therefore lies in a fixed finite subset of $R$ (Lemma~\ref{lem:proper}). Since
\[
 y=\sum_{j<n}\beta^je_{\kappa_j}+\beta^ny_n,
\]
the point $y$ is determined by the path, which records the digits, and by $y_n$. The number of points $y\in\supp G_0$ with $\norm y\le\rho_\beta^{-n}$ is therefore bounded by a constant times $(n+1)^{h_0}$. Taking $n=\lceil\log N/\log(\rho_\beta^{-1})\rceil$ for $N\ge1$, and adjusting $C$ for $0\le N<1$, gives the bound.
\end{proof}

\subsection{Doubly automatic sparse sets are finite}\label{ssec:sparse}
A sparse set can be infinite: for instance, the set $\{7^n:n\in\N\}$ in $\Z$ is sparse and $7$-automatic. Finiteness will follow from multiplicative independence of the bases. The following theorem uses neither (H) nor the existence of an excess prime, and holds in any degree.

\begin{theorem}[Sparse sets]\label{thm:sparse}
Let $K$ be a number field and let $R=\cO_K$. Let $\alpha,\beta\in R$ be expansive and multiplicatively independent ($\alpha^m\ne\beta^n$ for all $m,n\ge1$). A sparse subset $X\subseteq R$ which is both $\alpha$- and $\beta$-automatic is finite.
\end{theorem}
\begin{proof}
For $\mu\in\{\alpha,\beta\}$, fix an admissible digit set $D'_\mu\subseteq R$ containing $0$ (Lemma~\ref{lem:short}) and a total order on its digits. We show that, in each base, the elements of $X$ are the values of finitely many exponential polynomials, and then compare the two bases by means of a unit equation.

\emph{1. Unique representatives.} Two words $w,w'$ of the same length $n$ over $D'_\mu$ have the same value if and only if the carries $t_0=0$, $t_{i+1}=(t_i+w_i-w'_i)/\mu$ all lie in $R$ and $t_n=0$: if the values are equal, the partial differences of the values are divisible by the appropriate powers of $\mu$, so all the carries lie in $R$. These carries have norm at most $2b_D\rho_\mu/(1-\rho_\mu)$, where $b_D:=\max_{e\in D'_\mu}\norm e$ and $\rho_\mu$ is as in Definition~\ref{def:expansive}. There are therefore finitely many of them, and a finite automaton reading pairs of digits, least significant digit first, recognises the pairs of words of equal length and equal value.

The language $L_0$ of all expansions of elements of $X$ over $D'_\mu$ is regular (Lemma~\ref{lem:kernel} and Corollary~\ref{cor:windowcrit}). Since expansions are read least significant digit first, a value has an expansion of length less than $n\ge1$ if and only if it has an expansion of length $n$ ending in $0$: pad with trailing zeros, or delete the final zero. Let $L$ consist of the words $w\in L_0$ such that no expansion of $\val w_\mu$ of the same length as $w$ ends in $0$ or is lexicographically smaller than $w$. Each $x\in X$ is then the value of exactly one word of $L$, namely the lexicographically least of its shortest expansions (the empty word if $x=0$). The excluded words are the first coordinates of the pairs $(w,w')$ of equal length and equal value with $w'$ ending in $0$, or with $w'<_{\mathrm{lex}}w$; both conditions are recognised by finite automata, so $L$ is regular by closure under projection and difference.

\emph{2. Polynomial growth.} Since $\norm\mu>1$, every word of length $n$ has value of norm at most $C_\mu\norm\mu^n$, for a constant $C_\mu>0$ depending only on $D'_\mu$ and $\mu$. As $X$ is sparse and the value map is injective on $L$, the number of words of length $n$ in $L$ is bounded by a polynomial in $n$. By the characterisation of regular languages of polynomial growth \cite{SYZS}, $L$ is then a finite union of languages
\[
 w_0u_1^*w_1\cdots u_k^*w_k,
\]
with fixed words $w_i$ and non-empty words $u_i$ over $D'_\mu$. For completeness, we recall why. Take a trim deterministic automaton for $L$ (allowing partial transitions), that is, one in which every state is reachable from the initial state and can reach an accepting state. No state has two outgoing edges within its strongly connected component. Otherwise these edges would start closed paths with labels $v_1,v_2$ beginning with different letters. The words $v_1^{|v_2|}$ and $v_2^{|v_1|}$ have the same length $p=|v_1||v_2|$. Concatenating $m$ arbitrary choices of them, between fixed initial and accepting paths, would then give $2^m$ distinct accepted words of length $mp+O(1)$, contrary to polynomial growth. The strongly connected components are therefore simple cycles or single vertices without loops. Once the components visited, the edges between them and the entry and exit points on each cycle are chosen, among finitely many possibilities, only the numbers of complete turns remain free; this gives the displayed form.

Fix one of these languages and write $b_i=|u_i|\ge1$. The identities $\val{ab}_\mu=\val a_\mu+\mu^{|a|}\val b_\mu$ and
\[
 \val{u_i^{n_i}}_\mu=\val{u_i}_\mu\,\frac{1-\mu^{b_in_i}}{1-\mu^{b_i}}
 \qquad(n_i\in\N),
\]
where $\mu^{b_i}\ne1$ because $\mu$ is expansive, give
\[
 \val{w_0u_1^{n_1}w_1\cdots u_k^{n_k}w_k}_\mu
 =\sum_{T\subseteq\{1,\ldots,k\}}\lambda_T\,\mu^{\sum_{i\in T}b_in_i},
\]
with coefficients $\lambda_T\in K$ depending only on the fixed words and on $\mu$. Terms with zero coefficients are omitted.

\emph{3. The unit equation.} Fix one language of this form in base $\beta$ and one in base $\alpha$; there are finitely many such pairs. A common value $x$ has two expressions
\[
 x=\sum_T\lambda_T\beta^{\epsilon_T}
   =\sum_{T'}\lambda'_{T'}\alpha^{\epsilon'_{T'}},
\]
with fixed coefficients and exponents in $\N$. Subtracting gives a relation
\[
 \sum_{i=1}^{N_0}g_iX_i=0
\]
with fixed $g_i\in K^\times$, where each $X_i$ is a power of $\alpha$ or of $\beta$ and each index keeps its base even when $X_i=1$; if there are no terms, then $x=0$. Choosing minimal vanishing subsums successively, we partition the indices into blocks $J$ whose sum vanishes while no non-empty proper subsum does. There are finitely many such partitions, and no block is a singleton since $g_iX_i\ne0$. For a block $J$ and a fixed $i_J\in J$, division by $-g_{i_J}X_{i_J}$ gives
\[
 \sum_{i\in J\setminus\{i_J\}}
 \frac{-g_i}{g_{i_J}}\frac{X_i}{X_{i_J}}=1,
\]
where, by minimality of $J$, no non-empty subsum of the left-hand side vanishes. The theorem of Evertse, Schlickewei and Schmidt \cite[Thm.~1.1]{ESS} asserts that an equation $a_1z_1+\cdots+a_mz_m=1$ with $a_i\in K^\times$ has only finitely many solutions in a finite-rank subgroup of $(\overline K^\times)^m$ for which no non-empty proper subsum of the left-hand side vanishes. Applied to $\langle\alpha,\beta\rangle^{|J|-1}$, the indicated direct power of the subgroup of $K^\times$ generated by the bases, it shows that the ratios $X_i/X_{i_J}$ ($i\in J$) take only finitely many values.

The homomorphism
\[
 \Z^2\longrightarrow K^\times,\qquad (m,n)\longmapsto\alpha^m\beta^n,
\]
is injective: by expansivity at any archimedean place, a non-trivial relation $\alpha^m\beta^n=1$ has both exponents non-zero and of opposite signs, and would give an equality of positive powers of the two bases, excluded by multiplicative independence. Hence, in a block containing a power $\beta^\epsilon$ and a power $\alpha^{\epsilon'}$, the finitely many values of $\alpha^{-\epsilon'}\beta^\epsilon$ give finitely many exponents $\epsilon,\epsilon'$, hence finitely many values of these two powers and, through the ratios, of all the $X_i$ in the block.

Finally, $x$ is the sum of the $\beta$-terms of the relation: blocks containing only powers of $\beta$ contribute $0$, blocks containing only powers of $\alpha$ contain no $\beta$-term, and mixed blocks contribute finitely many values. Thus $x$ takes finitely many values for each pair of languages and each partition, and $X$ is finite.
\end{proof}

\subsection{The theorem for different supports}\label{ssec:differentmain}
\begin{theorem}[Different supports]\label{thm:different}
Let $K$ be a number field of degree $d\ge2$, let $R=\cO_K$, and let $\alpha,\beta\in R$ be expansive elements satisfying \textup{(H)}. Suppose that a prime ideal of $R$ contains one of the bases but not the other. For a finite alphabet $\Omega$, a function $f:R\to\Omega$ is automatic in both bases if and only if it is eventually periodic.
\end{theorem}
\begin{proof}
By symmetry, we may suppose that the prime ideal contains $\alpha$ but not $\beta$, and normalise as in \S\ref{ssec:completion}. The converse is Proposition~\ref{prop:converse}.

Let $f$ be doubly automatic. Each level set $\{x\in R:f(x)=\omega\}$, with $\omega\in\Omega$, is doubly automatic by Corollary~\ref{cor:closure}. Suppose that each of these sets agrees outside a finite set with a set $G_\omega$ periodic under a subgroup $\Lambda_\omega\subseteq R$ of rank $d$. Then $\Lambda:=\bigcap_{\omega\in\Omega}\Lambda_\omega$ has finite index in $R$. Each $G_\omega$ is a union of cosets of $\Lambda$, and every such coset is infinite. As the original level sets form a partition and there are only finitely many exceptional points, each coset of $\Lambda$ belongs to exactly one $G_\omega$. Thus the $G_\omega$ form a partition and define a $\Lambda$-periodic function agreeing with $f$ outside a finite set.

It therefore suffices to treat the case $f=\mathbf 1_X$, for $X\subseteq R$. By Theorem~\ref{thm:reduction}, it suffices to show that $\cI=\varnothing$. Suppose that $X\in\cI$, choose a terminal component $\cC$ and a vertex $G_0\in\cC$, and put $E_0:=\supp G_0\in\cI$. By Lemma~\ref{lem:nobranch} and Proposition~\ref{prop:sparsity}, $E_0$ is sparse. Since the bases are multiplicatively independent (step (4) above), Theorem~\ref{thm:sparse} implies that $E_0$ is finite, whereas the members of $\cI$ are infinite. This contradiction proves the theorem.
\end{proof}

\begin{remark}
The excess prime has three roles. It excludes aperiodic minimal subsystems (Theorem~\ref{thm:finiteorbits}). It gives a non-zero factor $\cA$ in which $\beta$ is a unit, and a non-unit $\delta$, as required by the rotation argument and the fibre theorem respectively. Finally, it ensures multiplicative independence of the bases. Condition (Ind) also follows from its existence.
\end{remark}
\section{Transfer between bases}\label{sec:transfer}

This section collects the tools used in the case of equal supports (Section~\ref{sec:equal}), where the method of Sections~\ref{sec:orbits} and~\ref{sec:different} is no longer available: transfer of automaticity between bases, and exact symmetries of kernels and of their sections. As in Section~\ref{sec:numeration}, $K$ is a number field of degree $d$, $R=\cO_K$ and $\Omega$ is a finite alphabet; neither (H) nor (Ind) is used.

\subsection{Gluing along residue classes}\label{ssec:glue}
Eventual periodicity can be checked on the sections of a finite partition into residue classes; this reduces the eventual periodicity of a function to that of suitable kernels.

\begin{lemma}[Gluing]\label{lem:glue}
Let $R=\bigsqcup_{j=1}^J(t_j+\eta_jR)$ be a finite partition into residue classes, with $t_j,\eta_j\in R$ and $\eta_j\ne0$, and let $F:R\to\Omega$. If each section $F_j(y)=F(t_j+\eta_jy)$ is eventually periodic, then so is $F$.
\end{lemma}
\begin{proof}
For each $j$, let $F_j=Q_j$ outside a finite set $X_j\subseteq R$, where $Q_j:R\to\Omega$ is periodic under an additive subgroup $L_j\subseteq R$ of rank $d$. Define $Q:R\to\Omega$ by $Q(x)=Q_j((x-t_j)/\eta_j)$ for $x\in t_j+\eta_jR$. Then $F=Q$ outside the finite set $\bigcup_j(t_j+\eta_jX_j)$.

The subgroup $L:=\bigcap_j\eta_jL_j$ has finite index in $R$, as a finite intersection of full-rank subgroups. For $l\in L\subseteq\eta_jR$, translation by $l$ preserves each class $t_j+\eta_jR$ and shifts the argument of $Q_j$ by $l/\eta_j\in L_j$. Hence $Q$ is $L$-periodic.
\end{proof}

\subsection{Cyclic kernels}\label{ssec:cyclic}
Every automatic function admits a finite partition of its domain into residue classes on which the associated kernels have exact affine symmetries.

Let $\beta\in R$ be expansive and let $E\subseteq R$ be a complete set of representatives of $R/\beta R$ containing $0$ (a window set by Lemma~\ref{lem:windows}). Let $f:R\to\Omega$ be $\beta$-automatic. Its \emph{$\beta$-kernel graph} has as vertices the kernels $y\mapsto f(c+\beta^my)$, with $m\in\N$ and $c\in E^{(m)}$, and an edge from $G$ to $y\mapsto G(b+\beta y)$ labelled by $b$ for each $b\in E$. This graph is finite. We call a kernel \emph{cyclic} if it lies on a closed path of length $p\ge1$. Reading the labels $b_1,\dots,b_p$ along such a path starting and ending at $G$ gives the \emph{exact symmetry}
\begin{equation}\label{eq:cyclesym}
G(e+\beta^py)=G(y)\quad(y\in R),\qquad e=\sum_{i=1}^p\beta^{i-1}b_i.
\end{equation}

\begin{lemma}[Cyclic kernels cover $R$]\label{lem:cyclic}
If $f$ is $\beta$-automatic, there is a finite partition of $R$ into classes $t_j+\beta^{m_j}R$, with $m_j\in\N$ and $t_j\in E^{(m_j)}$, such that each kernel $y\mapsto f(t_j+\beta^{m_j}y)$ is cyclic.
\end{lemma}
\begin{proof}
Label each word over $E$, of length $m$ and value $c$, by the kernel $y\mapsto f(c+\beta^my)$, and stop each branch of the tree of words as soon as its label is cyclic. Before stopping, a branch cannot repeat a label: a repetition would give a closed path through the repeated label, which would then be cyclic. Hence every branch stops at depth less than the number of kernels, and the stopped tree is finite.

For $x\in R$ and $m\in\N$, there is a unique word $w_m$ of length $m$ whose value is congruent to $x$ modulo $\beta^mR$ (Lemma~\ref{lem:windows}\,(ii)), and $w_m$ is the beginning of $w_{m'}$ for $m'\ge m$. These words follow a single branch and meet exactly one leaf, so the classes of the leaves form the required partition.
\end{proof}

\begin{lemma}[Reconstruction]\label{lem:reconstruct}
Let $f$ be $\beta$-automatic. \textup{(i)} If $\lambda\in R$ is expansive and every cyclic $\beta$-kernel of $f$ is $\lambda$-automatic, then $f$ is $\lambda$-automatic. \textup{(ii)} If every cyclic $\beta$-kernel of $f$ is eventually periodic, then $f$ is eventually periodic.
\end{lemma}
\begin{proof}
Use the partition of Lemma~\ref{lem:cyclic}, with kernels $G_j(y)=f(t_j+\beta^{m_j}y)$. For (i), Lemma~\ref{lem:insert} gives $\lambda$-automatic functions $\hat G_j:R\to\Omega\cup\{\bot\}$, where $\bot\notin\Omega$, taking the value $\bot$ outside their respective classes. By the partition, exactly one entry of $(\hat G_1(x),\dots,\hat G_J(x))$ is different from $\bot$ at each $x\in R$, and it equals $f(x)$. Thus $f=\Psi(\hat G_1,\dots,\hat G_J)$, where $\Psi$ returns the unique entry different from $\bot$ when there is exactly one, and is arbitrary otherwise; Corollary~\ref{cor:closure} gives (i). Part (ii) is Lemma~\ref{lem:glue}.
\end{proof}

\subsection{Complementary base and global transfer}\label{ssec:complement}
A symmetry whose multiplier is a power of $\delta$ allows us to pass from automaticity in base $\delta\gamma$ to automaticity in base $\gamma$; with cyclic kernels, this gives automaticity in expansive quotients of bases.

\begin{lemma}[Complementary base]\label{lem:complement}
Let $\delta,\gamma\in R$ be expansive and let $f:R\to\Omega$ be $\delta\gamma$-automatic. Suppose that $f(e+\delta^px)=f(x)$ for every $x\in R$, with $e\in R$ and an integer $p\ge1$. Then $f$ is $\gamma$-automatic. Neither $\delta$-automaticity nor coprimality of $\delta$ and $\gamma$ is assumed.
\end{lemma}
\begin{proof}
Put $\lambda=\delta\gamma$, $q=e/(1-\delta^p)\in K$ and, for $\ell\in\N$,
\[
 a_\ell=q(1-\delta^{p\ell})=e\sum_{0\le i<\ell}\delta^{pi}\in R;
\]
iterating the symmetry gives $f(a_\ell+\delta^{p\ell}x)=f(x)$ for $x\in R$. Fix window sets $D_\gamma,D_\lambda$ with size constants $C_\gamma,C_\lambda$ (Lemma~\ref{lem:fixedlength}). Let $n\in\N$ and $w\in D_\gamma^{(n)}$, and put $\ell=\lceil n/p\rceil$ and $s=p\ell-n\in[0,p)$. Since $\delta^{p\ell}\gamma^n=\delta^s\lambda^n$,
\[
f(w+\gamma^ny)=f(a_\ell+\delta^{p\ell}w+\delta^s\lambda^ny)\qquad(y\in R).
\]
Writing $a_\ell+\delta^{p\ell}w=c+\lambda^nv$ with $c\in D_\lambda^{(n)}$ and $v\in R$, this $\gamma$-kernel is $y\mapsto Z(v+\delta^sy)$, where the $\lambda$-kernel $Z$ is given by $Z(z)=f(c+\lambda^nz)$. Using $\delta^{p\ell}=\delta^s\lambda^n/\gamma^n$ and the bounds $|\sigma w|\le C_\gamma|\sigma\gamma|^n$, $|\sigma c|\le C_\lambda|\sigma\lambda|^n$, we obtain at each archimedean place
\begin{align*}
 |\sigma v|&\le |\sigma q|\bigl(|\sigma\lambda|^{-n}+|\sigma\delta|^s|\sigma\gamma|^{-n}\bigr)
       +C_\gamma|\sigma\delta|^s+C_\lambda\\
 &\le |\sigma q|(1+L_\sigma)+C_\gamma L_\sigma+C_\lambda,
 \qquad L_\sigma:=|\sigma\delta|^{p-1}.
\end{align*}
Thus $v$ lies in a finite set (Lemma~\ref{lem:proper}); so do $Z$ and $s$, and Corollary~\ref{cor:windowcrit} shows that $f$ is $\gamma$-automatic.
\end{proof}
The centre $q$ need not lie in $R$; only the iterates $x\mapsto a_\ell+\delta^{p\ell}x$, which map $R$ into itself, are used.

\begin{theorem}[Global transfer]\label{thm:transfer}
Let $\delta,\gamma\in R$ be expansive and let $f:R\to\Omega$. If $f$ is both $\delta\gamma$-automatic and $\delta$-automatic, then it is $\gamma$-automatic.
\end{theorem}
\begin{proof}
Fix a complete set of representatives of $R/\delta R$ containing $0$. A cyclic $\delta$-kernel $G$ of $f$ with respect to this set satisfies a symmetry of the form \eqref{eq:cyclesym} in base $\delta$. As a section of $f$, it is $\delta\gamma$-automatic (Proposition~\ref{prop:affine}), and hence $\gamma$-automatic by Lemma~\ref{lem:complement}. We conclude by Lemma~\ref{lem:reconstruct}\,(i).
\end{proof}
Together with Lemma~\ref{lem:powers}, this shows that the set of expansive elements in which $f$ is automatic is closed under positive powers, roots in $R$ and expansive quotients (when the quotient belongs to $R$).

\subsection{A second symmetry on a section}\label{ssec:sections}
When a factor of an expansive base is not expansive, Theorem~\ref{thm:transfer} does not apply. The following lemma still gives finiteness up to translation, from which a second exact symmetry can be obtained.

\begin{lemma}[Division of a symmetry]\label{lem:division}
Let $\delta,\gamma\in R\setminus\{0\}$ be such that $\lambda=\delta\gamma$ is expansive, and let $G:R\to\Omega$ be $\lambda$-automatic. Suppose that $G(e+\delta^py)=G(y)$ for every $y\in R$, with $e\in R$ and an integer $p\ge1$. Put $\mu=\gamma^p$ and fix a window set $D$ for $\lambda$. Every section $y\mapsto G(r+\mu^ky)$, with $r\in R$ and $k\in\N$, is a translate $y\mapsto h(t+y)$, with $t\in R$, of a $\lambda$-kernel $h\in\cK_{\lambda,D}(G)$. In particular, up to translation these sections form a finite family.
\end{lemma}
\begin{proof}
For $k\in\N$, put $a_k=e\sum_{0\le i<k}\delta^{pi}\in R$; iterating the symmetry gives $G(a_k+\delta^{pk}y)=G(y)$. Given $r\in R$ and $k\in\N$, choose $c\in D^{(pk)}$ with $c\equiv a_k+\delta^{pk}r\pmod{\lambda^{pk}R}$. Since $\lambda^{pk}=\delta^{pk}\mu^k$, the element $r':=(c-a_k)/\delta^{pk}$ lies in $R$ and $r-r'\in\mu^kR$. The symmetry gives
\[
 G(r'+\mu^ky)=G(a_k+\delta^{pk}(r'+\mu^ky))=G(c+\lambda^{pk}y),
\]
a $\lambda$-kernel $h$ of $G$. With $t=(r-r')/\mu^k\in R$, we get $G(r+\mu^ky)=h(t+y)$. Finally, $\cK_{\lambda,D}(G)$ is finite because $G$ is $\lambda$-automatic.
\end{proof}
Unlike the carry in Lemma~\ref{lem:complement}, the translation $t$ need not be bounded: the sections are finite in number only \emph{up to translation}, which suffices for symmetries but not for automaticity.

\begin{lemma}[Repetition along a chain]\label{lem:repetition}
Let $a,b\in R\setminus\{0\}$ be such that $\lambda=ab$ is expansive, and let $\Xi:R\to\Omega$ be $\lambda$-automatic. Suppose that $\Xi(s_0+ax)=\Xi(x)$ for every $x\in R$, with $s_0\in R$. Consider a sequence of sections $F_j:R\to\Omega$ ($j\in\N$) defined by
\[
 F_0(y)=\Xi(r_0+b^{i_0}y),\qquad F_{j+1}(y)=F_j(\varepsilon_j+by),
\]
where $r_0\in R$, $i_0\in\N$ and $\varepsilon_j\in R$ for each $j\in\N$. There exist $i,j\in\N$ with $i<j$ and $w\in R$ such that $F_i(w+b^{j-i}y)=F_i(y)$ for every $y\in R$.
\end{lemma}
\begin{proof}
Each $F_j$ is a section $y\mapsto\Xi(r^{(j)}+b^{i_0+j}y)$ with $r^{(j)}\in R$. By Lemma~\ref{lem:division} with $p=1$, $F_j(y)=h_j(t_j+y)$, where $t_j\in R$ and $h_j$ lies in the finite set $\cK_{\lambda,D}(\Xi)$ for a fixed window set $D$ for $\lambda$. Choose $i<j$ with $h_i=h_j$; then $F_j(y)=F_i(y+t)$ with $t=t_j-t_i$.

On the other hand, the chain gives $F_j(y)=F_i(c+b^\ell y)$ with $\ell=j-i$ and $c=\sum_{i\le m<j}b^{m-i}\varepsilon_m$. Hence $F_i(c+b^\ell y)=F_i(y+t)$, and substituting $z=y+t$ gives $F_i(c-b^\ell t+b^\ell z)=F_i(z)$, with translation term $w=c-b^\ell t\in R$.
\end{proof}
\section{Equal supports}\label{sec:equal}

It remains to treat bases contained in the same prime ideals, as in the running example $(14,14u)$ of the Introduction. The method of Sections~\ref{sec:orbits}--\ref{sec:different} then fails: for all $n\in\N$ and $t\ge1$ we have $\beta^kR\subseteq\alpha^{tn}R$ for large $k$, so the image of $\beta^kR$ in $R/\alpha^{tn}R$ is trivial and Proposition~\ref{prop:factorbound} gives no contradiction; moreover, the factor $\mathcal A$ of the completion is absent. Instead, we use exact symmetries of kernels (Section~\ref{sec:transfer}) and compare two symmetries whose multipliers behave oppositely at an archimedean place.

Throughout this section, $K$ is a number field of degree $d\ge2$, $R=\cO_K$, and the expansive bases $\alpha,\beta\in R$ satisfy (H) and (Ind). We write $S$ for their common support and take $f:R\to\Omega$, with $\Omega$ a finite alphabet.

\emph{Strategy.} The pairs $(r,s)\in\Z^2$ with $\alpha^r\beta^s\in R$ are the integer points of a cone, whose extremal rays give elements $a_0,b_0\in R$ (\S\ref{ssec:cone}). Neither has a conjugate of absolute value $1$, so each is either expansive or contracting at some archimedean place.
\begin{enumerate}[label=\textup{(\arabic*)},nosep]
\item \emph{An expansive ray element.} Global transfer makes $f$ automatic in this base, whose support differs from that of $\alpha$, and Section~\ref{sec:different} applies.
\item \emph{A ray element $a_0$ contracting at $\sigma_0$.} A function which is not eventually periodic yields another such function with two exact symmetries, whose multipliers are powers of $\alpha$ and $a_0$. Their centres coincide, because suitable powers of the multipliers are complementary and one contracts at $\sigma_0$. After recentring, we obtain a function $F$ invariant under $x\mapsto ax$ and $x\mapsto bx$, where $ab$ is a power of $\alpha$. We prove exact periodicity of certain deep sections and deduce that $F$ is eventually periodic, a contradiction. The steps are Lemma~\ref{lem:twosym} (two symmetries), Lemma~\ref{lem:recentre} (recentring), Lemmas~\ref{lem:deepchildren} and~\ref{lem:deepsec} (deep sections), Lemma~\ref{lem:ending} (conclusion) and Proposition~\ref{prop:contracting}.
\end{enumerate}

\subsection{The integer cone and its ray elements}\label{ssec:cone}
For $Q\in S$, put $h_Q=v_Q(\beta)/v_Q(\alpha)$. Let $h_{\min}$ and $h_{\max}$ be the least and greatest of these rational numbers. They are distinct: if all $h_Q$ were equal to $m/n$, with $m,n\ge1$, the principal ideals $\beta^nR$ and $\alpha^mR$ would coincide, contrary to (Ind). For $(r,s)\in\Z^2$, the monomial $\alpha^r\beta^s$ belongs to $R$ if and only if $r+sh_Q\ge0$ for all $Q\in S$. Write $h_{\min}=m_0/n_0$ and $h_{\max}=m_1/n_1$ in lowest terms, with positive numerators and denominators, and set
\[
 a_0=\frac{\beta^{n_0}}{\alpha^{m_0}},\qquad
 b_0=\frac{\alpha^{m_1}}{\beta^{n_1}}.
\]
These are the elements associated with the primitive integer vectors $(-m_0,n_0)$ and $(m_1,-n_1)$ on the two extremal rays.

\begin{lemma}[Ray elements]\label{lem:rays}
The elements $a_0,b_0$ belong to $R$ and are not units. For $Q\in S$, we have $v_Q(a_0)=0$ if and only if $h_Q=h_{\min}$, and $v_Q(b_0)=0$ if and only if $h_Q=h_{\max}$. Moreover,
\[
 a_0^{n_1}b_0^{n_0}=\alpha^\kappa,\qquad
 \kappa=n_0m_1-m_0n_1\ge1.
\]
\end{lemma}
\begin{proof}
For $Q\in S$,
\[
 v_Q(a_0)=v_Q(\alpha)(n_0h_Q-m_0)\ge0,\qquad
 v_Q(b_0)=v_Q(\alpha)(m_1-n_1h_Q)\ge0.
\]
Outside $S$, both valuations are zero. Thus both elements are integral, and the strict inequality $h_{\min}<h_{\max}$ ensures that neither is a unit. The identity follows by multiplying the two monomials.
\end{proof}

\begin{lemma}[Absolute value one]\label{lem:modone}
A non-zero algebraic integer with a conjugate of absolute value $1$ is a unit. Consequently each ray element is either expansive or contracting at some archimedean place.
\end{lemma}
\begin{proof}
Let $P\in\Z[X]$ be the monic minimal polynomial of the algebraic integer, of degree $n$ and constant term $c_0$, and let $z$ be a root with $|z|=1$. Then $1/z=\bar z$ is also a root of $P$, so $z$ is a root of the reciprocal polynomial $P^*(X)=X^nP(1/X)$, whose constant term is $1$. By irreducibility, $P$ divides $P^*$ in $\Q[X]$. Comparing leading and constant coefficients gives $P^*=c_0P$ and $c_0^2=1$. Thus the algebraic integer is a unit. The assertion about the ray elements follows from Lemma~\ref{lem:rays}.
\end{proof}

\begin{example}
Return to $R=\Z[\sqrt2]$, $\alpha=14$, $\beta=14u$, where $u=8+5\sqrt2$. With $P=(3+\sqrt2)R$ and $\bar P=(3-\sqrt2)R$, we have $S=\{\sqrt2R,P,\bar P\}$. The values of $h_Q$ at these three primes are respectively $3/2,2,1$. The ray elements are $a_0=u$ and $b_0=\bar u$, and $\kappa=1$, since $a_0b_0=u\bar u=14=\alpha$. At the two real places $\sigma_+,\sigma_-$ (sending $\sqrt2$ to $\sqrt2$ and $-\sqrt2$, respectively), the absolute values of $u$ are approximately $15.07$ and $0.929$, and those of $\bar u$ are interchanged. Each ray element therefore contracts at one place.
\end{example}

\subsection{An expansive ray element}\label{ssec:exray}
\begin{proposition}[Expansive ray element]\label{prop:exray}
If either $a_0$ or $b_0$ is expansive, every doubly automatic function is eventually periodic.
\end{proposition}
\begin{proof}
Suppose that $a_0$ is expansive (otherwise interchange $\alpha$ and $\beta$, which interchanges the ray elements). The function $f$ is automatic in the bases $\alpha^{m_0}$ and $\beta^{n_0}=\alpha^{m_0}a_0$. By global transfer (Theorem~\ref{thm:transfer}), it is $a_0$-automatic. The pair $(\alpha,a_0)$ has different supports, since a prime $Q$ with $h_Q=h_{\min}$ contains $\alpha$ but not $a_0$. It also satisfies (H): for $r,s\ge1$, the field $\Q(\alpha^r,a_0^s)$ contains $\alpha^r$ and
\[
 (\alpha^r)^{m_0s}(a_0^s)^r=\beta^{n_0rs}.
\]
By (H) for the original pair this field is $K$. We conclude by Theorem~\ref{thm:different}.
\end{proof}

\subsection{Exact symmetries and their centres}\label{ssec:returns}
We collect the properties of exact symmetries used below: inheritance of eventual periodicity, a bound on the carries of exact returns, equality of centres for complementary multipliers, and a criterion for passing from eventual to exact periodicity.

\begin{lemma}[Two inheritance properties]\label{lem:hereditary}
\textup{(i)} If $G:R\to\Omega$ is doubly automatic and is not eventually periodic, then $\Per(G)=\{0\}$. \textup{(ii)} If $G$ is eventually periodic, so is every affine section $y\mapsto G(t+\eta y)$, with $t,\eta\in R$ and $\eta\ne0$. Thus if a section is not eventually periodic, neither is $G$.
\end{lemma}
\begin{proof}
For (i), a non-zero exact period is a period at infinity, so Theorem~\ref{thm:onefar} applies. For (ii), suppose that $G$ agrees outside a finite set $X_0\subseteq R$ with a function $g$ periodic under a full-rank subgroup $L\subseteq R$. Choose an integer $m\ge1$ with $mR\subseteq L$. Since $\eta\ne0$, the section of $G$ agrees with $y\mapsto g(t+\eta y)$ outside a finite set. The latter is periodic under the subgroup $\{y\in R:\eta y\in L\}\supseteq mR$.
\end{proof}

\begin{lemma}[Carries of exact symmetries]\label{lem:carries}
Let $\lambda\in R$ be expansive and let $G:R\to\Omega$ be $\lambda$-automatic, with $\Per(G)=\{0\}$. Suppose that $d_n\in R$ for $n$ in a non-empty set $\mathcal N\subseteq\N$ and that
\[
 G(d_n+\lambda^ny)=G(y)\qquad(y\in R,\ n\in\mathcal N).
\]
Then, at every archimedean place $\sigma$, we have $\sup_{n\in\mathcal N}|\sigma d_n|/|\sigma\lambda|^n<\infty$.
\end{lemma}
\begin{proof}
Fix a window set $D$ for $\lambda$. Write $d_n=c_n+\lambda^nq_n$, where $c_n\in D^{(n)}$ and $q_n\in R$. The symmetry gives $G(c_n+\lambda^ny)=G(y-q_n)$. The functions on the left are kernels and form a finite family. Since $\Per(G)=\{0\}$, distinct translations of $G$ are distinct functions, so the $q_n$ take only finitely many values. The size estimate for $c_n$ (Lemma~\ref{lem:fixedlength}) now gives the bound.
\end{proof}

\begin{lemma}[Complementary multipliers]\label{lem:compcentres}
Let $\lambda\in R$ be expansive and let $G:R\to\Omega$ be $\lambda$-automatic, with $\Per(G)=\{0\}$. Suppose that the maps
\[
 \varphi(y)=z_1+A(y-z_1),\qquad
 \theta(y)=z_2+B(y-z_2)
\]
with $z_1,z_2\in K$ and $A,B\in R$, send $R$ into $R$ and satisfy $G\circ\varphi=G\circ\theta=G$. If $AB=\lambda$ and $|\sigma_0 A|<1$ at some archimedean place $\sigma_0$, then $z_1=z_2$.
\end{lemma}
\begin{proof}
For each $n\ge1$, the symmetries $\varphi^n\circ\theta^n$ and $\theta^n\circ\varphi^n$ both have multiplier $\lambda^n$, with translation terms in $R$ given by
\begin{align*}
 d_n&=(1-A^n)z_1+A^n(1-B^n)z_2,\\
 d'_n&=(1-B^n)z_2+B^n(1-A^n)z_1.
\end{align*}
Their difference is $(1-A^n)(1-B^n)(z_1-z_2)$. By Lemma~\ref{lem:carries}, the quotient $(d_n-d'_n)/\lambda^n$, namely
\[
 (A^{-n}-1)(B^{-n}-1)(z_1-z_2),
\]
is bounded at $\sigma_0$. But $|\sigma_0 A^{-n}|\to\infty$ and $|\sigma_0 B^{-n}|\to0$, since $|\sigma_0 B|=|\sigma_0\lambda|/|\sigma_0 A|>1$. It follows that $|\sigma_0(z_1-z_2)|=0$, hence $z_1=z_2$.
\end{proof}
In the running example, with $A=a_0$ and $B=\alpha=14$, the scalar factor has absolute value asymptotic to $(8-5\sqrt2)^{-n}$ at $\sigma_-$.

\begin{lemma}[A symmetry centred outside $R$]\label{lem:offcentre}
Let $G:R\to\Omega$ be eventually periodic and suppose that $G(e+cy)=G(y)$ for every $y\in R$, with $e,c\in R$ and $|\sigma c|>1$ at some archimedean place $\sigma$. If the centre $z=e/(1-c)$ does not belong to $R$, then $G$ is periodic under a full-rank subgroup, without exceptions.
\end{lemma}
\begin{proof}
Write $G=G'$ outside a finite set $X_0\subseteq R$, where $G'$ is periodic under a full-rank subgroup $L\supseteq mR$ for some integer $m\ge1$, and put $\psi(y)=e+cy$. The set on which $G'\circ\psi\ne G'$ is contained in $X_0\cup\psi^{-1}(X_0)$, which is finite, and is periodic under
\[
 L\cap\{y\in R:cy\in L\}\supseteq mR.
\]
It is therefore empty. For any $y\in R$, we have $y\ne z$ and $\psi^k(y)=z+c^k(y-z)$. Consequently $|\sigma\psi^k(y)|\to\infty$, and $\psi^k(y)\notin X_0$ for all sufficiently large $k$. It follows that $G(y)=G(\psi^k(y))=G'(\psi^k(y))=G'(y)$, as required.
\end{proof}

\subsection{Two symmetries of a function}\label{ssec:twosym}
Suppose that $|\sigma_0a_0|<1$ at some archimedean place $\sigma_0$ (if it is $b_0$ which contracts, interchange $\alpha$ and $\beta$, thereby interchanging $a_0$ and $b_0$ and leaving $\kappa$ unchanged). The two ray elements never contract at the same place: $|\sigma a_0|^{n_1}|\sigma b_0|^{n_0}=|\sigma\alpha|^\kappa>1$.

\begin{lemma}[Two symmetries]\label{lem:twosym}
If $f$ is doubly automatic and is not eventually periodic, there exists a doubly automatic function $H:R\to\Omega$, also not eventually periodic, such that
\[
 H(e'+\alpha^{p_1}y)=H(y),\qquad H(d'+a_0^{p_2}y)=H(y)\qquad(y\in R),
\]
where $p_1,p_2\ge1$ and $e',d'\in R$.
\end{lemma}
\begin{proof}
\emph{1. An initial symmetry.} By Lemma~\ref{lem:reconstruct}\,(ii) in base $\alpha$, some cyclic $\alpha$-kernel $G$ of $f$ is not eventually periodic. It is doubly automatic (Proposition~\ref{prop:affine}) and satisfies $G(e+\alpha^py)=G(y)$ with $e\in R$ and $p\ge1$. Iterating $m_0$ times gives $G(e_1+\Delta y)=G(y)$, where $e_1\in R$ and $\Delta=\alpha^{m_0p}$.

\emph{2. A chain of sections.} Put $\mu=a_0^p$, so that $\Delta\mu=\beta^{n_0p}=:\lambda_1$ is expansive. The support of $\mu$ is contained in $S$ and $\Delta$ has positive valuation at every prime in $S$, so $\nu:=\Delta^j/\mu\in R$ for some integer $j\ge1$. For $k\in\N$, set
\[
 A_k=e_1\sum_{0\le i<k}\Delta^i,\qquad H_k(z)=G(A_{jk}+\mu^kz)\quad(z\in R).
\]
As $\Delta^{jk}=\mu^k\nu^k$, iterating the initial symmetry gives $G(z)=H_k(\nu^kz)$. Thus $G$ is a section of $H_k$, so $H_k$ is not eventually periodic (Lemma~\ref{lem:hereditary}\,(ii)); and $H_k$, a section of $G$, is doubly automatic. The identity $A_{j(k+1)}=A_{jk}+\Delta^{jk}A_j$ gives
\[
 H_{k+1}(z)=H_k(\nu^k A_j+\mu z),
\]
a chain of the kind considered in Lemma~\ref{lem:repetition}.

\emph{3. A symmetry with multiplier $\Delta$.} For $m\in\N$, the identity $A_{jm}+\Delta^{jm}e_1=e_1+\Delta A_{jm}$ gives
\begin{align*}
 H_m(e_1\nu^m+\Delta z)
 &=G(A_{jm}+\Delta^{jm}e_1+\Delta\mu^m z)\\
 &=G(e_1+\Delta(A_{jm}+\mu^m z))=H_m(z).
\end{align*}

\emph{4. A second symmetry.} Lemma~\ref{lem:repetition}, applied with $\Xi=G$, $a=\Delta$ and $b=\mu$ ($G$ is $\lambda_1$-automatic), gives $k<\ell$ and $w\in R$ with $H_k(w+\mu^{\ell-k}y)=H_k(y)$. Take $H=H_k$, $e'=e_1\nu^k$, $d'=w$, $p_1=m_0p$ and $p_2=p(\ell-k)$. Only the factorisation $\beta^{n_0p}=\alpha^{m_0p}a_0^p$ has been used, not automaticity in base $a_0$.
\end{proof}

\subsection{Recentring}\label{ssec:recentre}
The two symmetries of Lemma~\ref{lem:twosym} have the same centre; moving this centre to the origin turns them into invariances under multiplication by elements of $R$.

\begin{lemma}[Recentring]\label{lem:recentre}
Let $H$ be as in Lemma~\ref{lem:twosym}.
\begin{enumerate}[label=\textup{(\roman*)},nosep]
\item The two symmetries have the same centre $z\in K$, integral at every prime in $S$.
\item Put $\eta=1-\alpha^{p_1}$ and choose a symbol $\bot\notin\Omega$. The function $F:R\to\Omega\cup\{\bot\}$ defined by
\[
 F(x)=\begin{cases}
 H(z+x/\eta)&\text{if }x\in-\eta z+\eta R,\\
 \bot&\text{otherwise}
 \end{cases}
\]
is doubly automatic, is not eventually periodic and satisfies $F(gx)=F(x)$ for every $x\in R$ and every $g\in R\cap\langle a_0^{p_2},\alpha^{p_1}\rangle$, where the multiplicative group is taken in $K^\times$.
\item There exist positive powers $a,b$ of $a_0,b_0$, respectively, in this group, and an integer $N\ge1$, such that $ab=\alpha^N$ and $F(ax)=F(x)=F(bx)$ for all $x\in R$.
\end{enumerate}
\end{lemma}
\begin{proof}
Write $\varphi(y)=d'+a_0^{p_2}y$ and $\theta(y)=e'+\alpha^{p_1}y$ for the two symmetries of $H$.

(i) The iterates $\varphi^{p_1}$ and $\theta^{m_0p_2}$ have the same centres as $\varphi$ and $\theta$, and multipliers
\[
 A=a_0^{p_1p_2},\qquad B=\alpha^{m_0p_1p_2},
\]
with $AB=\beta^{n_0p_1p_2}$ expansive and $|\sigma_0 A|<1$. As $H$ is automatic in base $AB$ and $\Per(H)=\{0\}$ (Lemma~\ref{lem:hereditary}\,(i)), Lemma~\ref{lem:compcentres} gives equality of the centres. Since $z=e'/(1-\alpha^{p_1})$ and $1-\alpha^{p_1}$ is a unit at every prime in $S$, the centre $z$ is integral at these primes.

(ii) We have $\eta z=e'\in R$, and $\eta R$ is coprime to every prime in $S$. The function $F$ is obtained by inserting $H$ on the class $-\eta z+\eta R$ (Lemma~\ref{lem:insert}), so it is doubly automatic; and $H(y)=F(\eta y-\eta z)$ expresses $H$ as a section of $F$, which is therefore not eventually periodic.

Let $g=a_0^{ip_2}\alpha^{jp_1}$ with $i,j\in\N$. The map $y\mapsto z+g(y-z)$ is the symmetry $\varphi^i\circ\theta^j$ of $H$; in particular $(1-g)z\in R$. Multiplication by $g$ permutes the classes modulo $\eta R$, since the support of $g$ is contained in $S$, and fixes the class of $-\eta z$, since $(1-g)z\in R$. For $x=-\eta z+\eta w$ with $w\in R$, we obtain
\[
 F(gx)=H(z+g(w-z))=H(w)=F(x);
\]
outside this class both values are $\bot$.

For an arbitrary $g\in R$ in the indicated group, choose $h=a_0^{ip_2}\alpha^{jp_1}$ with $i,j\in\N$ large enough that $hg$ is also of this form. Then $F(gx)=F(hgx)=F(x)$ for $x\in R$.

(iii) Set $M=p_1p_2$, $a=a_0^{n_1M}$ and $b=b_0^{n_0M}$. By Lemma~\ref{lem:rays},
\[
 ab=\alpha^{\kappa M},\qquad b=\alpha^{\kappa M}a_0^{-n_1M},
\]
so $a$ and $b$ belong to the group in (ii). Take $N=\kappa M$.
\end{proof}

\emph{Summary.} We now have a doubly automatic function $F$, not eventually periodic, such that
\[
 F(ax)=F(x)=F(bx),\qquad ab=\lambda=\alpha^N,
\]
where $a,b$ are positive powers of the ray elements $a_0,b_0$, with $|\sigma_0 a|<1<|\sigma_0 b|$ at the contracting place $\sigma_0$. The function $F$ is automatic in base $\lambda$, not necessarily in base $a$ or $b$: only the two invariances are used. It remains to show that $F$ is eventually periodic, contradicting Lemma~\ref{lem:recentre}\,(ii): we first establish exact periodicity of certain deep sections $G_{r,i}$ (\S\ref{ssec:deep}), and then use it near the axes in a product of local rings (\S\ref{ssec:ending}).

\subsection{Deep sections are periodic}\label{ssec:deep}
Partition $S$ into the sets $S_a$, $S_b$ and $S_{ab}$ of primes containing only $a$, only $b$, and both, respectively. The set $S_b$ consists of the primes with $h_Q=h_{\min}$, and $S_a$ of those with $h_Q=h_{\max}$ (Lemma~\ref{lem:rays}); both are non-empty, and since $h_{\min}<h_{\max}$ every prime in $S$ contains $a$ or $b$. For $r\in R$ and $i\in\N$, put $G_{r,i}(y)=F(r+b^iy)$ for $y\in R$. We say that $r$ is \emph{deep at level $i$} if
\[
 v_Q(r)\ge i\,v_Q(b)\qquad(Q\in S_{ab}),
\]
with the convention $v_Q(0)=+\infty$. The condition is vacuous if $S_{ab}=\varnothing$. It is not vacuous in the running example of \S\ref{ssec:cone}, where $S_b=\{\bar P\}$, $S_a=\{P\}$ and $S_{ab}=\{\sqrt2R\}$.

\begin{lemma}[Deep children]\label{lem:deepchildren}
Let $r$ be deep at level $i$ and put $G=G_{r,i}$.
\textup{(i)} There exists an integer $T=T(i)\ge1$, independent of $r$, such that $G$ admits the symmetry $\psi(y)=z_i+a^T(y-z_i)$, where $z_i=-r/b^i$ and $\psi(R)\subseteq R$. The integer $T$ may be replaced by any positive multiple.
\textup{(ii)} If $G$ is not eventually periodic, there exists $e^*\in R$ such that $r+b^ie^*$ is deep at level $i+1$ and $y\mapsto G(e^*+by)$ is not eventually periodic.
\end{lemma}
\begin{proof}
(i) Since $a$ is a unit at the primes in $S_b$, there is $T\ge1$ such that $a^T\equiv1$ modulo the part of $b^iR$ supported on $S_b$. Then $(a^T-1)r\in b^iR$: at the primes in $S_b$ by the choice of $T$, at those in $S_{ab}$ by depth, and at all other primes because $b$ is a unit there. Thus $\psi(R)\subseteq R$, and
\[
 G(\psi(y))=F(a^T(r+b^iy))=G(y).
\]

(ii) Replacing $T$ by a multiple, we may also assume that $v_Q(a^T)\ge v_Q(b)$ for $Q\in S_{ab}$. Let $D_b$ be a complete set of representatives of $R/bR$. For $e\in D_b$, choose $e^*\in D_b$ with $e^*\equiv\psi(e)\pmod{bR}$ and put $q=(\psi(e)-e^*)/b\in R$. Then
\[
 G(e+by)=G(\psi(e)+a^Tby)=G(e^*+b(q+a^Ty)),
\]
so the section on the class of $e$ is an affine section of the section on the class of $e^*$. For $Q\in S_{ab}$, the element $z_i$ is integral at $Q$ by depth, and $\psi(e)-z_i=a^T(e-z_i)$; hence $v_Q(e^*-z_i)\ge v_Q(b)$ and $v_Q(r+b^ie^*)\ge(i+1)v_Q(b)$, so $r+b^ie^*$ is deep at level $i+1$. If all the sections on the classes of these $e^*$ were eventually periodic, so would be all the sections on classes modulo $bR$ (Lemma~\ref{lem:hereditary}\,(ii)), and then $G$ (Lemma~\ref{lem:glue}), a contradiction.
\end{proof}

\begin{lemma}[Exact periodicity of deep sections]\label{lem:deepsec}
If $r$ is deep at level $i$ and $r\notin b^iR$, then $G_{r,i}$ is periodic under a full-rank subgroup, without exceptions. The corresponding statement holds with $a$ and $b$ interchanged (using depth defined by the valuations of $a$).
\end{lemma}
\begin{proof}
\emph{1. A chain of deep sections.} Suppose that $G_{r,i}$ is not eventually periodic. Applying Lemma~\ref{lem:deepchildren}\,(ii) repeatedly, we obtain $r_0=r$ and $r_{k+1}=r_k+b^{i+k}e_k$ ($e_k\in R$), with $r_k$ deep at level $i+k$, such that the sections $\Phi_k=G_{r_k,i+k}$ are not eventually periodic. Then $\Phi_{k+1}(y)=\Phi_k(e_k+by)$ and $r_k\equiv r\pmod{b^iR}$.

\emph{2. Two centres.} Lemma~\ref{lem:repetition}, applied with $\Xi=F$ and $s_0=0$, gives $k\in\N$, $\ell\ge1$ and $w\in R$ such that $\Phi=\Phi_k$ satisfies $\Phi(w+b^\ell y)=\Phi(y)$. Its centre $w/(1-b^\ell)$ is integral at the primes in $S_{ab}\cup S_b$, where $1-b^\ell$ is a unit. With $k$ now fixed, Lemma~\ref{lem:deepchildren}\,(i) gives $\Phi$ a symmetry with multiplier $a^{T'}$, where $T'=T(i+k)$, and centre $-r_k/b^{i+k}$.

Put $M'=T'\ell$. Iterating the two symmetries $T'$ and $\ell$ times gives symmetries with multipliers $b^{M'}$ and $a^{M'}$, unchanged centres and product of multipliers $\lambda^{M'}$; both maps send $R$ into $R$. The function $\Phi$ is $\lambda$-automatic as a section of $F$ (Proposition~\ref{prop:affine}), hence $\lambda^{M'}$-automatic (Lemma~\ref{lem:powers}); being doubly automatic and not eventually periodic, it has $\Per(\Phi)=\{0\}$ (Lemma~\ref{lem:hereditary}\,(i)). As $|\sigma_0 a^{M'}|<1$, Lemma~\ref{lem:compcentres} gives
\[
 -r_k/b^{i+k}=w/(1-b^\ell).
\]
The left-hand side is therefore integral at the primes in $S_b$; it is integral at those in $S_{ab}$ by depth, and $b$ is a unit at all other primes. Hence $r_k\in b^{i+k}R$, and $r_k\equiv r\pmod{b^iR}$ gives $r\in b^iR$, a contradiction.

\emph{3. Exact periodicity.} Thus $G_{r,i}$ is eventually periodic. Its symmetry with multiplier $a^{T(i)}$ has centre $-r/b^i\notin R$, and $|\sigma a^{T(i)}|>1$ at some archimedean place $\sigma$, since $a$ is a non-zero non-unit and so $|N_{K/\Q}(a)|\ge2$. Lemma~\ref{lem:offcentre} gives exact periodicity under a full-rank subgroup.

With $a$ and $b$ interchanged, the proof is the same, except that in Step~2 the multiplier contracting at $\sigma_0$ is $a^{M'}$, obtained by repetition rather than from Lemma~\ref{lem:deepchildren}, and that Step~3 uses that $b$ too is a non-unit.
\end{proof}

\subsection{The argument near the axes}\label{ssec:ending}
For $V\subseteq S$, write $\hat R_V=\prod_{Q\in V}\cO_Q$, where $\cO_Q$ is the ring of integers of the completion of $K$ at $Q$. In particular,
\[
 \hat R_S=\hat R_{S_a}\times\hat R_{S_{ab}}\times\hat R_{S_b}.
\]
An empty product here denotes the trivial ring. We shall use the two axes on which only the $S_a$ or the $S_b$ component can be non-zero.

\begin{lemma}[Conclusion near the axes]\label{lem:ending}
The function $F$ of Lemma~\ref{lem:recentre} is eventually periodic.
\end{lemma}
\begin{proof}
Fix a window set $D$ for $\lambda$. Near either axis, the value of $F$ is read off from finitely many periodic deep sections; multiplication by suitable powers of $a$, or of $b$, carries every non-zero element towards one axis without changing the value of $F$.

\emph{1. Periodic kernels near the axes.} Let $n\in\N$ and $c\in D^{(n)}\setminus\lambda^nR$ satisfy $v_Q(c)\ge v_Q(\lambda^n)$ for $Q\in S_a\cup S_{ab}$. Since $a$ is a unit at every other prime, $c=a^nc_1$ with $c_1\in R$. As $v_Q(\lambda)=v_Q(a)+v_Q(b)$, the element $c_1$ is deep at level $n$, and $c_1\notin b^nR$ because $c\notin\lambda^nR$. By invariance under multiplication by $a$,
\[
 F(c+\lambda^nz)=F(c_1+b^nz)=G_{c_1,n}(z),
\]
which is periodic under a full-rank subgroup by Lemma~\ref{lem:deepsec}. The same holds, with $a$ and $b$ interchanged, when the valuation condition is imposed at the primes in $S_{ab}\cup S_b$.

\emph{2. A common period subgroup.} Since $F$ is $\lambda$-automatic, these kernels form a finite family, so they are all periodic under $mR$ for some integer $m\ge1$. Write $mR=I_SI_{\rm out}$, where $I_S$ is supported on $S$ and $I_{\rm out}$ is coprime to every prime in $S$, and put $\mathcal Z=\hat R_S\times R/I_{\rm out}$, in which the diagonal image of $R$ is dense by the Chinese remainder theorem.

\emph{3. Values near the $S_b$-axis.} Let $y\in\hat R_{S_b}\setminus\{0\}$ and $\omega\in R/I_{\rm out}$, and choose $n$ with $y\notin b^n\hat R_{S_b}$, which is possible since these ideals have intersection $0$. Let $U\subseteq\mathcal Z$ be the open set of points whose first three components are congruent to $(0,0,y)$ modulo $\lambda^nI_S\hat R_S$ and whose last component is $\omega$.

We show that $F$ is constant on $R\cap U$. Take $x,x'\in R\cap U$ and $c\in D^{(n)}$ congruent to $x$ modulo $\lambda^nR$. Then $c$ satisfies the valuation condition of Step~1 at the primes in $S_a\cup S_{ab}$, and $c\notin\lambda^nR$, because the $S_b$-component of $x$ is congruent to $y$ modulo $b^n\hat R_{S_b}=\lambda^n\hat R_{S_b}$. Moreover $x-x'\in\lambda^nI_S\cap I_{\rm out}=\lambda^nmR$ by coprimality, so $x'$ lies in the class of $x$ modulo $\lambda^nR$ and the two arguments of the corresponding kernel differ by an element of $mR$. By Step~2, $F(x)=F(x')$.

Let $V_b(y,\omega)$ be this common value. It is well defined, that is, independent of the choice of $n$: two such neighbourhoods of the same axis point meet in an open set, which contains a point of $R$ by density. Moreover, $y\mapsto V_b(y,\omega)$ is locally constant on $\hat R_{S_b}\setminus\{0\}$: axis points sufficiently close to $(0,0,y,\omega)$ lie in $U$, and their neighbourhoods meet $U$ in an open set containing a point of $R$. In the same way we obtain a locally constant function $y\mapsto V_a(y,\omega)$ on $\hat R_{S_a}\setminus\{0\}$.

\emph{4. Reading along either axis.} Let $x\in R\setminus\{0\}$ and let $\omega$ be its class modulo $I_{\rm out}$. Choose $n_j\to\infty$ in $\N$ such that $n_j$ is a multiple of the order of $a$ in the unit groups modulo $Q^j$ for all $Q\in S_b$, and modulo $I_{\rm out}$ (the latter condition is vacuous if $I_{\rm out}=R$). Since $a$ is a unit at all these primes and is divisible by the primes in $S_a\cup S_{ab}$,
\[
 a^{n_j}x\longrightarrow(0,0,x_{S_b},\omega)\qquad\text{in }\mathcal Z,
\]
where the image $x_{S_b}$ of $x$ in $\hat R_{S_b}$ is non-zero, since $S_b$ is non-empty and $R$ embeds in each local ring. For large $j$, the point $a^{n_j}x$ lies in a neighbourhood as in Step~3, so that
\[
 F(x)=F(a^{n_j}x)=V_b(x_{S_b},\omega).
\] Interchanging $a$ and $b$ gives $F(x)=V_a(x_{S_a},\omega)$.

\emph{5. Conclusion.} Let $x,x'\in R\setminus\{0\}$ have the same class $\omega$ modulo $I_{\rm out}$. By density, choose $y\in R$ of class $\omega$ whose $S_b$ and $S_a$ components lie in neighbourhoods of $x_{S_b}$ and $x'_{S_a}$, not containing $0$, on which, for this $\omega$, $V_b$ and $V_a$ respectively are constant; then $y\ne0$. Step~4, applied to $x$, $y$ and $x'$, gives
\[
 F(x)=V_b(x_{S_b},\omega)=F(y)=V_a(x'_{S_a},\omega)=F(x').
\]
Thus on $R\setminus\{0\}$ (Step~4 requires $x\ne0$) the value of $F$ depends only on the class modulo $I_{\rm out}$. Since every class has a non-zero representative, the resulting function on $R/I_{\rm out}$ lifts to an $I_{\rm out}$-periodic function on $R$ agreeing with $F$ except possibly at $0$; hence $F$ is eventually periodic.
\end{proof}

\begin{proposition}[Contracting ray element]\label{prop:contracting}
If a ray element contracts at an archimedean place, every doubly automatic function is eventually periodic.
\end{proposition}
\begin{proof}
If $f$ is not eventually periodic, Lemmas~\ref{lem:twosym} and~\ref{lem:recentre} give a function $F$ which is not eventually periodic, in contradiction with Lemma~\ref{lem:ending}.
\end{proof}

\subsection{The theorem for equal supports}\label{ssec:equalmain}
\begin{theorem}[Equal supports]\label{thm:equal}
Let $K$ be a number field of degree $d\ge2$, let $R=\cO_K$ and let $\alpha,\beta\in R$ be expansive elements satisfying \textup{(H)} and \textup{(Ind)}. Suppose that $\alpha$ and $\beta$ have the same support. For a finite alphabet $\Omega$, a function $f:R\to\Omega$ is both $\alpha$- and $\beta$-automatic if and only if it is eventually periodic.
\end{theorem}
\begin{proof}
The converse is Proposition~\ref{prop:converse}. If one of the ray elements is expansive, apply Proposition~\ref{prop:exray}. Otherwise, by Lemma~\ref{lem:modone}, each ray element contracts at some archimedean place, and Proposition~\ref{prop:contracting} applies.
\end{proof}

\begin{remark}[Use of (Ind)]
Condition (Ind) enters through the strict inequality $h_{\min}<h_{\max}$. It ensures that the ray elements are non-units (Lemma~\ref{lem:rays}, used in Lemma~\ref{lem:modone}), that the centres are well defined ($a_0^{p_2}\ne1$, $b^\ell\ne1$), that $|N_{K/\Q}(a)|\ge2$, that $\kappa\ge1$ (and hence that $\lambda=\alpha^N$ is expansive), and that $S_a,S_b$ are non-empty. Without (Ind) the conclusion is false (Proposition~\ref{prop:idealdep}); a particularly simple example is the following. For an expansive $\alpha$ such that $\Q(\alpha^r)=K$ for all $r\ge1$, the pair $(\alpha,\alpha^2)$ satisfies (H) and has equal supports. The set $\{\alpha^n:n\in\N\}$ is $\alpha$-automatic. Indeed, take an admissible digit set containing $0,1$. The powers are represented by the words of the regular language $0^*10^*$, read least significant digit first. Using the equal-value automaton of the proof of Theorem~\ref{thm:sparse}, the language of words having the same value as a word of $0^*10^*$ of the same length is regular. Taking the right quotient of this language by $0^*$ removes the padding when the matching word is longer. We thus obtain the regular language of all expansions of the powers. Thus the set is also $\alpha^2$-automatic by Lemma~\ref{lem:powers}. It has no period at infinity: for $m\ne n$ and any archimedean place $\sigma$,
\[
 |\sigma(\alpha^m-\alpha^n)|\ge|\sigma\alpha|^{\max(m,n)-1}(|\sigma\alpha|-1),
\]
so any fixed non-zero difference can occur for only finitely many pairs $(m,n)$. In particular a non-zero translation cannot preserve this set outside a finite set.
\end{remark}
\section{Proof of the main theorem}\label{sec:main}

\subsection{The theorem}\label{ssec:mainthm}
\begin{theorem}[Main theorem]\label{thm:main}
Let $K$ be a number field of degree $d\ge2$, let $R=\cO_K$ and let $\alpha,\beta\in R$ be expansive. Assume that
\begin{enumerate}[label=\textup{(\roman*)},nosep]
\item \textup{(H)}: $\Q(\alpha^r,\beta^s)=K$ for all $r,s\ge1$;
\item \textup{(Ind)}: $\alpha^rR\ne\beta^sR$ for all $r,s\ge1$.
\end{enumerate}
Then, for a finite alphabet $\Omega$, a function $f:R\to\Omega$ is both $\alpha$- and $\beta$-automatic if and only if it is eventually periodic.
\end{theorem}
\begin{proof}
Eventually periodic functions are automatic in every expansive base (Proposition~\ref{prop:converse}). Conversely, if a prime ideal contains one base but not the other, apply Theorem~\ref{thm:different}. Otherwise the supports are equal, and Theorem~\ref{thm:equal} applies.
\end{proof}

\emph{The role of the hypotheses.} Condition (H) enters through the principle of periods at infinity (Lemma~\ref{lem:good}), which is used in Lemmas~\ref{lem:bdstab} and~\ref{lem:hereditary}, Theorem~\ref{thm:reduction} and Proposition~\ref{prop:stability}. It also enters directly in Theorem~\ref{thm:ranks} and Proposition~\ref{prop:centres}. Condition (Ind) is needed only in the case of equal supports (it follows from the existence of an excess prime in the other case). In both cases the hypothesis $d\ge2$ is used to pass from a full-rank group of periods at infinity to eventual periodicity (Proposition~\ref{prop:rank2}).

\subsection{Periods at infinity in every degree}\label{ssec:degreeone}
For $d=1$ the conclusion fails: $\mathbf 1_{\Z_{\ge0}}$ is automatic in every expansive base but is not eventually periodic in our sense (Remark~\ref{rem:degreeone}). What remains true is that the group of periods at infinity has full rank. In degree one this follows from Cobham's theorem \cite{Cobham}, stated at the beginning of the Introduction and used here without proof.

\begin{corollary}[Periods at infinity in every degree]\label{cor:farperiods}
Let $K$ be a number field of degree $d\ge1$, let $R=\cO_K$ and let $\alpha,\beta\in R$ be expansive and satisfy (H) and (Ind). Then $f:R\to\Omega$, with $\Omega$ finite, is both $\alpha$- and $\beta$-automatic if and only if $\Perf(f)$ has rank $d$.
\end{corollary}
\begin{proof}
Let $d\ge2$. If $f$ is doubly automatic, it is eventually periodic by Theorem~\ref{thm:main}, and an eventually periodic function has a group of periods at infinity of rank $d$ (\S\ref{ssec:periods}). Conversely, if $\Perf(f)$ has rank $d$, then $f$ is eventually periodic by Proposition~\ref{prop:rank2}, hence automatic in every expansive base (Proposition~\ref{prop:converse}). For $d=1$, we have $R=\Z$, and (Ind) means that $|\alpha|$ and $|\beta|$ are multiplicatively independent. By Lemma~\ref{lem:powers}, we may replace the bases by $k=\alpha^2$ and $l=\beta^2$. The digit sets $\{0,\ldots,k-1\}$ and $\{0,-1,\ldots,-(k-1)\}$ are window sets for $k$ (Lemma~\ref{lem:windows}\,(ii)), and similarly for $l$. By Proposition~\ref{prop:indep} and the kernel criterion, both $f|_{\N}$ and the sequence $n\mapsto f(-n)$ on $\N$ are $k$- and $l$-automatic in the classical sense. Applying Cobham's theorem to their level sets, and taking a common period $m\ge1$ for these finitely many sets, we obtain $m\in\Perf(f)$.

Conversely, suppose that $m\Z\subseteq\Perf(f)$ for some integer $m\ge1$. There is an integer $N\ge0$ such that outside $[-N,N]$ the function $f$ agrees with a function taking the values of an $m$-periodic function $g_1:\Z\to\Omega$ on $\Z_{\ge0}$ and of an $m$-periodic function $g_2:\Z\to\Omega$ on $\Z_{<0}$. Periodic functions and characteristic functions of singletons are automatic in every expansive base by Proposition~\ref{prop:converse}. The same holds for $\mathbf 1_{\Z_{\ge0}}$: for a fixed window set $D$ and an expansive base $\gamma\in\Z$, its kernels are characteristic functions of sets of the form $\{y\in\Z:y\ge-c/\gamma^n\}$ or $\{y\in\Z:y\le-c/\gamma^n\}$. Since $|c/\gamma^n|\le C_D$ (Lemma~\ref{lem:fixedlength}), there are only finitely many such sets. Closure under pointwise composition (Corollary~\ref{cor:closure}) now gives automaticity of $f$.
\end{proof}

\subsection{The hypotheses cannot be dropped}\label{ssec:necessity}
Neither (H) nor (Ind) can be omitted: in each case we construct a doubly automatic set that is not eventually periodic.

\begin{proposition}[Dependent ideals]\label{prop:idealdep}
Let $\alpha,\beta\in R$ be expansive. If $\alpha^rR=\beta^sR$ for some $r,s\ge1$, there exists a subset $Y\subseteq R$ which is both $\alpha$- and $\beta$-automatic but has no period at infinity and hence is not eventually periodic.
\end{proposition}
\begin{proof}
Let $\gamma_1=\alpha^r$, $\gamma_2=\beta^s=\gamma_1u$, where $u\in R^\times$, and put $I=\gamma_1R$. For $x\in R\setminus\{0\}$ define $\nu(x)=\max\{n\in\N:x\in I^n\}$, which is finite since $|N_{K/\Q}(\gamma_1)|\ge2$. Set $\nu(0)=+\infty$. If $\nu(x)\ne\nu(z)$, then $\nu(x+z)=\min(\nu(x),\nu(z))$. Also $\nu(\gamma^nz)=n+\nu(z)$ for $\gamma\in\{\gamma_1,\gamma_2\}$. Let $Y=\{x\in R\setminus\{0\}:\nu(x)\text{ is even}\}$.

\emph{Automaticity.} Fix $\gamma\in\{\gamma_1,\gamma_2\}$ and a window set $D$ for $\gamma$. For $c\in D^{(n)}$, if $c\notin I^n$ then $\nu(c)<n\le\nu(\gamma^ny)$ for all $y\in R$, and the corresponding kernel is constant. Otherwise $c=\gamma^nw$ with $w\in R$ belonging to a fixed finite set by the size bound (Lemma~\ref{lem:fixedlength}). Since $\nu(c+\gamma^ny)=n+\nu(w+y)$, the kernel is determined by $w$ and the parity of $n$ (with value zero at $y=-w$). The power invariance of automaticity (Lemma~\ref{lem:powers}) gives the conclusion.

\emph{No period at infinity.} Given $v\in R\setminus\{0\}$, put $j_0=\nu(v)$ and choose an integer $t>j_0$ of opposite parity. For $k\in\N$, the elements $x_k=\gamma_1^t(1+\gamma_1k)$ satisfy $\nu(x_k)=t$ and $\nu(x_k+v)=j_0$. Thus $\mathbf 1_Y(x_k)\ne\mathbf 1_Y(x_k+v)$ for infinitely many distinct $x_k$.
\end{proof}
For example, in $\Z[\sqrt2]$ the bases $3$ and $3(1+\sqrt2)$ satisfy (H) and generate the same prime ideal $3R$. The set of non-zero elements with even $3R$-valuation is automatic in both bases but is not eventually periodic.

\begin{proposition}[Common subfield]\label{prop:subfield}
Let $\alpha,\beta\in R$ be expansive. If a proper subfield $F\subsetneq K$ contains $\alpha^r$ and $\beta^s$ for some $r,s\ge1$, then $X=R\cap F$ is both $\alpha$- and $\beta$-automatic but is not eventually periodic.
\end{proposition}
\begin{proof}
Put $b=\alpha^r$ and let $\pi:K_{\R}\to K_{\R}/(F\otimes_{\Q}\R)$ be the quotient map. The image $\pi(R)$ is a lattice, since $F$ is a rational subspace of $K$, and multiplication by $b$ induces an invertible linear map on the quotient. For a window set $D$ of $b$ and $c\in D^{(n)}$, we have, for $y\in R$,
\[
 c+b^ny\in F\quad\Longleftrightarrow\quad \pi(y)=-\pi(b^{-n}c).
\]
By the size bound (Lemma~\ref{lem:fixedlength}), $\norm{b^{-n}c}\le C_D$. If the kernel is not identically zero, $\pi(b^{-n}c)$ belongs to the discrete group $\pi(R)$ and lies in a fixed bounded set. There are only finitely many such values, each determining the characteristic function of a coset of $R\cap F$ in $R$. Hence $X$ is $b$-automatic and thus $\alpha$-automatic by Lemma~\ref{lem:powers}. The argument for $\beta$ is identical.

For any $v\in R\setminus F$, the sets $X$ and $v+X$ are disjoint, and $X$ is infinite. Thus $\Perf(\mathbf 1_X)=R\cap F$ has rank $[F:\Q]<d$, which rules out eventual periodicity.
\end{proof}
For example, in $\Z[i]$ the bases $3$ and $5$ do not satisfy (H), and $\Z$ is automatic in both bases without being eventually periodic.

\begin{corollary}[Characterisation]\label{cor:characterisation}
For a number field $K$ of degree $d\ge2$ and expansive bases $\alpha,\beta\in R=\cO_K$, conditions (H) and (Ind) hold if and only if every subset of $R$ automatic in both bases is eventually periodic.
\end{corollary}

\begin{remark}[Weil bases]\label{rem:weil}
Let $K$ be a CM field and suppose that $\alpha\bar\alpha=q_\alpha$ and $\beta\bar\beta=q_\beta$, where $q_\alpha,q_\beta$ are integers greater than $1$, as in \cite{CMW} (equivalently, $|\sigma\alpha|$ and $|\sigma\beta|$ do not depend on the archimedean place $\sigma$). Multiplicative independence of $\alpha$ and $\beta$ implies (Ind). Indeed, if $\alpha^rR=\beta^sR$, then $\alpha^r=u\beta^s$ for a unit $u\in R^\times$. Multiplying by the complex conjugates gives $q_\alpha^r=u\bar u\,q_\beta^s$, and taking absolute norms gives $q_\alpha^r=q_\beta^s$. Hence $u\bar u=1$; as complex conjugation commutes with every embedding of the CM field $K$, we obtain $|\sigma u|=1$ at every archimedean place, and $u$ is a root of unity by Kronecker's theorem. For some integer $m\ge1$ we therefore have $\alpha^{rm}=\beta^{sm}$.
\end{remark}

\subsection{Determining which case applies}\label{ssec:examples}
For a pair of expansive bases in degree $d\ge2$, the following tests determine whether the conclusion of the main theorem holds and which part of the proof applies.
\begin{enumerate}[label=\textup{(\arabic*)},nosep]
\item (H) holds if and only if $K_\alpha K_\beta=K$ (Lemma~\ref{lem:fieldpowers}). For example, $K_\alpha=\Q(\sqrt2)$ and $K_\beta=\Q(\sqrt3)$ for the pair $(22(1+\sqrt2),2(4+\sqrt3))$ in $\Q(\sqrt2,\sqrt3)$, which therefore satisfies (H).
\item (Ind) fails if and only if the valuation vectors $(v_Q(\alpha))_Q$ and $(v_Q(\beta))_Q$ are proportional. If either condition fails, there is a doubly automatic set that is not eventually periodic (\S\ref{ssec:necessity}). We therefore assume both conditions in the remaining two steps.
\item If a prime ideal contains one base but not the other, Theorem~\ref{thm:different} applies, as for $(7,3+\sqrt2)$ in $\Z[\sqrt2]$.
\item Otherwise, compute the ray elements $a_0,b_0$. None of their archimedean absolute values is $1$ (Lemma~\ref{lem:modone}), and Proposition~\ref{prop:exray} or Proposition~\ref{prop:contracting} applies, according as one of them is expansive or not.
\end{enumerate}

The ideals alone do not determine whether the ray elements are expansive. In $K=\Q(\sqrt2,i\sqrt3)$, put $z=2+\sqrt2+i\sqrt3$ and $\varepsilon=1+\sqrt2$. The pair $(2z^2\bar z,2z\bar z^2)$ has two expansive ray elements, whereas the pair $(2z^2\bar z,2z\bar z^2\varepsilon^{-3})$, whose bases generate the same respective ideals, has contracting ray elements. Indeed, the ray elements are $2\bar z^3$ and $2z^3$ for the first pair, and $2\bar z^3\varepsilon^{-6}$ and $2z^3\varepsilon^3$ for the second; each of the latter has absolute value less than $1$ (about $0.57$ and $0.87$, respectively) at one of the two complex places.

\begin{remark}[Why two methods]\label{rem:twomethods}
Assume that (H) and (Ind) hold. The cone of \S\ref{ssec:cone} makes sense for any pair of expansive bases: for $r,s\in\Z$, $\alpha^r\beta^s\in R$ if and only if $r\,v_Q(\alpha)+s\,v_Q(\beta)\ge0$ at every prime ideal $Q$ containing $\alpha$ or $\beta$. If each base has an excess prime, the cone is $\N^2$ and its extremal elements are $\alpha$ and $\beta$. If only $\alpha$ has excess primes, they are $\beta$ and $\alpha^{m_1}/\beta^{n_1}$, where $m_1/n_1$, in lowest terms, is the largest value of $v_Q(\beta)/v_Q(\alpha)$ at the common primes. With equal supports, they are the ray elements $a_0,b_0$.

The argument of Section~\ref{sec:equal} compares the centres of two exact symmetries at an archimedean place where one multiplier contracts and the other expands (Lemma~\ref{lem:compcentres}); it therefore needs an extremal element which is not expansive. When each base has an excess prime, every monomial $\alpha^r\beta^s\ne1$ in $R$ is expansive, so no such comparison is available; nor does global transfer help, since there is no proper quotient to transfer to. Section~\ref{sec:different} works instead at an excess prime $Q$, where $\alpha$ is topologically nilpotent and $\beta$ is a unit. The powers of a unit stay in a compact group, so a $\beta$-symmetry yields rotations (Lemma~\ref{lem:rotations}) rather than a contradiction with a bound. Equality of centres (Proposition~\ref{prop:centres}) then comes from an open-mapping argument against the emptiness of the interior (Proposition~\ref{prop:emptyint}), which holds in the rigid class and therefore requires the reduction of Section~\ref{sec:orbits}.

Thus the contracting-ray argument of Section~\ref{sec:equal} is used precisely when the supports are equal and neither ray element is expansive. In all other cases the proof goes through an excess prime, directly or after the transfer of Proposition~\ref{prop:exray}, which replaces $\beta$ by the expansive ray element $a_0$ and creates an excess prime for $\alpha$. In the intermediate case where only $\alpha$ has excess primes and $\alpha^{m_1}/\beta^{n_1}$ contracts at some archimedean place, a contracting multiplier exists, but the argument of Section~\ref{sec:equal} does not apply as it stands: step~2 of the proof of Lemma~\ref{lem:twosym} needs a power of the base to be divisible by the contracting factor, whereas, with $b=\alpha^{m_1}/\beta^{n_1}$, no power of $\beta$ is divisible by $b$, because of the primes dividing $\alpha$ but not $\beta$. This case is covered by Section~\ref{sec:different}.
\end{remark}

\section{Orders, lattices and matrix bases}\label{sec:orders}

Throughout this section, $K$ is a number field of degree $d=[K:\Q]$ and $R=\cO_K$. The main theorem is stated for $R$. We extend it to any lattice in $K$ stable under both bases, in particular to any order. A \emph{lattice} is a free additive subgroup $M\subseteq K$ of rank $d$ spanning $K$ over $\Q$. Its \emph{multiplier ring} is $\mathfrak o(M)=\{x\in K:xM\subseteq M\}$. Fix an integer $m\ge1$ such that $mR\subseteq M\subseteq m^{-1}R$.

\begin{lemma}[Multiplier ring]\label{lem:multiplier}
The ring $\mathfrak o(M)$ is an order in $K$, and $m^2R\subseteq\mathfrak o(M)\subseteq R$. In particular, every order in $K$ is contained in $R$.
\end{lemma}
\begin{proof}
An element $x\in\mathfrak o(M)$ acts on $M\simeq\Z^d$ by an integer matrix. The characteristic polynomial of this matrix is monic with integer coefficients and annihilates $x$, so $x\in R$. Conversely, if $x\in m^2R$, then $xM\subseteq m^2R\cdot m^{-1}R=mR\subseteq M$. Thus $\mathfrak o(M)$ is a subring of $R$ containing $1$ and a full-rank subgroup, hence an order.
\end{proof}

\subsection{Numeration on a lattice}\label{ssec:ordernumeration}
Let $\gamma\in\mathfrak o(M)$ be expansive. The definitions of Section~\ref{sec:numeration} apply with $M$ in place of $R$: digit sets are contained in $M$, kernels are functions $y\mapsto f(c+\gamma^ny)$ on $M$, and eventual periodicity is defined using subgroups of $M$ of rank $d$.

The results of Section~\ref{sec:numeration} and Proposition~\ref{prop:converse} hold on $M$ with the same proofs, with the following minor changes. In Lemma~\ref{lem:short}, use the covering radius $c_M$ of the image of $M$ in $K_{\R}$ in place of $c_R$, and round the quotient in $M$; the digit then belongs to $M$ since $\gamma M\subseteq M$. Multiplication by $\gamma$ is injective and
\[
 [M:\gamma^nM]=|N_{K/\Q}(\gamma)|^n.
\]
In Proposition~\ref{prop:converse}, use $q:=[M:\Lambda]$. The affine maps $y\mapsto t+\eta y$ of Proposition~\ref{prop:affine} and Lemma~\ref{lem:insert} map $M$ into itself when $t\in M$ and $\eta\in\mathfrak o(M)\setminus\{0\}$. Example~\ref{ex:parity} and the valuations $v_P$ concern $R$ alone and are not needed here.

\subsection{From a lattice to the ring of integers}\label{ssec:ordertransfer}
For a function $f:M\to\Omega$ and $a\in M$, define $g_a:R\to\Omega$ by $g_a(y)=f(a+my)$; this makes sense since $mR\subseteq M$.

\begin{lemma}[Transfer to $\cO_K$]\label{lem:ordertransfer}
If $\gamma\in\mathfrak o(M)$ is expansive and $f:M\to\Omega$ is $\gamma$-automatic, then $g_a$ is $\gamma$-automatic on $R$ for every $a\in M$.
\end{lemma}
\begin{proof}
Let $D$ and $D'$ be window sets for $\gamma$ on $M$ and on $R$, respectively. Given $n\in\N$ and $c\in D'^{(n)}$, write $a+mc=c^*+\gamma^nw$ with $c^*\in D^{(n)}$ and $w\in M$. Then
\[
 g_a(c+\gamma^ny)=\Phi(w+my),\qquad \Phi(z)=f(c^*+\gamma^nz).
\]
The function $\Phi$ belongs to the finite kernel family of $f$. At each archimedean place,
\[
 |\sigma w|\le |\sigma a|+mC_{D'}+C_D.
\]
Since $mw\in R$, the possible values of $w$ form a finite set. Thus $g_a$ has only finitely many kernels.
\end{proof}

\begin{lemma}[(Ind) on a lattice]\label{lem:latticeind}
Let $\alpha,\beta\in\mathfrak o(M)\setminus\{0\}$. The following three conditions, each required for all $r,s\ge1$, are equivalent:
\[
 \alpha^rM\ne\beta^sM,\qquad
 \alpha^r\mathfrak o(M)\ne\beta^s\mathfrak o(M),\qquad
 \alpha^rR\ne\beta^sR.
\]
\end{lemma}
\begin{proof}
Put $u=\beta^s/\alpha^r$. The equality $\alpha^rM=\beta^sM$ holds if and only if $uM=M$, that is, if and only if $u\in\mathfrak o(M)^\times$; this is also equivalent to $\alpha^r\mathfrak o(M)=\beta^s\mathfrak o(M)$.

For an order $\mathfrak o$ and an integer $m'\ge1$ with $m'R\subseteq\mathfrak o$, the kernel of the reduction map $R^\times\to(R/m'R)^\times$ lies in $\mathfrak o^\times$: if $v\equiv1\pmod{m'R}$, then $v\in\mathfrak o$ and $v^{-1}=1-v^{-1}(v-1)\in\mathfrak o$. Hence $R^\times/\mathfrak o^\times$ is finite. If $\alpha^rR=\beta^sR$, then $u\in R^\times$ and $u^k\in\mathfrak o(M)^\times$ for some $k\ge1$, so $\alpha^{rk}\mathfrak o(M)=\beta^{sk}\mathfrak o(M)$. The reverse implication follows by extending ideals to $R$.
\end{proof}
Thus (H) and (Ind) for a pair in $\mathfrak o(M)$ are the same as in $R$; the choice of lattice merely restricts the available bases. For instance, $3$ and $3i$ generate the same ideal of $\Z[i]$ but not of $\Z[3i]$, whereas their squares generate the same ideal in both rings.

\subsection{The theorem for lattices}\label{ssec:ordertheorem}
\begin{theorem}[Lattices and orders]\label{thm:orders}
Let $K$ be a number field of degree $d\ge2$, let $M$ be a lattice in $K$ and let $\alpha,\beta\in\mathfrak o(M)$ be expansive and satisfy \textup{(H)} and \textup{(Ind)}. Then $f:M\to\Omega$, with $\Omega$ finite, is both $\alpha$- and $\beta$-automatic if and only if it is eventually periodic. In particular, the theorem holds on every order in $K$.
\end{theorem}
\begin{proof}
The converse is Proposition~\ref{prop:converse} on $M$. Suppose that $f$ is doubly automatic and let $T\subseteq M$ be a complete set of representatives of $M/mR$. By Lemma~\ref{lem:ordertransfer}, each $g_a$ ($a\in T$) is doubly automatic on $R$. By Theorem~\ref{thm:main}, it agrees outside a finite set $X_a\subseteq R$ with a function $h_a$ periodic under a full-rank subgroup $L_a\subseteq R$. Put $L=\bigcap_{a\in T}L_a$ and define $h:M\to\Omega$ by $h(a+my)=h_a(y)$ for $a\in T$ and $y\in R$. Then $h$ is periodic under the full-rank subgroup $mL$ of $M$, and $f=h$ outside the finite set $\bigcup_{a\in T}(a+mX_a)$.
\end{proof}

\begin{corollary}[Characterisation for lattices]\label{cor:charorders}
For $K$ of degree $d\ge2$, a lattice $M$ in $K$ and expansive bases $\alpha,\beta\in\mathfrak o(M)$, conditions (H) and (Ind) hold if and only if every subset of $M$ automatic in both bases is eventually periodic.
\end{corollary}
\begin{proof}
Sufficiency is Theorem~\ref{thm:orders}. If (Ind) fails, on $M$ or equivalently in $R$ (Lemma~\ref{lem:latticeind}), we have $\alpha^rM=\beta^sM$ for some $r,s\ge1$. The proof of Proposition~\ref{prop:idealdep} applies on $M$, with $\gamma_1:=\alpha^r$, $\gamma_2:=\beta^s$ and the subgroups $\gamma_1^jM=\gamma_2^jM$ in place of the ideals. Three points deserve mention. First, $u=\gamma_2/\gamma_1$ is a unit of $\mathfrak o(M)$ by Lemma~\ref{lem:latticeind}, so multiplication by $u^n$ preserves each $\gamma_1^jM$. Second, $\nu(x)$ is finite for $x\in M\setminus\{0\}$: if $x\in\gamma_1^nM$, then $m\gamma_1^{-n}x\in R\setminus\{0\}$, whence
\[
 |N_{K/\Q}(\gamma_1)|^n\le m^d|N_{K/\Q}(x)|.
\]
Finally, $M$ need not contain $1$. In the proof that no period at infinity exists, use instead $x_k=\gamma_1^t(m_0+\gamma_1k)$ with $k\in M$ and a fixed $m_0\in M\setminus\gamma_1M$. Such an $m_0$ exists since $[M:\gamma_1M]=|N_{K/\Q}(\gamma_1)|\ge2$. The elements $m_0+\gamma_1k$ lie in $M$, since $\gamma_1\in\mathfrak o(M)$, and have the same non-zero residue class modulo $\gamma_1M$, so the distinct elements $x_k$ all satisfy $\nu(x_k)=t$.

If (H) fails, a proper subfield $F$ contains $\alpha^r$ and $\beta^s$ for suitable $r,s\ge1$. The proof of Proposition~\ref{prop:subfield} applies with $X=M\cap F$, which is a full-rank lattice in $F$ (it contains $m\cO_F$) and is saturated in $M$. This gives the required counterexample.
\end{proof}

\subsection{Matrix bases}\label{ssec:matrices}
A matrix $A\in M_k(\Z)$ is \emph{expansive} if all its eigenvalues have absolute value greater than $1$. For a finite set $D\subseteq\Z^k$ containing $0$, put $D^{(n)}=\{\sum_{j<n}A^jw_j:w_j\in D\}$. We call $D$ a window set if $D^{(n)}$ meets every residue class of $\Z^k/A^n\Z^k$ for every $n\in\N$. A function $f:\Z^k\to\Omega$ is $A$-automatic if the kernels $y\mapsto f(c+A^ny)$ ($n\in\N$, $c\in D^{(n)}$) form a finite family for a window set $D$.

\begin{corollary}[Commuting matrix bases]\label{cor:matrices}
Let $d\ge2$ and let $A,B\in M_d(\Z)$ be commuting expansive matrices such that $\Q[A,B]$ is a field of degree $d$ (for instance, if the characteristic polynomial of $A$ is irreducible). Assume that \textup{(i)} $\dim_{\Q}\Q[A^r,B^s]=d$ for all $r,s\ge1$ and \textup{(ii)} $A^r\Z^d\ne B^s\Z^d$ for all $r,s\ge1$. Then $f:\Z^d\to\Omega$, with $\Omega$ finite, is both $A$- and $B$-automatic if and only if it is eventually periodic. If either (i) or (ii) fails, there is a subset of $\Z^d$ automatic in both bases but not eventually periodic.
\end{corollary}
\begin{proof}
Let $K=\Q[A,B]$. The space $\Q^d$ is a one-dimensional vector space over $K$, so for any $z_0\in\Z^d\setminus\{0\}$ the map $x\mapsto xz_0$ is a $K$-linear isomorphism $\varphi:K\to\Q^d$. The inverse image $M=\varphi^{-1}(\Z^d)$ is a lattice in $K$, and the elements $\alpha=A$ and $\beta=B$ of $K$ belong to $\mathfrak o(M)$. Under $\varphi$, window sets, kernels and subgroups of rank $d$ correspond. The matrix $A$ is expansive if and only if $\alpha$ is expansive, since its eigenvalues are the conjugates of $\alpha$; likewise for $B$. Condition (i) is (H) and (ii) is (Ind) on $M$ (Lemma~\ref{lem:latticeind}). We conclude by Theorem~\ref{thm:orders} and Corollary~\ref{cor:charorders}.

If the characteristic polynomial of $A$ is irreducible, $\Q^d$ is a one-dimensional vector space over $\Q[A]$. Every matrix commuting with $A$ is then $\Q[A]$-linear and belongs to $\Q[A]$, so $\Q[A,B]=\Q[A]$.
\end{proof}
Every lattice of $K$ arises in this way, by taking matrices of multiplication by suitable elements in a $\Z$-basis of the lattice. By the theorem of Latimer and MacDuffee \cite{LatimerMacDuffee}, the $\mathrm{GL}_d(\Z)$-conjugacy classes of integer matrices with a given monic irreducible characteristic polynomial correspond to classes of lattices stable under a root of that polynomial, up to multiplication by non-zero field elements. Such lattices need not be orders, which explains the formulation of our theorem for lattices.

\section{Further results and questions}\label{sec:questions}

We keep the notation of the previous sections: $K$ is a number field of degree $d$, $R=\cO_K$ and $\Omega$ is a finite alphabet. The main theorem describes double automaticity on $R$ for two bases satisfying (H) and (Ind). This section proves two further results, Proposition~\ref{prop:farF} on periods at infinity without (H) and Proposition~\ref{prop:factormap} on prime ideals under a morphism, and collects the questions which seem most natural to us: modules over $R$ (\S\ref{ssec:qmodules}), several bases (\S\ref{ssec:farF}), more general expansive endomorphisms (\S\ref{ssec:groups}), the arithmetic determined by the dynamics (\S\ref{ssec:factormaps}) and decidability (\S\ref{ssec:qdecidability}). Nothing in this section is used elsewhere in the paper.

\subsection{Modules over \texorpdfstring{$\cO_K$}{the ring of integers}}\label{ssec:qmodules}
Let the bases act diagonally on $R^e$, with $e\ge1$. Since $R^e$ is a lattice in $(K\otimes_{\Q}\R)^e$ on which these actions are expansive, the automaticity of \S\ref{ssec:matrices} applies to the block-diagonal matrices of multiplication by the bases, and the constructions of Section~\ref{sec:numeration} extend to this setting. For $e>1$, eventual periodicity is no longer the right conclusion: $R\times\{0\}^{e-1}$ is doubly automatic. More generally, let $V\subseteq K^e$ be a $K$-vector subspace, $\Lambda\subseteq R^e$ a subgroup of finite index and $a\in R^e$. The set $R^e\cap V$ is automatic in every expansive scalar base, by the argument of Proposition~\ref{prop:subfield} applied in $(K\otimes_{\Q}\R)^e$ to the real subspace $V\otimes_{\Q}\R$, which is rational with respect to $R^e$. The same therefore holds for $a+(\Lambda\cap V)$, by stability under intersections with periodic sets and under translations. Condition (H) gives an algebraic reason for the choice of directions: a $\Q$-subspace of $K^e$ stable under multiplication by $\alpha^r$ and $\beta^s$, for some $r,s\ge1$, is stable under $\Q[\alpha^r,\beta^s]=K$ (Lemma~\ref{lem:fieldpowers}), hence $K$-linear.

\begin{question}\label{q:modules}
Let $K$ be a number field of degree $d\ge2$, and let $\alpha,\beta\in R=\cO_K$ be expansive and satisfy (H) and (Ind), acting diagonally on $R^e$ for an integer $e\ge1$. Is every subset of $R^e$ automatic in both bases a finite Boolean combination of sets $a+(\Lambda\cap V)$, where $a\in R^e$, $\Lambda\subseteq R^e$ is an additive subgroup of rank $de$, and $V\subseteq K^e$ is a $K$-vector subspace?
\end{question}
For $e=1$, the only $K$-subspaces of $K$ are $0$ and $K$, the finite Boolean combinations in the question are exactly the eventually periodic sets, and Theorem~\ref{thm:main} gives a positive answer.

\emph{Half-spaces.} The Presburger description in the Cobham--Semenov theorem allows half-spaces; here, for $d\ge2$ and under (H) and (Ind), they are not doubly automatic. The restriction of an $\alpha$-automatic function on $R^e$ to an affine line $t+Rw$, with $t,w\in R^e$ and $w\ne0$, is $\alpha$-automatic as a function of the parameter in $R$, by the carry argument of Proposition~\ref{prop:affine}. A non-trivial half-space of $(K\otimes_{\Q}\R)^e$ meets a suitable such line in a non-trivial half-space of $R$, which is not eventually periodic: every subgroup of full rank contains a period transverse to the bounding hyperplane, so each of its cosets has points of arbitrarily large norm on both sides. If the half-space were doubly automatic, so would be this trace, contrary to Theorem~\ref{thm:main}. The difference between rational bases and the scalar actions considered here is thus already visible on lines.

\emph{Periods.} The principle of Theorem~\ref{thm:onefar} suggests that periods should also play a role on $R^e$ for $e>1$. The difficulty is that lower-dimensional sets can no longer be treated as finite exceptions.

\emph{The case where (H) fails.} Let $F=K_\alpha K_\beta=\Q(\alpha^r,\beta^s)$ for suitable $r,s\ge1$ (Lemma~\ref{lem:fieldpowers}), and assume $F\ne K$. If $w_1,\ldots,w_e\in R$ is an $F$-basis of $K$, the subgroup $\cO_Fw_1+\cdots+\cO_Fw_e\simeq\cO_F^{\,e}$ has finite index in $R$, and $\alpha^r,\beta^s\in\cO_F$ act on it as scalars. Passing to these powers and to the finitely many cosets of this subgroup, we are led to automaticity on $\cO_F^{\,e}$ with diagonal scalar actions. The pair $(\alpha^r,\beta^s)$ satisfies (Ind) over $\cO_F$ whenever $(\alpha,\beta)$ does, since an equality of principal ideals of $\cO_F$ persists after extension to $R$. For $F=\Q$ and multiplicatively independent bases, this is the setting of the Cobham--Semenov theorem \cite{Semenov} (see \cite{Muchnik} for $\Z^n$ and decidability, and \cite{DurandSemenov} for a dynamical proof). Together with that theorem, a classification on $\cO_F^{\,e}$ for scalar bases satisfying (Ind) and generating $F$ after every passage to powers would address, through this reduction, the excluded case where (H) fails but (Ind) holds in $K$. The same question may be asked when the coordinates carry different powers of the bases (\S\ref{ssec:groups}), and for more than two bases (\S\ref{ssec:farF}).

\subsection{Periods at infinity without (H)}\label{ssec:farF}
Under (H), Theorem~\ref{thm:onefar} shows that a single non-zero period at infinity forces a full-rank group of periods at infinity. Without (H), the same argument still imposes an algebraic structure on the space of periods: the relevant field is now the compositum $F$ introduced in \S\ref{ssec:qmodules}. The following proposition assumes neither (H) nor (Ind).

\begin{proposition}[Periods at infinity without \textup{(H)}]\label{prop:farF}
Let $F=K_\alpha K_\beta$, the smallest of the fields $\Q(\alpha^r,\beta^s)$ for $r,s\ge1$ (Lemma~\ref{lem:fieldpowers}), and let $f:R\to\Omega$ be both $\alpha$- and $\beta$-automatic. Then $\Pi(f)=\operatorname{span}_{\Q}\Perf(f)$ is an $F$-vector subspace of $K$. In particular, the rank of $\Perf(f)$ is zero or a multiple of $[F:\Q]$.
\end{proposition}
\begin{proof}
We return to the proof of Theorem~\ref{thm:onefar}. Condition (H) is used only in the last step of Lemma~\ref{lem:good}. Without it, the field generated in that lemma still satisfies
\[
 \Q[\eta_1,\eta_2]\supseteq\Q(\alpha^{|\Delta|},\beta^{|\Delta|})\supseteq F,
\]
by Lemma~\ref{lem:fieldpowers}.

Take $v\in\Perf(f)\setminus\{0\}$ and choose $(k,l)=tw_0$ at distance greater than $|\cG|$ from $\cR$, as in \S\ref{ssec:proofonefar}. Put $\Gamma=\alpha^k\beta^l$. For each residue class $c$ in the complete set of representatives given by Lemma~\ref{lem:graph}\,(c), take a path $\pi_c$ of exponent $(k,l)$ ending at $g_c:y\mapsto f(c+\Gamma y)$. Each such path meets a two-dimensional component or two one-dimensional components of different directions: otherwise the path removal argument of Lemma~\ref{lem:bad} would place its exponent within distance $|\cG|$ of $\cR$. We may therefore follow the proof of Lemma~\ref{lem:good} with $\pi=\pi_c$ and $\Gamma_0=\Gamma$, up to the inclusion $\Pi(g_c)\supseteq\Gamma^{-1}v\,\Q[\eta_1,\eta_2]$. Crucially, the scale $\Gamma$ is the same for all $c$. The resulting inclusion is
\[
 \Gamma^{-1}vF\subseteq\Pi(g_c).
\]
Since $\Perf(g_c)$ has finite index in $R\cap\Pi(g_c)$, it follows that $\Perf(g_c)\cap\Gamma^{-1}vF$ has finite index in the lattice $R\cap\Gamma^{-1}vF$. Intersecting these subgroups over the finitely many $c$, we obtain a full-rank subgroup $L'$ of that lattice. By the last step of the proof of Theorem~\ref{thm:onefar}, $\Gamma L'\subseteq\Perf(f)$. Hence $vF\subseteq\Pi(f)$. Since this holds for every $v\in\Perf(f)\setminus\{0\}$, the $\Q$-span $\Pi(f)$ is stable under multiplication by $F$.
\end{proof}
Both rank zero and the smallest positive rank allowed by the proposition occur. In $K=\Q(\sqrt2,\sqrt3)$, take $\alpha=4+\sqrt2$ and $\beta=5+\sqrt2$. Then $F=\Q(\sqrt2)$ and $\Perf(\mathbf 1_{\cO_F})=\cO_F$, of rank $2=[F:\Q]$. In $\Z[i]$ with bases $2$ and $3$ (so $F=\Q$), the quadrant $\{a+bi:a,b\in\Z_{\ge0}\}$ has a single kernel in each base for the window sets $\{a+bi:0\le a,b<k\}$, $k=2,3$, and has no period at infinity. The proposition thus restricts the geometry of the periods, not the functions themselves: as the quadrant shows, there may be no non-zero period at all.

\emph{Several bases.} Applied to pairs, Proposition~\ref{prop:farF} shows that, for a function automatic in expansive bases $\alpha_1,\ldots,\alpha_t\in R$ ($t\ge2$), the $\Q$-span of the periods at infinity is stable under the compositum $K_{\alpha_1}\cdots K_{\alpha_t}$. When this compositum is $K$, a single non-zero period at infinity forces full rank, and hence eventual periodicity for $d\ge2$. This can happen even though no pair satisfies (H), as in $\Q(\sqrt2,\sqrt3,\sqrt5)$ with bases $4+\sqrt2$, $4+\sqrt3$ and $4+\sqrt5$. Since propagation is then established, the issue in the following question is to obtain one non-zero period at infinity; the condition on principal ideals excludes common filtrations.

\begin{question}\label{q:severalbases}
Let $d\ge2$ and let $\alpha_1,\ldots,\alpha_t\in R$ ($t\ge2$) be expansive, with $\Q(\alpha_1^{n_1},\ldots,\alpha_t^{n_t})=K$ for all $n_i\ge1$, and $\alpha_i^rR\ne\alpha_j^sR$ whenever $i\ne j$ and $r,s\ge1$. Is every function $f:R\to\Omega$ automatic in all these bases eventually periodic, even if no pair satisfies (H)?
\end{question}

\subsection{Groups with expansive endomorphisms}\label{ssec:groups}
The setting of Sections~\ref{sec:numeration} and~\ref{sec:orders} is that of a group isomorphic to $\Z^k$ equipped with an expansive endomorphism $A$ (an integer matrix whose eigenvalues have absolute value greater than $1$). Numeration, kernels, the kernel criterion and independence of the choice of digits may be formulated as in \S\ref{ssec:matrices}, using a norm adapted to $A$ for which $\norm{A^{-1}}<1$. Automatic sets with respect to an endomorphism of a finitely generated abelian group have been studied by Bell and Moosa \cite{BellMoosa} and by Hawthorne \cite{Hawthorne}, who asks for an analogue of Cobham's theorem for two endomorphisms \cite[Ch.~7, p.~119]{HawthorneThesis}; see also \cite{vonHaeseler} for sequences indexed by a group with an expansive endomorphism.

Corollary~\ref{cor:matrices} treats two commuting endomorphisms generating a field of degree $k$. Two natural directions extend the setting of this corollary: reducible actions, in which invariant subspaces may give doubly automatic sets of lower rank, as in \S\ref{ssec:qmodules}; and irreducible actions of endomorphisms which do not commute. For the latter, our hypotheses have natural analogues. In the field case, (H) says that, for any positive powers of the two bases, the $\Q$-algebra they generate acts irreducibly on the $\Q$-vector space $K$; this makes sense without commutativity. The analogue of (Ind) compares the images of the lattice, not the matrices themselves, as is already necessary for elements of a number field, which may generate the same ideal without being equal.

\begin{question}\label{q:noncommuting}
Let $k\ge2$ and let $A,B\in M_k(\Z)$ be expansive matrices such that, for all $r,s\ge1$, the algebra $\Q[A^r,B^s]$ acts irreducibly on $\Q^k$ and $A^r\Z^k\ne B^s\Z^k$. Is every subset of $\Z^k$ which is both $A$- and $B$-automatic eventually periodic?
\end{question}
If some positive power $A^n=cI$ is scalar, the answer is yes. The hypothesis with $r=n$ says that the commutative algebra $\Q[B^s]$ acts irreducibly on $\Q^k$ for every $s\ge1$; by Schur's lemma it is a field of degree $k$. The pair $(cI,B)$ then falls within Corollary~\ref{cor:matrices}, and an $A$-automatic set is $A^n$-automatic (the argument of Lemma~\ref{lem:powers} applies to matrices). The same holds if a positive power of $B$ is scalar. In general, commutativity is used in the construction of mixed kernels and of their graph (Proposition~\ref{prop:mixed}, Lemmas~\ref{lem:graph}\,(b) and~\ref{lem:good}), in the comparison of symmetries (Proposition~\ref{prop:centres} and Lemma~\ref{lem:compcentres}), and in the integer cone of \S\ref{ssec:cone} and the unit equations of Theorem~\ref{thm:sparse}.

\emph{Reducible actions.} The scalar actions of \S\ref{ssec:qmodules} are, for $e\ge2$, a particular instance of reducible commuting actions. A simple extension is to consider $(\alpha^k,\alpha^l)$ and $(\beta^k,\beta^l)$ on $R^m\times R^n$, with positive integers $k\ne l$. It is natural to ask how the geometry of the doubly automatic sets reflects the structure of the two endomorphisms.

\emph{Beyond abelian groups.} Automaticity with respect to expansive endomorphisms also leads to questions on non-abelian groups. Heisenberg groups over $R$ provide a natural example, in which the natural dilations give different weights to the horizontal and central coordinates. In this setting, the question is not merely whether a Cobham-type theorem holds, but what form its conclusion should take.

The common-filtration obstruction of \S\ref{ssec:necessity} also raises the question of appropriate arithmetic independence conditions for more general endomorphisms.

\subsection{Morphisms between substitution subshifts}\label{ssec:factormaps}
Cabezas and Leroy \cite[\S7.2]{CabezasLeroy} recall a dynamical form of Cobham's theorem: by a result of Fagnot \cite{Fagnot}, if the subshift of a primitive aperiodic constant-length substitution is a factor of that of another such substitution, their lengths have a common positive power. They observe that, by their Theorem~4.1 and a result of Durand \cite{DurandSemenov}, the same holds in higher dimension when the expansion matrices are multiples of the identity, and they ask for an analogue for general constant-shape substitutions (\cite[Question~7.2]{CabezasLeroy}). We obtain a necessary arithmetic condition for such a factor map when the expansions are given by multiplication in a number field.

In the case of different supports, the odometer factor of an aperiodic substitution subshift is incompatible with automaticity in a base which is invertible at one of the primes of its expansion (\S\S\ref{ssec:odometer} and~\ref{ssec:rankzero}). The argument applies to a morphism: a local rule sends an automatic configuration of the source to a configuration which is still automatic in the source base, while the target subshift has the odometer factors of its own expansion. Neither (H) nor (Ind) is needed.

\begin{proposition}[Prime ideals under a morphism]\label{prop:factormap}
Let $\alpha,\beta\in R$ be expansive, and let $\theta,\theta'$ be primitive constant-shape substitutions with expansions $\alpha^a,\beta^b$ ($a,b\ge1$) and shapes $E_\theta,E_{\theta'}$ such that $(E_\theta,\alpha^a)$ and $(E_{\theta'},\beta^b)$ are adapted windows. Assume that $X_\theta$ and $X_{\theta'}$ are aperiodic. If there is a morphism $X_\theta\to X_{\theta'}$ (a continuous, translation-equivariant map), every prime ideal containing $\beta$ contains $\alpha$. In particular, if $X_\theta$ and $X_{\theta'}$ are conjugate, $\alpha$ and $\beta$ have the same prime support.
\end{proposition}
\begin{proof}
Let $\varphi:X_\theta\to X_{\theta'}$ be the morphism and choose $x\in X_\theta$ fixed by some positive power $\theta^M$ (Lemma~\ref{lem:germ}). The kernel argument of Lemma~\ref{lem:generator}, with $\pr=\mathrm{id}$, shows that $x$ is $\alpha^a$-automatic and hence $\alpha$-automatic by Lemma~\ref{lem:powers}. By the Curtis--Hedlund--Lyndon theorem \cite[\S2.2]{CabezasLeroy}, there are a finite set $P\subseteq R$ and a local rule $\Phi$ such that $y:=\varphi(x)$ satisfies $y(q)=\Phi((x(q+w))_{w\in P})$ for all $q\in R$. Thus $y$ is a pointwise function of finitely many translates of $x$ and is $\alpha$-automatic by Proposition~\ref{prop:affine} and Corollary~\ref{cor:closure}. Since $X_{\theta'}$ is minimal (Lemma~\ref{lem:minimal}), it equals $\cX_y$.

The odometer factor of $X_{\theta'}$ (Corollary~\ref{cor:odometer}, with $\beta$ in place of $\alpha$) gives a positive power $\nu$ of $\beta$ and factor maps $\psi_n:X_{\theta'}\to R/\nu^nR$. Write $\chi_n:R\to R/\nu^nR$ for the quotient map. Applying Proposition~\ref{prop:factorbound} to $y$ in base $\alpha$, with $Z=X_{\theta'}$, the group $R/\nu^nR$ and $\psi=\psi_n$, gives a bound $N^*$ depending only on $y$, $\alpha$ and a window set, and hence independent of $n$, such that for each $n\ge1$,
\[
 |\chi_n(\alpha^kR)|\le N^*
\]
for all $k$ beyond a threshold depending on $n$ (see the paragraph following Theorem~\ref{thm:finiteorbits}). If a prime ideal $Q$ contains $\beta$ but not $\alpha$, then $R/Q^n$ is a quotient of $R/\nu^nR$. As $\alpha$ is invertible modulo $Q^n$, the subgroup $\alpha^kR$ maps onto $R/Q^n$. Consequently,
\[
 |\chi_n(\alpha^kR)|\ge |R/Q^n|\ge2^n
\]
for all $k\in\N$. Choosing $n$ with $2^n>N^*$ and then $k$ beyond the corresponding threshold gives a contradiction.
\end{proof}
The proposition reads the primes of the expansion off the dynamics alone, without any combinatorial relation between the two substitutions; a morphism gives an inclusion of supports, a conjugacy their equality. It is natural to ask how much more of the arithmetic of the expansion, in particular of the filtration by its powers, is determined by the dynamics.

The conclusion concerns prime support, not exponents. For matrices whose algebra is not a field, the analogue $A^r\Z^d=B^s\Z^d$ of the equality $\alpha^rR=\beta^sR$ in Question~\ref{q:factormap} below may fail: if $\theta,\theta'$ are primitive aperiodic substitutions of lengths $2$ and $3$, then their product substitution and the product of $\theta^2$ and ${\theta'}^3$ define the same subshift $X_\theta\times X_{\theta'}$ on $\Z^2$, with respective expansion matrices $\mathrm{diag}(2,3)$ and $\mathrm{diag}(4,27)$. No positive powers of the two matrices have equal images on $\Z^2$. This reducible example (\S\ref{ssec:groups}) shows that such an equality cannot be expected for arbitrary expansion matrices. In the field setting there is a reason to ask for it: under (H) and (Ind), the main theorem excludes doubly automatic configurations whose orbit closure is minimal and aperiodic. But the image of the automatic point in the proof above is automatic in the source base only, and membership in the target subshift does not make it automatic in the target base. The main theorem therefore does not immediately answer the following question.

\begin{question}\label{q:factormap}
In the setting of Proposition~\ref{prop:factormap}, and under (H) for $(\alpha,\beta)$, does the existence of a morphism $X_\theta\to X_{\theta'}$ imply that $\alpha^rR=\beta^sR$ for some $r,s\ge1$? Does this hold at least when the two subshifts are conjugate?
\end{question}
A positive answer would give a form of the rigidity sought in \cite[Question~7.2]{CabezasLeroy} in the arithmetic setting of this paper: a factor map would force not only an inclusion of prime supports, but also proportionality of the prime-ideal valuations of the two bases.

\subsection{Decidability}\label{ssec:qdecidability}
The main theorem, which involves two bases, suggests an algorithmic question about a single automaton: given a function automatic in one base, decide whether it is eventually periodic. Under (H) and (Ind), by Corollary~\ref{cor:farperiods}, a function given by an automaton in base $\alpha$ is $\beta$-automatic if and only if $\Perf(f)$ has rank $d$; deciding this would decide automaticity in a second base, with no automaton in that base as input.

For integer bases, eventual periodicity of an automatic subset of $\N$ and Presburger definability of a recognisable subset of $\Z^m$ are decidable \cite{Honkala,Muchnik,Leroux,MarsaultSakarovitch}; see also \cite{BHMV}. Cabezas and Leroy \cite[Question~7.5]{CabezasLeroy} ask, for integer expansion matrices, for analogues of the equivalences between automaticity, substitution representations and first-order definability, and for the decidability of the corresponding theory; in our setting, the kernel substitution of \S\ref{ssec:kernelsubst} gives a substitution representation, through an adapted window and a positive power of the base.

For the narrower problem above, automatic structures \cite{KhoussainovNerode,BlumensathGraedel} seem well suited. Let $\gamma\in R$ be expansive with an admissible digit set. Equality of values and addition are recognised by carry automata (as in the proof of Theorem~\ref{thm:sparse}), so that $R$, equipped with addition, a regular well-order and the level sets of a $\gamma$-automatic function $f$, is an automatic structure; its first-order theory with the quantifier $\exists^\infty$ (``there exist infinitely many'') is decidable \cite[Cor.~5.4]{BlumensathGraedel}. The group of periods at infinity, whose definition does not refer to a second base, is expressed by
\[
 v\in\Perf(f)\quad\Longleftrightarrow\quad
 \neg\,\exists^\infty x\ \bigl(f(x+v)\ne f(x)\bigr).
\]
One may therefore expect that it is decidable whether $\Perf(f)$ has finite index, with computable index when it is. For $d\ge2$, this would also decide eventual periodicity, and one may further ask for an effective determination of the finite exceptional set. We do not pursue these algorithmic questions here.

\section*{Acknowledgements}
{\emergencystretch=2em We thank Jean-Paul Allouche for his encouragement and for suggesting the reference to von Haeseler's book.\par}

\section*{Statement on the use of artificial intelligence}
Part of this article grew out of extended exchanges with several artificial intelligence tools. The original arguments have since been substantially reworked and rewritten. The authors take full responsibility for the content.

\bigskip
\noindent\begin{minipage}[t]{0.48\textwidth}\small
Francesco L. Chiera\\ Sceaux, France\\ \texttt{fchiera@ac-versailles.fr}
\end{minipage}\hfill
\begin{minipage}[t]{0.48\textwidth}\small
Taoufik Safer\\ Paris, France\\ \texttt{Taoufik.Safer@ac-versailles.fr}
\end{minipage}
\end{document}